\documentclass[11pt]{article}
 \usepackage{amsfonts,amssymb,enumerate,epsf,amsmath,bbm}
 \usepackage[mathscr]{eucal}
 \usepackage{graphicx}
 \usepackage{xcolor}
 \usepackage{standalone}
\date{}
\usepackage{hyperref}
\usepackage{verbatim}
\usepackage{tikz}
\usepackage{pgfplots}
\usepackage{pgfplotstable}
\usepgfplotslibrary{colorbrewer}
\usepgfplotslibrary{groupplots}
\usetikzlibrary{patterns}
\usetikzlibrary{arrows}
\usepackage{authblk}
\usepackage{mathrsfs}

\definecolor{grigio}{RGB}{232,232,232}
\definecolor{grigio2}{RGB}{182,182,182}

 \usepackage[english]{babel}
 \def\botcaption#1#2{\medskip\centerline{{\scshape #1.}\kern8pt
 {\rm #2}}\bigskip}

  \newtheorem{ittheorem}{Theorem}[section]
 \newtheorem{itlemma}[ittheorem]{Lemma}
 \newtheorem{itproposition}[ittheorem]{Proposition}
 \newtheorem{itdefinition}[ittheorem]{Definition}
 \newtheorem{itremark}[ittheorem]{Remark}
 \newtheorem{itclaim}[ittheorem]{Claim}
 \newtheorem{itcorollary}[ittheorem]{Corollary}
 \numberwithin{equation}{section}

 \newenvironment{theorem}{
 \begin{ittheorem}}{\end{ittheorem}}

 \newenvironment{lemma}{
 \begin{itlemma}}{\end{itlemma}}

 \newenvironment{proposition}{
 \begin{itproposition}}{\end{itproposition}}

 \newenvironment{definition}{
 \begin{itdefinition}}{\end{itdefinition}}

 \newenvironment{remark}{
 \begin{itremark}}{\end{itremark}}

 \newenvironment{claim}{
 \begin{itclaim}}{\end{itclaim}}

\newenvironment{corollary}{
	\begin{itcorollary}}{\end{itcorollary}}

 \newenvironment{proof}{\noindent {\bf Proof.\,}
 }{\hspace*{\fill}$\qed$\medskip}

 \newenvironment{proof*}{\noindent {\bf Proof\,}
}{\hspace*{\fill}$\qed$\medskip}

 \newcommand{\be}[1]{\begin{equation}\label{#1}}
 \newcommand{\ee}{\end{equation}}

 \newcommand{\bl}[1]{\begin{lemma}\label{#1}}
 \newcommand{\el}{\end{lemma}}

 \newcommand{\br}[1]{\begin{remark}\label{#1}}
 \newcommand{\er}{\end{remark}}

 \newcommand{\bt}[1]{\begin{theorem}\label{#1}}
 \newcommand{\et}{\end{theorem}}

 \newcommand{\bd}[1]{\begin{definition}\label{#1}}
 \newcommand{\ed}{\end{definition}}

 \newcommand{\bcl}[1]{\begin{claim}\label{#1}}
 \newcommand{\ecl}{\end{claim}}

 \newcommand{\bp}[1]{\begin{proposition}\label{#1}}
 \newcommand{\ep}{\end{proposition}}

 \newcommand{\bc}[1]{\begin{corollary}\label{#1}}
 \newcommand{\ec}{\end{corollary}}

 \newcommand{\bpr}{\begin{proof}}
 \newcommand{\epr}{\end{proof}}

 \newcommand{\bi}{\begin{itemize}}
 \newcommand{\ei}{\end{itemize}}

 \newcommand{\ben}{\begin{enumerate}}
 \newcommand{\een}{\end{enumerate}}

\def\ss{{\cX^s}}
\def\sm{{\cX^m}}

 \def \ba {\begin{array}}
 \def \ea {\end{array}}

 \def \qed {{\square\hfill}}
 \def \Z {{\mathbb Z}}
 \def \R {{\mathbb R}}
 
 \def \N {{\mathbb N}}
 \def \P {{\mathbb P}}
 \def \E {{\mathbb E}}

 \def \ra {\rightarrow}

 \def \cS {{\cal S}}
 
 \def \cF {{\cal F}}
 \def \cA {{\cal A}}
 \def \cR {{\cal R}}
  \def \cG {{\cal G}}
\def \cI {{\cal I}}
\def \cL{{\cal L}}
 \def \cQ {{\cal Q}}
 \def \cC {{\cal C}}
 \def \cD {{\cal D}}
 \def \cX {{\cal X}}
 \def \cB {{\cal B}}
 
 \def \cW {{\cal W}}
 
 \def \cP {{\cal P}}
 \def \cL {{\cal L}}
 \def \cV {{\cal V}}
   \def \cZ {{\cal Z}}

 \def \G {{\Gamma}}
 \def \L {{\Lambda}}
 
 \def \a {{\alpha}}
 \def \b {{\beta}}
 \def \e {{\varepsilon}}
 \def \D {{\Delta}}
 \def \r {{\rho}}
\def \m {{\mu}}

 \def \h {{\eta}}
 \def \s {{\sigma}}
 \def \z {{\zeta}}
 
 \def \t {{\tau}}
 \def \o {{\omega}}
 \def \d {{\delta}}
\def \p {{\pi}}
\def \x {{\xi}}

 \def \CAPA {{\hbox{\footnotesize\rm CAP}}}

\def\pieno{{\blacksquare}}
\def\vuoto{{\square}}

\def\csgeo{{\cal C}_{sa}^*}

\def\gs{\G_{sa}^*}
\def\gw{\G_{wa}^*}

\begin{document}

\title{Critical Droplets and sharp asymptotics for Kawasaki dynamics \\ with strongly anisotropic interactions}
	
	\author[
	{}\hspace{0.5pt}\protect\hyperlink{hyp:email1}{1},\protect\hyperlink{hyp:affil1}{a}
	]
	{\protect\hypertarget{hyp:author1}{Simone Baldassarri}}
	
	\author[
	{}\hspace{0.5pt}\protect\hyperlink{hyp:email2}{2},\protect\hyperlink{hyp:affil1}{a,b}
	]
	{\protect\hypertarget{hyp:author2}{Francesca R.\ Nardi}}

	\affil[ ]{\centering
		\small\parbox{365pt}{\centering
			\parbox{5pt}{\textsuperscript{\protect\hypertarget{hyp:affil2}{a}}}Dipartimento di Matematica e Informatica ``Ulisse Dini", Universit\`{a} degli Studi di Firenze, Firenze, Italy.
		}
	}
	
	\affil[ ]{\centering
		\small\parbox{365pt}{\centering
			\parbox{5pt}{\textsuperscript{\protect\hypertarget{hyp:affil2}{b}}}Department of Mathematics and Computer Science, Eindhoven University of Technology, Eindhoven, the Netherlands.
		}
	}
	
	\affil[ ]{\centering
		\small\parbox{365pt}{\centering
			\parbox{5pt}{\textsuperscript{\protect\hypertarget{hyp:email1}{1}}}\texttt{\footnotesize\href{mailto:simone.baldassarri@unifi.it}{simone.baldassarri@unifi.it}},
			\parbox{5pt}{\textsuperscript{\protect\hypertarget{hyp:email2}{2}}}\texttt{\footnotesize\href{mailto:francescaromana.nardi@unifi.it}{francescaromana.nardi@unifi.it}},
		}
	}
	
	\maketitle
	
	\vspace{-1.25cm}
	
	\begin{center}
		{\it Unfortunately my coauthor Francesca Nardi passed away on 21 October 2021 \\ during the review process of the paper. I wish to thank her for the bright \\ person and talented mathematician she was.}
		
	\end{center}
	
	\begin{abstract}
			In this paper we analyze metastability and nucleation in the context of the Kawasaki dynamics for the two-dimensional Ising lattice gas at very low temperature. Let $\Lambda\subset\mathbb{Z}^2$ be a finite box. Particles perform simple exclusion on $\Lambda$, but when they occupy neighboring sites they feel a binding energy $-U_1<0$ in the horizontal direction and $-U_2<0$ in the vertical one. Thus the Kawasaki dynamics is conservative inside the volume $\Lambda$. Along each bond touching the boundary of $\Lambda$ from the outside to the inside, particles are created with rate $\rho=e^{-\Delta\beta}$, while along each bond from the inside to the outside, particles are annihilated with rate $1$, where $\beta>0$ is the inverse temperature and $\Delta>0$ is an activity parameter. Thus, the boundary of $\Lambda$ plays the role of an infinite gas reservoir with density $\rho$. We consider the parameter regime $U_1>2U_2$ also known as the strongly anisotropic regime. We take $\Delta\in{(U_1,U_1+U_2)}$, so that the empty (respectively full) configuration is a metastable (respectively stable) configuration. We consider the asymptotic regime corresponding to finite volume in the limit as $\beta\rightarrow\infty$.  We investigate how the transition from empty to full takes place with particular attention to the critical configurations that asymptotically have to be crossed with probability 1. The derivation of some geometrical properties of the saddles allows us to identify the full geometry of the minimal gates and their boundaries for the nucleation in the strongly anisotropic case. We observe very different behaviors for this case with respect to the isotropic ($U_1=U_2$) and weakly anisotropic ($U_1<2U_2$) ones. Moreover, we derive mixing time, spectral gap and sharp estimates for the asymptotic transition time for the strongly anisotropic case.

		\medskip
		{\it AMS} 2020 {\it subject classifications.} 60J10; 60K35; 82C20; 82C22; 82C26
		
		\medskip
		{\it Key words and phrases.} Lattice gas, Kawasaki dynamics, metastability, critical droplet, large deviations, pathwise approach, potential theory.
		
		\medskip
		{\it Acknowledgment.} F.R.N. was partially supported by the Netherlands Organisation for Scientific Research (NWO) [Gravitation Grant number 024.002.003--NETWORKS]. This work was supported by the ``Gruppo Nazionale per l'Analisi Matematica e le loro Applicazioni" (GNAMPA-INdAM). The authors are grateful to Anna Gallo, Vanessa Jacquier and Cristian Spitoni for useful and fruitful discussions.
		
			\end{abstract}

	\newpage
	\tableofcontents

\section{Introduction}
Metastability is a dynamical phenomenon that occurs when a thermodynamic system is close to a first order phase transition, that takes place when some physical parameter such as the temperature, pressure or magnetic field changes. The phenomenon of metastability for very low temperature dynamics is characterized by the tendency of the system to remain for a long time in a state (the metastable state $m$) different from the stable states $\cX^s$. Moreover, the system leaves this apparent equilibrium at some random time performing a sudden transition to the stable state. This transition is called {\it metastability} or {\it metastable behavior}. Metastability is an ubiquitous phenomenon with many examples from physical systems such as supersatured vapour, superheated and supercooled water, magnetic hysteresis loop, and from wireless networks. In the study of metastablity there are three main issues that are tipically investigated. The first one is the study of the {\it typical transition time} from the metastable to the stable state. The second and third issues, that are physically more interesting, concern the geometrical description of the {\it gate configurations} (also called {\it critical configurations}) and the {\it tube of typical trajectories}, that we will discuss in the sequel. A central role in these descriptions is played by the {\it gates} $\cW(m,\cX^s)$ from $m$ to $\cX^s$, that are sets of configurations typically visited during the last excursion from $m$ to $\cX^s$ (see Section \ref{modinddef} point 4 for the precise definition). A {\it minimal gate} has the physical meaning of ``collection of critical configurations" and it is defined as a gate such that removing any configuration from it, the new set has not the gate property. Because of this, the characterization of the union of minimal gates $\cG(m,\cX^s)$ (see (\ref{defg})) is important. The third issue concerns the identification of the so-called tube of typical trajectories. This is the set of typical paths followed by the system during the transition from the metastable to the stable state. We note that the hypotheses needed to discuss the gates are weaker than the ones necessary to completely characterize the tube of typical paths. The geometrical characterization of the union of minimal gates $\cG(m,\cX^s)$ is a central issue both from a probabilistic and from a physical point of view and it is a crucial point in the description of the typical trajectories. We remark that in several models proposed to describe ferromagnetic systems and analyzed in the literature in the context of Freidlin-Wentzell Markov chains evolving under Glauber dynamics, the minimal gate was unique but, in general, there may exist many minimal sets with the gate property, either distinct or overlapping. In order to model mathematically phenomena such as superheated or supercooled water is often proposed the use of lattice gas models evolving according to Kawasaki dynamics since the dynamics conserves the number of particles.

In this paper we consider the metastable behavior of the two-dimensional Ising lattice gas with strongly anisotropic interactions at very low temperature and low density that evolves under Kawasaki dynamics, i.e., a discrete time Markov chain defined by the Metropolis algorithm with transition probabilities given in (\ref{defkaw}). Let $\b>0$ be the inverse temperature and let $\L\subset\Z^2$ be a finite box with open boundary conditions. Particles live and evolve in a conservative way inside $\L$, but when they occupy neighboring sites they feel a binding energy $-U_1<0$ in the horizontal direction and $-U_2<0$ in the vertical one (see Section \ref{S1.1} for more details). Without loss of generality we may assume $U_1\geq U_2$. Along each bond touching the boundary of $\Lambda$ from the outside to the inside, particles are created with rate $\rho=e^{-\Delta\beta}$, while along each bond from the inside to the outside, particles are annihilated with rate $1$, where $\Delta>0$ is an activity parameter. Thus, the boundary of $\Lambda$ plays the role of an infinite gas reservoir with density $\rho$. We fix the parameters $U_1$, $U_2$ and $\D$ such that $U_1>2U_2$ in what we call the {\it strongly anisotropic case}. We take $\Delta\in{(U_1,U_1+U_2)}$, so that the empty (respectively full) configuration is the metastable (respectively stable) state. We consider the asymptotic regime corresponding to finite volume $\L$ in the limit of large inverse temperature $\beta$. We investigate how the system {\it nucleates}, i.e., how it reaches $\pieno$ (box full of particles) starting from $\vuoto$ (empty box).

One of the main goals of the paper is to investigate this two-dimensional model in the strongly anisotropic case giving the geometrical description of the set $\cG(\vuoto,\pieno)$ in Theorem \ref{gstrong}. We will prove that there are many distinct minimal gates that we will geometrically characterize together with their union. Let us explain the strategy we adopt in our paper. In \cite[Theorem 5.1]{MNOS} there is a characterization of the set $\cG(m,\cX^s)$ in terms of {\it essential saddles} (see Section \ref{modinddef} point 4 for the definition of essential saddle). Thanks to this equivalence, we reduce our study to the identification of the set of all the essential saddles that has to be crossed during the transition between the metastable state $\vuoto$ and the stable state $\pieno$. We apply the model-independent strategy carried out in \cite[Section 3.1]{BN2} to the strongly anisotropic two-dimensional case, where $m=\vuoto$ and $\cX^s=\{\pieno\}$, in order to eliminate some unessential saddles. Thus we need to verify that the required model-dependent inputs are valid in our case. This study together with the characterization of the essential saddles rely on a detailed analysis of the motion of particles along the border of the droplet. On the one hand this is a typical feature of the Kawasaki dynamics, on the other hand this is peculiar in the strongly anisotropic case. Indeed in the metastable regime, particles move along the border of a droplet more rapidly than they arrive from the boundary of the box. More precisely, before the arrival of the next particle, we have that single particles attached to one side of a droplet tipycally detach (because $e^{U_1\b}\ll e^{\D\b}$ and $e^{U_2\b}\ll e^{\D\b}$), while bars of two or more particles tipycally do not detach (because $e^{\D\b}\ll e^{(U_1+U_2)\b}$). Roughly speaking, we will investigate the saddles that are crossed ``just before visiting" and ``just after visiting" the gate of the transition $\cG(\vuoto,\pieno)$. Additionally, we prove sharp asymptotics for the transition time in Theorem \ref{sharptimestrong} and we investigate the spectral gap and mixing time in Theorem \ref{gapspettrale}.

Some properties of the metastable behavior for the strongly anisotropic case have been already derived in the literature. More precisely, in \cite{BN} the authors derived the asymptotic behavior of the transition time in probability, law and expectation. Additionally, they gave the geometrical description of a gate in \cite[Theorem 2.4]{BN}: we improve that statement in Theorem \ref{thgate}. Moreover, our Theorem \ref{gstrong} gives a more detailed description of the geometry of the minimal gates, their union and the entrance in it with respect to their results.

For the isotropic interactions, i.e., $U_1=U_2$, in \cite{HOS} the authors investigated the asymptotic properties of the transition time together with an intrinsic description of a gate (see Section \ref{modinddef} point 4 for the precise definition). This paper initiated the study of this model that we describe in the following discussion. In \cite{BHN} a geometric characterization of a subset of $\cG(\vuoto,\pieno)$ is given and this result is improved in \cite[Theorem 4.2]{BN2} with the identification of the minimal gates and their union $\cG(\vuoto,\pieno)$. For the three-dimensional lattice gas we refer to \cite{HNOS}, where the authors investigated the asymptotic properties of the transition time and an intrinsic description of a gate. Moreover, for both two and three-dimensional isotropic case, using the {\it potential theoretic approach} the authors investigated in \cite{BHN} the sharp asymptotics of the mean transition time, the so-called {\it pre-factor}. They proved that it is a constant that asymptotically depends only on the size of the box and the cardinality of the gate that they identified, but not on the parameter $\b$. In the framework of the {\it pathwise approach} it is natural to study the third issue of metastability, namely the tube of typical trajectories realizing the transition between $\vuoto$ and $\pieno$.
This has been analyzed only in \cite{GOS} for two dimensions, indeed actually there are no known results about the tube for the three-dimensional isotropic case and for the anisotropic one. Concerning the weakly anisotropic case, i.e., $U_1<2U_2-2\e$ with $\e:=U_1+U_2-\D$, the asymptotic behavior of the transition time in probability, law and expectation has been derived in \cite{NOS}. \cite[Theorems 4.5, 4.6]{BN2} give a more detailed description of the geometry of the minimal gates, their union and the entrance in it with respect to the geometric description of a gate given in \cite{NOS}

The geometrical analysis of the union of minimal gates with their boundary for the isotropic and weakly anisotropic cases is given in \cite{BN2} and we discuss the differences and similarities below. Despite the structure of the gate is similar for the three cases, we emphasize that the entrance in them is very different. In particular, for the strongly anisotropic case there are two different mechanisms to enter the gate (see Lemma \ref{entratastrong}), while for the other two cases there is a unique one (see \cite[Lemma 7.13]{BN2}). This is a consequence of a larger rigidity of the dynamics in the strongly anisotropic case: an important part of the regularizing motions of particles along the border of the clusters is lost in such a way that a new mechanism of entering in the critical configurations set appears, with some impact on the prefactor of the mean transition time to stability. On the other hand, it is clear that the properties that are strictly related to the horizontal and vertical interactions are the same for both weakly and strongly anisotropic cases. While some properties that involve the motion of particles along the border of the droplet are very different. Intuitively, one may think of the weakly anisotropic case as an interpolation between the isotropic and strongly anisotropic ones. Indeed, it has some properties similar to the first, others to the latter. This specific difference between these cases motivates together with applications the rigorous investigation of the anisotropic cases. Moreover, we highlight this difference in the description of the set $\cG(\vuoto,\pieno)$, indeed for the isotropic case more motions along the border are allowed and thus a totally explicit geometric description of the set is more difficult (see \cite[Theorem 4.2]{BN2}), but for the anisotropic cases we fully obtain it, since the condition $U_1\neq U_2$ makes more difficult the sliding of particles along the border of the droplet. Among the anisotropic cases, by \cite[Theorem 4.6]{BN2} and Theorem \ref{gstrong} it is clear that the structure of the set $\cG(\vuoto,\pieno)$ strongly depends on how large is $U_1$ with respect to $U_2$, indeed in the case $U_1>2U_2$ less slidings along the border are allowed and thus the structure of the union of minimal gates is less rich than the weak anisotropic case.

{\bf State of the art.} The first dynamical approach, known as \emph{pathwise approach}, was initiated in 1984 in \cite{CGOV}, developed in \cite{OS,OS2} and summerized in the monograph \cite{OV}. For Metropolis chains associated with statistical mechanics systems, metastability has been described by this approach in an elegant way in terms of the energy landscape associated to the Hamiltonian of the system. This approach focuses on the {\it dynamics} of the transition from metastable to stable states and it is so flexible that has been later developed to treat the tunnelling, namely the transition from a stable state to another stable state or stable states. Independently, a graphical approach was introduced in \cite{CC} and later used for Ising-like models \cite{CaTr}. Using the pathwise approach it is possible to obtain a detailed description of metastable behavior of the system and it made possible to answer all the three questions of metastability. A modern version of the pathwise approach can be found in \cite{{MNOS},{CNbc},{CNS2},{BNZ}}. In particular, in \cite{MNOS}, for the Metropolis markov chains, there are model-independent results concerning the transition time in probability, expectation and distribution, and concerning minimal gates and their union disentangled with respect to the tube of typical trajectories. In \cite{MNOS} the results on hitting times are obtained with minimal model-dependent knowledge,
i.e., find all the metastable states and the minimal energy barrier which separates them from the stable states. In \cite[Sections 2-3]{CNbc} the authors prove model-independent results to treat systems with multiple metastable states and give a sufficient condition to identify them. In \cite{CNS2} the authors extend the results of \cite{MNOS} to general Markov chains (reversible and non reversible) with rare transitions setup, also called Freidlin-Wentzel Markov chains. These results are a useful tool to approach metastability for non-Metropolis systems such as Probabilistic Cellular Automata. In \cite[Section 3]{BNZ} the authors extended the model-independent framework of \cite{MNOS} to study the first hitting times from any starting configuration (not necessarily metastable) to any target subset of configurations (not necessarily the set of stable configurations). This approach developed over the years has been extensively applied to study metastability in Statistical Mechanics lattice models. In this context, this approach and the one that follows (\cite{BEGK,MNOS,OV}) have been developed with the aim of finding answers valid with maximal generality and to reduce as much as possible the number of model dependent inputs necessary to describe the metastable behavior of any given system. 

Another approach is the \emph{potential-theoretic approach}, initiated in \cite{BEGK}. We refer to \cite{BH} for an extensive discussion and applications to different models. In this approach, the metastability phenomenon is interpreted as a sequence of visits of the path to different metastable sets. This method focuses on a precise analysis of hitting times of these sets with the help of \emph{potential theory}. In the potential-theoretic approach the mean transition time is given in terms of the so-called \emph{capacities} between two sets. Crucially capacities can be estimated by exploiting powerful variational principles. This means that the estimates of the average crossover time that can be derived are much sharper than those obtained via the pathwise approach. 

These mathematical approaches, however, are not equivalent as they rely on different definitions of metastable states (see \cite[Section 3]{CNbc} for a comparison) and thus involve different properties of hitting and transition times. The situation is particularly delicate for evolutions of infinite-volume systems, for irreversible systems, and degenerate systems, i.e., systems where the energy landscape has configurations with the same energy (as discussed in \cite{CNbc,CNS2,CNS2017}). More recent approaches are developed in \cite{BL1,BL2,BiGa}.

Statistical mechanical models for magnets deal with dynamics that do not conserve the total number of particles or the total magnetization. They include single spin-flip Glauber dynamics and many probabilistic cellular automata (PCA), that is a parallel dynamics. The pathwise approach was applied in finite volume at low temperature in \cite{CGOV,NevSchbehavdrop,CaTr,KOd,KOs,CO,NO,CL,AJNT,BGNneg2021,BGNpos2021} for single-spin-flip Glauber dynamics and in \cite{CN,CNSp,CNS22,CNS2016} for parallel dynamics. The potential theoretic approach was applied to models at finite volume and at low temperature in \cite{BM,BHN,HNT1,HNT2,NS1,HNTA2018,BJN}. The more involved infinite volume limit at low temperature or vanishing magnetic field was studied in \cite{DS1,DS2,S2,SS,MO1,MO2,HOS,GHNOS,GN,BHS,CeMa,GMV} for Ising-like models under single-spin-flip Glauber and Kawasaki dynamics.

The outline of the paper is as follows. In Section \ref{S2} we define the model with open boundary conditions and the Kawasaki dynamics. In Section \ref{modinddef} we give some model-independent definitions and in Section \ref{S4} we give some geometric definitions valid for Kawasaki dynamics (see Section \ref{moddepdef}). We state our main results concerning the gates in Section \ref{sani} and about the sharp asymptotics in Section \ref{sharpestimates}. In Section \ref{dependentdef} we apply the model-independent strategy carried out in \cite[Section 3.1]{BN2}. In Section \ref{sitigood} we give some model-dependent definitions, in Section \ref{lemmi3modelli} some tools that are useful in Section \ref{S6.4} for our model-dependent strategy. In Section \ref{proofstrong} we give the proof of the main results for the strongly anisotropic case regarding the description of the gate (see Theorem \ref{thgate}) and the geometric characterization of the union of all the minimal gates (see Theorem \ref{gstrong}). In Section \ref{sharpasymptotics} we give the proof of the main theorems about the sharp asymptotics (see Theorems \ref{sharptimestrong} and \ref{gapspettrale}).

\section{Definition of the model}
\label{S2}

\subsection{The model with open boundary conditions}
\label{S1.1}

Let $\Lambda=\{0,..,L\}^2\subset \Z^2$ be a finite box centered at the
origin. The side length $L$ is fixed, but arbitrary, and later we will require $L$ to be sufficiently large. Let
\be{inbd}
\partial^- \L:= \{x\in\L: \exists\; y \notin\L\: |y-x|=1\},
\ee

\noindent be   the interior  boundary of $ \Lambda$ and let
$ \Lambda_0:= \Lambda\setminus\partial^- \Lambda$ be the interior of $\L$.
With each $x\in \Lambda$ we associate an occupation variable
$\eta(x)$, assuming values 0 or 1. A lattice configuration is
denoted by $\eta\in {\cal X} =\{ 0,1\} ^{ \Lambda }$. Each configuration $\h\in \cX$ has an energy given by the following Hamiltonian:

\be{hamilt} H(\eta):= -U_1 \sum_{(x,y)\in   \Lambda_{0,h}^{*}}
\eta(x)\eta(y) -U_2\sum_{(x,y)\in   \Lambda_{0,v}^{*}} \eta(x)
\eta(y)+ \D \sum _{x\in \L} \eta (x), \ee

\noindent
where $ \Lambda_{0,h}^{*}$ (resp.\ $ \Lambda_{0,v}^{*}$) is the set of
the horizontal (resp.\ vertical) unoriented  bonds joining nearest-neighbors points in
$ \Lambda_0$. Thus the interaction is acting only inside $
\Lambda_0$; the binding energy associated to a horizontal
(resp.\ vertical) bond is $-U_1<0$ (resp.\ $-U_2<0$). We may assume without
loss of generality that $U_1\ge U_2$. 

The grand-canonical Gibbs measure associated with $H$ is
\be{misura} \m(\eta):= {  e^{- \b H(\eta) }\over Z} \qquad \h\in
\cX, \ee

\noindent
where
\be{partfunc} Z:=\sum_{\eta\in {\cal X}}e^{-\b H(\eta)}
\ee

\noindent
is the so-called {\it partition function}.

\subsection{Local Kawasaki dynamics}
\label{S1.2}

Next we define Kawasaki dynamics on $\L$ with boundary
conditions that mimic the effect of an infinite gas reservoir
outside $\L$ with density $ \r = e^{-\D\b}.$ Let $b=(x \to y)$ be
an oriented bond, i.e., an {\it ordered} pair of nearest neighbour
sites, and define

\be{Loutindef}
\ba{lll}
\partial^* \L^{out} &:=& \{b=(x \to y): x\in\partial^- \L,
y\not\in\L\},\\
\partial^* \L^{in}  &:=& \{b=(x \to y): x\not\in
\L, y\in\partial^-\L\},\\
\L^{*, orie} &:=& \{b=(x \to y): x,y\in\L\},
\ea
\ee

\noindent and put $ \bar\L^{*, orie}:=\partial^* \L ^{out}\cup
\partial^* \L ^{in}\cup\L^{*,\;orie}$.
Two configurations $  \eta,
\eta'\in {\cal X}$ with $ \eta\ne \eta'$ are said to be {\it
	communicating states} if there exists a bond
$b\in  \bar\L^{*,orie}$ such that $ \eta' = T_b \eta$, where $T_b   \eta$ is the
configuration obtained from $ \eta$ in any of these ways:

\begin{itemize}
	
	\item
	for $b=(x \to y)\in\L^{*,\;orie}$, $T_b \eta$ denotes the
	configuration obtained from $ \eta$ by interchanging particles
	along $b$:
	\be{Tint}
	T_b \h(z) =
	\left\{\ba{ll}
	\h(z) &\mbox{if } z \ne x,y,\\
	\h(x) &\mbox{if } z = y,\\
	\h(y) &\mbox{if } z = x.
	\ea
	\right.
	\ee
	
	\item
	For  $b=(x \to y)\in\partial^*\L^{out}$ we set:
	\be{Texit}
	T_b \h(z) =
	\left\{\ba{ll}
	\h(z) &\mbox{if } z \ne x,\\
	0     &\mbox{if } z = x.
	\ea
	\right.
	\ee
	
	\noindent
	This describes the annihilation of a particle along the border;
	
	\item
	for  $b=(x  \to y)\in\partial^*\L^{in}$ we set:
	\be{Tenter}
	T_b \h(z) =
	\left\{\ba{ll}
	\h(z) &\mbox{if } z \ne y,\\
	1     &\mbox{if } z=y.
	\ea
	\right.
	\ee
	
	\noindent
	This describes the creation of a particle along the border.
	
\end{itemize}

\noindent
The Kawasaki dynamics is  the discrete time Markov chain
$(\eta_t)_{t\in \mathbb{N}}$ on state space $ {\cal X} $ given by
the following transition  probabilities: for  $  \eta\not= \eta'$:
\be{defkaw}
P( \eta,  \eta'):=\left\{\ba{ll}
{ |\bar\L^{*,\;orie}|}^{-1} e^{-\b[H( \eta') - H( \eta)]_+}
&\mbox{if }  \exists b\in \bar\L ^{*, orie}: \eta' =T_b \eta   \\
0   &\mbox{ otherwise }  \ea \right.
\ee

\noindent
where $[a]_+ =\max\{a,0\}$ and $P(\h,\h):=1-\sum_{\h'\neq\h}P(\h,\h')$.
This describes a standard Metropolis dynamics with open boundary conditions: along each
bond touching $\partial^-\L$ from the outside, particles are created with
rate $\rho=e^{-\D\b}$ and are annihilated with rate 1, while inside $\L_0$ particles
are conserved.
Note that an exchange of occupation
numbers $\h(x)$ for any $x$ inside the ring $ \L\setminus  \L_0$
does not involve any change in energy.\par

\br{p1}
The stochastic dynamics defined by
(\ref{defkaw}) is reversible w.r.t.\ Gibbs measure (\ref{misura}) corresponding
to $ H$.
\er

\br{simpmodel}
The analysis of the fully conservative model, namely Kawasaki dynamics inside a large box $\L^\b\subset\Z^2$, with periodic boundary conditions and $\L\subset\L^\b$, such that $\lim_{\b\ra\infty}\frac{1}{\b}\log|\L^\b|=\infty$, is out of the scope of this paper. An extension of the model considered in Section \ref{S1.1} that goes in this direction is what we call \emph{simplified model}, in which we consider interactions only inside $\L_0=\L\setminus\partial^-\L$, where $\partial^-\L$ is defined in (\ref{inbd}), while in $\L\setminus\L_0$ we remove interactions and in $\L^\b\setminus\L$ we remove both interactions and exclusion so that the dynamics of the gas outside $\L$ is that of independent random walks. 

Following the strategy proposed in \cite{HOS} for $int=is$, we are able to derive results concerning the transition time, the gate and supercritical and subcritical rectangles for the strongly anisotropic simplified model similar to the one derived in \cite[Theorem 1.53]{HOS} for the isotropic case. In \cite[Section 2]{HOS} the authors give several large deviation estimates concerning exponential clocks, that hold also for the anisotropic cases. In \cite[Section 3]{HOS} the authors give several large deviation estimates concerning random walks. All these results are valid for the anisotropic cases (both strong and weak) without changes except for \cite[Proposition 3.13]{HOS}, in which we have to replace $U$ with $U_1$. The recurrence property for the anisotropic simplified model is obtained with similar arguments carried out in \cite[Section 6]{HOS}. To this end, we modify the definition of the set $\bar\cX_2$ given in \cite[eq.\ (5.8)]{HOS} by replacing $U$ with $U_1$. Therefore also the definition of the set $\cX_2$ given in \cite[eq.\ (6.1)]{HOS} should be modified accordingly. Thus, if we define for the anisotropic model $T_1=e^{0\b}$, $T_2=e^{U_1\b}$ and $T_3=e^{\D\b}$, \cite[Proposition 6.2]{HOS} holds also for the anisotropic cases. Concerning the reduction, we follow the strategy proposed in \cite[Section 7]{HOS}. In particular, we have to study the behavior of the gas and its interaction with the dynamics in the box $\L$. There are two classes of gas particles with different behavior: particles that have been in $\L^\b\setminus\L$ for a long time (say of order $T_3$), which we call green particles; and particles that exit from $\L$ and afterwards return to $\L$ in a short time (say of order 1), which we call red particles. The effect of green (resp.\ red) particles is studied in \cite[Section 7.6]{HOS} (resp.\ \cite[Section 7.7]{HOS}) and can be extended to the anisotropic cases by modifying the times $T_1$, $T_2$ and $T_3$, and the sets $\cX_2$ and $\bar\cX_2$ as above. In the strongly anisotropic case, from this discussion and \cite[Theorems 2.3, 2.4 and 2.5]{BN} the desired results follow. For the geometrical characterization of the critical droplets in the simplified model we refer to Remark \ref{simpgates}.
\er

\section{Model-independent definitions and notations}
\label{modinddef}
We will use italic capital letters for subsets of $\L$, script
capital letters for subsets of $\cX$, and boldface capital letters for events
under the Kawasaki dynamics. We use this convention in order to keep the
various notations apart. We will denote by $\P_{\h_0}$ the probability law of the Markov
process  $(\eta_t)_{t\geq 0}$ starting at $\h_0$ and by
$\E_{\h_0}$ the corresponding expectation.

\medskip
\noindent
{\bf 1. Paths and hitting times.}
\bi

\item
A {\it path\/} $\o$ is a sequence $\o=(\o_1,\dots,\o_k)$, with
$k\in\N$, $\o_i\in\cX$ and $P(\o_i,\o_{i+1})>0$ for $i=1,\dots,k-1$.
We write $\o\colon\;\h\to\h'$ to denote a path from $\h$ to $\h'$,
namely with $\o_1=\h,$ $\o_k=\h'$. A set
$\cA\subset\cX$ with $|\cA|>1$ is {\it connected\/} if and only if for all
$\h,\h'\in\cA$ there exists a path $\o:\h\to\h'$ such that $\o_i\in\cA$
for all $i$. We indicate with $\o_1\circ\o_2$ the composition of two paths $\o_1$ and $\o_2$, namely if $\o_1=(\o_1^1,...,\o_k^1)$ and $\o_2=(\o_1^2,...,\o_m^2)
$ then $\o_1\circ\o_2=(\o_1^1,...,\o_k^1,\o_1^2,...,\o_m^2)$.

\item[$\bullet$]
Given a non-empty
set $\cA\subset\cX$, define the {\it first-hitting time of} $\cA$
as
\be{tempo}
\t_{\cA}:=\min \{t\geq 0:  \eta_t \in \cA \}.
\ee
\ei
\noindent

\medskip
\noindent
{\bf 2. Min-max and communication height}
\bi

\item Given a function $f:\cX\ra\R$ and a subset $\cA\subseteq\cX$, we denote by 
\be{defargmax}
\arg \hbox{max}_{\cA}f:=\{\h\in\cA: f(\h)=\max_{\z\in\cA}f(\z)\}
\ee

\noindent
the set of points where the maximum of $f$ in $\cA$ is reached. If $\o=(\o_1,...,\o_k)$ is a path, in the sequel we will write $\arg \max_{\o}H$ to indicate $\arg \max_{\cA}H$, with $\cA=\{\o_1,...,\o_k\}$ e $H$ the Hamiltonian.

\item 
The {\it bottom} $\cF(\cA)$ of a  non-empty
set $\cA\subset\cX$ is the
set of {\it global minima} of the Hamiltonian $H$ in  $\cA$:
\be{Fdef}
\cF(\cA):=\arg \hbox{min}_{\mathcal{A}}H= \{\h\in\cA: H(\h)=\min_{\z\in\cA} H(\z)\}.
\ee
For a set $\cA\subset\cX$ such that all the configurations have the same energy, with an abuse of notation we denote this energy by $H(\cA)$.

\item 
The {\it communication height} between a pair $\h$, $\h'\in\cX$ is
\be{}
\Phi(\h,\h'):= \min_{\o:\h\ra\h'}\max_{\z\in\o} H(\z).
\ee
\noindent
Given $\cA\subset\cX$, we define the {\it restricted communication height} between $\h,\h'\in\cA$ as
\be{}
\Phi_{|\cA}(\h,\h'):= \min_{\o:\h\ra\h'\atop\o\subseteq\cA}
\max_{\z\in\o} H(\z),
\ee
\noindent
where $(\o_1,...,\o_k)=\o\subseteq\cA$ means $\o_i\in\cA$ for every $i$.
\ei

\medskip
\noindent
{\bf 3. Stability level, stable and metastable states}
\bi

\item
We call
{\it stability level} of
a state $\z \in \cX$
the energy barrier
\be{stab}
V_{\z} :=
\Phi(\z,\cI_{\z}) - H(\z),
\ee

\noindent where $\cI_{\z}$ is the set of states with
energy below $H(\z)$:
\be{iz} \cI_{\z}:=\{\eta \in \cX: H(\eta)<H(\z)\}. \ee

\noindent 
We set $V_\z:=\infty$ if $\cI_\z$ is empty.

\item
We call {\it $V$-irreducible states}
the set of all states with stability level larger than  $V$:
\be{xv} 
\cX_V:=\{\h\in\cX: V_{\h}>V\}. 
\ee

\item
The set of {\it stable states} is the set of the global minima of
the Hamiltonian:
\be{st.st.}
\cX^s:=\cF(\cX). 
\ee

\item
The set of {\it metastable states} is given by
\be{st.metast.} \sm:=\{\h\in\cX:
V_{\h}=\max_{\z\in\cX\backslash \ss}V_{\z}\}. \ee
\noindent
We denote by $\G_m$ the stability level of the states in $\cX^m$.

\ei

\medskip
\noindent
{\bf 4. Optimal paths, saddles and gates}
\bi

\item 
We denote by $(\h\to\h')_{opt} $ the {\it set of optimal paths\/} as the set of all
paths from $\h$ to $\h'$ realizing the min-max in $\cX$, i.e.,
\be{optpath}
(\h\to\h')_{opt}:=\{\o:\h\to\h'\; \hbox{such that} \; \max_{\xi\in\o} H(\xi)=  \Phi(\h,\h') \}.
\ee

\item
The set of {\it minimal saddles\/} between
$\h,\h'\in\cX$
is defined as
\be{minsad}
\cS(\h,\h'):= \{\z\in\cX\colon\;\; \exists\o\in (\h\to\h')_{opt},
\ \o\ni\z \hbox{ such that } \max_{\xi\in\o} H(\xi)= H(\z)\}.
\ee

\item A saddle $\x\in\cS(\h,\h')$ is called {\it unessential} if for any $\o\in(\h\ra\h')_{opt}$ such that $\o\cap\x\neq\emptyset$ we have $\{\arg\max_{\o}H\}\setminus\{\x\}\neq\emptyset$ and there exists $\o'\in(\h\ra\h')_{opt}$ such that $\{\arg\max_{\o'}H\}\subseteq\{\arg\max_{\o}H\}\setminus\{\x\}$.

\item A saddle $\x\in\cS(\h,\h')$ is called {\it essential} if it is not unessential, i.e., if either

\bi
\item[(i)] there exists $\o\in(\h\ra\h')_{opt}$ such that $\{\hbox{arg max}_{\o}H\}=\{\x\}$ or
\item[(ii)] there exists $\o\in(\h\ra\h')_{opt}$ such that $\{\hbox{arg max}_{\o}H\}\supset\{\x\}$ and $\{\hbox{arg max}_{\o'}H\}\nsubseteq\{\hbox{arg max}_{\o}H\}\setminus\{\x\}$ for all $\o'\in(\h\ra\h')_{opt}$.
\ei

\item
Given a pair $\h,\h'\in\cX$,
we say that $\cW\equiv\cW(\h,\h')$ is a {\it gate\/}
for the transition $\h\to\h'$ if $\cW(\h,\h')\subseteq\cS(\h,\h')$
and $\o\cap\cW\neq\emptyset$ for all $\o\in (\h\to\h')_{opt}$. In words, a gate is a subset of $\cS(\h,\h')$ that is visited by all optimal paths.

\item
We say that $\cW(\h,\h')$ is a {\it minimal gate\/} for the transition
$\h\to\h'$ if it is a gate and for any $\cW'\subsetneq \cW(\h,\h')$ there
exists $\o'\in (\h\to \h')_{opt}$ such that $\o'\cap\cW'
=\emptyset$. In words, a minimal gate is a minimal subset of $\cS(\h,\h')$ by inclusion that is visited by all optimal paths.

\item
For a given pair of configurations $\h,\h'$, we denote by $\cG(\h,\h')$ the union of all minimal gates:
\be{defg}
\cG(\h,\h'):=\displaystyle\bigcup_{\cW(\h,\h') \hbox{ minimal gate}} \cW(\h,\h')
\ee

\ei

\section{Main results: the gates for our model}
\label{S4}
In this Section we state our main results: in Section \ref{sani} we obtain the geometrical characterization of the union of all minimal gates for the strongly anisotropic case. In order to do this we need some model-dependent definitions for the Kawasaki dynamics (see Section \ref{moddepdef}) and some specific ones for the strongly anisotropic case (see Section \ref{sani}). In Section \ref{sharpestimates} we derive sharp estimates for the asymptotic transition time. Moreover, we derive the mixing time and spectral gap. For the corresponding results obtained in the isotropic and weakly anisotropic cases, i.e., in the parameter regime $U_1=U_2$ and $U_1<2U_2-2\e$ respectively, where $\e$ is defined in (\ref{defepsilon}), we refer to \cite[Sections 4.2,4.3]{BN2} for the results concerning the gates and union of minimal gates and to \cite[Section 4.4]{BN2} for the results concerning the asymptotic transition time, mixing time and spectral gap.

\subsection{Geometric definitions for Kawasaki dynamics}
\label{moddepdef}
We give some model-dependent definitions and notations in order to state our main theorems.

\medskip
\noindent
{\bf 1. Free particles and clusters}

\bi

\item[$\bullet$]
For $x\in\L_0$, let $ \hbox{nn}(x):=\{ y\in \L_0\colon\;d(y,x)=1\}$ be the set of
nearest-neighbor sites of $x$ in $\L_0$, where $d$ in the entire paper denotes the lattice distance.
\item[$\bullet$]
A {\it free particle\/} in $\h\in\cX$ is a site $x$, with $\h(x)=1$, such that either $x\in\partial^-\L$, or $x\in\L_0$ and $\sum_{y\in nn(x)\cap\L_0}\eta(y)$ $=0$. We denote by $\h_{fp}$ the union of free particles in $\partial^-
\L$ and  free particles in $\L_0$. We denote by $n(\h)$ the number of free particles in $\h$.

We denote by $\h_{cl}$ the clusterized
part of the occupied sites of $\h$:
\be{hcl}\h_{cl} :=\{x\in\L_0: \ \h(x)=1\}\setminus\h_{fp}.
\ee

\item[$\bullet$]
We denote by $\h^{fp}$ the addition of a free particle anywhere in $\L$ to the configuration $\h$.

\item[$\bullet$]
Given a configuration $\h \in\cX$, consider the subset $C(\h_{cl})$ of $\R^2$ defined as the union of the $1\times 1$ closed squares centered at the occupied sites of $ \h_{cl}$ in $ \L_0$ and call the maximal connected components of this set the clusters of $\h_{cl}$.

\item[$\bullet$]
Given a set $A\subset\R^2$, we define as $|A|$ the number of $1\times1$ closed occupied squares in $A$ and as $||A||$ the numbers of $1\times1$ closed squares in $A$. Note that $||\cdot||$ takes into account the possibility that the squares are occupied or not.

\ei

\medskip
\noindent
{\bf 2. Projections, semi-perimeter and vacancies}

\bi

\item[$\bullet$]
For $\h\in\cX$, we denote by $g_1(\eta)$ (resp.\ $g_2(\eta)$) one half of the horizontal (resp.\ vertical) length of the Euclidean boundary of $C(\h_{cl})$. Recall the definition of $n(\h)$ given in Section \ref{moddepdef} point 1. Then the energy associated with $\h$ is given by
\be{Hcont}
H(\h) = - (U_1+U_2-\D)|C(\eta_{cl})| + {U_1}
g_2(\eta)+
{U_2} g_1(\eta)
+\D n(\h).
\ee

\item[$\bullet$]
Let $p_1(\h)$
and $p_2(\h)$ be the total lengths of  horizontal and vertical
projections of $ C(\h_{cl})$ respectively. More precisely, let
$r_{j,1}=\{x \in \Z^2:(x)_1=j\}$ be the $j$-th column and
$r_{j,2}=\{x \in \Z^2:(x)_2=j\}$ be the $j$-th row, where
$(x)_1$ or $(x)_2$ denote the first or second component of $x$. Let

\be{proie1} \p_1(\h):=\{j \in \Z:\, r_{j,1}\cap
C(\h_{cl})\not=\emptyset\} \ee

\noindent and $p_1(\h):=|\p_1(\h)|$. In a similar way we define
the vertical projection $\p_2(\h)$ and $p_2(\h)$.

\item[$\bullet$]
We define $g'_i(\h):= g_i(\h)- p_i(\h)\ge 0$; we call {\it monotone} a configuration such that $g_i(\h)= p_i(\h)$ for $i=1,2$.

\item[$\bullet$]
We define the {\it semi-perimeter} $s(\h)$ and the {\it vacancies} $v(\h)$ as
\be{defdep1}
\ba{lll}
s(\h)&:=& p_1(\h)+p_2(\h),\\
v(\h)&:=& p_1(\h)p_2(\h)- |C(\h_{cl})|.
\ea
\ee

\ei

\medskip
\noindent
{\bf 3. $\boldsymbol{n}$-manifold, rectangles and corners}

\bi 
\item[$\bullet$]
The configuration space $\cX$ can be partitioned as
\be{}
\cX=\displaystyle\bigcup_{n}\cV_n,
\ee
\noindent
where $\cV_n:=\{\h\in\cX: |C(\h_{cl})|+n(\h)=n\}$ is the set of configurations with $n$ particles, called the {\it $n$-manifold}.

\item[$\bullet$]
We denote by $\cR(l_1,l_2)$ the set of configurations that have no free particle and a single cluster such that $C(\h_{cl})$ is a rectangle $R(l_1,l_2)$, with $l_1,l_2\in\N$. For any $\h,\h'\in\cR (l_1,l_2)$ we have immediately:
\be{enerettan} H(\h)=H(\h')=H(\cR(l_1,l_2))={U_1}l_2 +{U_2} l_1
-\varepsilon l_1 l_2, 
\ee
\noindent
where
\be{defepsilon}
\e:= U_1+U_2 -\D.
\ee

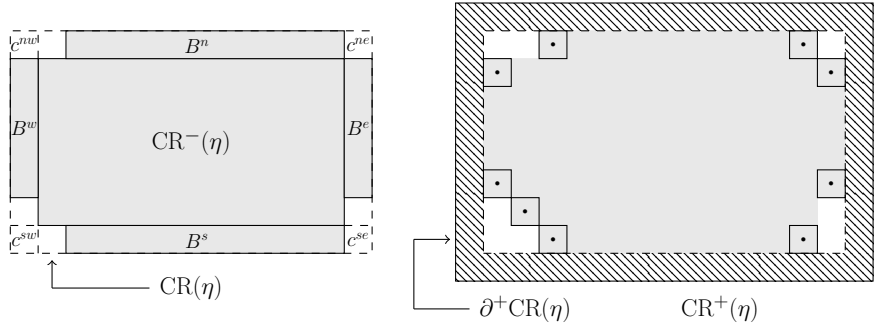
\begin{figure}
	\centering
	\begin{tikzpicture}[scale=0.37,transform shape]

		\draw [fill=grigio] (6,3) rectangle (17,9);
		\draw [fill=grigio] (17,4) rectangle (18,9);
		\node at (17.5,6.5){\huge{$B^e$}};
		\node at (11.5,6){\Huge{$\hbox{CR}^-(\h)$}};
		\draw [fill=grigio] (7,9) rectangle (17,10);
		\node at (11.7,9.45){\huge{$B^n$}};
		\draw [fill=grigio] (5,4) rectangle (6,9);
		\node at (5.5,6.5){\huge{$B^w$}};
		\draw [fill=grigio] (7,3) rectangle (17,2);
		\node at (11.7,2.45){\huge{$B^s$}};
		\draw [dashed] (5,2) rectangle (18,10);
		\draw [dashed] (5,3)--(6,3);
		\draw [dashed] (6,2)--(6,3);
		\node at (5.5,2.5){\huge{$c^{sw}$}};
		\draw [dashed] (17,3)--(18,3);
		\node at (17.5,2.5){\huge{$c^{se}$}};
		\node at (17.5,9.5){\huge{$c^{ne}$}};
		\draw [dashed] (6,9)--(6,10);
		\node at (5.5,9.5){\huge{$c^{nw}$}};
		\node at (11.5,0.75){\Huge{$\hbox{CR}(\h)$}};
		\draw (10,0.75) -- (6.5,0.75);
		\draw[->](6.5,0.75)--(6.5,1.75);

		\draw [color=grigio,fill=grigio] (23,3) rectangle (34,9);
		\draw [color=grigio,fill=grigio] (34,4) rectangle (35,9);
		\draw [color=grigio,fill=grigio] (24,9) rectangle (34,10);
		\draw [color=grigio,fill=grigio] (22,4) rectangle (23,9);
		\draw [color=grigio,fill=grigio] (24,3) rectangle (34,2);
		\draw [dashed] (22,2) rectangle (35,10);
		\draw (21,1) rectangle (36,11);
		\draw (22,4) rectangle (23,5);
		\node at (22.5,4.5) {$\bullet$};
		\draw (23,3) rectangle (24,4);
		\node at (23.5,3.5) {$\bullet$};
		\draw (24,2) rectangle (25,3);
		\node at (24.5,2.5) {$\bullet$};
		\draw (33,2) rectangle (34,3);
		\node at (33.5,2.5) {$\bullet$};
		\draw (34,4) rectangle (35,5);
		\node at (34.5,4.5) {$\bullet$};
		\draw (34,8) rectangle (35,9);
		\node at (34.5,8.5) {$\bullet$};
		\draw (33,9) rectangle (34,10);
		\node at (33.5,9.5) {$\bullet$};
		\draw (24,9) rectangle (25,10);
		\node at (24.5,9.5) {$\bullet$};
		\draw (22,8) rectangle (23,9);
		\node at (22.5,8.5) {$\bullet$};
		\fill[pattern=north west lines] (21,1) rectangle (22,11);
		\fill[pattern=north west lines] (22,1) rectangle (36,2);
		\fill[pattern=north west lines] (35,2) rectangle (36,11);
		\fill[pattern= north west lines] (22,10) rectangle (35,11);
		\node at (30.5,0){\Huge{$\hbox{CR}^+(\h)$}};
		\node at (23.5,0){\Huge{$\partial^+\hbox{CR}(\h)$}};
		\draw (21.5,0) -- (19.5,0);
		\draw (19.5,0)--(19.5,2.5);
		\draw[->] (19.5,2.5)--(20.8,2.5);

	\end{tikzpicture}
	
	\vskip 0 cm
	\caption{Here we depict the same configuration $\h$ on the left and on the right to emphasize different geometrical definitions. The grey area in both pictures represents $C(\h_{cl})$. In particular, on the left-hand side we stress the frame-angles $c^{\a\a'}(\h )$, the bars $B^\a(\h)$, $\hbox{CR}^-(\h)$ and the circumscribing rectangle $\hbox{CR}(\h)$ (respresented with a dashed line). While on the right-hand side we stress the sites that are in a corner (represented with a dot), $\hbox{CR}^+(\h)$ and the external frame $\partial^+\hbox{CR}(\h)$ (the dashed area).}
	\label{fig:figesempio}
\end{figure}

\item[$\bullet$]
A {\it corner} in $\h\in\cX$ is a closed $1\times1$ square centered in an occupied site $x\in\L_0$ such that, if we order clockwise its four nearest neighbors $x_1,x_2,x_3,x_4$, then $\sum_{y\in \hbox{nn}(x)}\h(y)=2$, with $\h(x_i)=\h(x_{i+1})=1$, with $i=1,...,4$ and the convention that $x_5=x_1$ (see Figure \ref{fig:figesempio} on the right-hand side).
\ei

\medskip
\noindent
{\bf 4. Circumscribed rectangle, frames and bars}
\bi

\item[$\bullet$]
If $\h$ is a configuration with a single cluster then we denote by CR$(\h)$ the rectangle {\it circumscribing} $C(\h_{cl})$.

We denote $\partial^+\hbox{CR}(\h)$ the {\it external frame of} $\hbox{CR}(\h)$ as the union of squares $1\times1$ centered at sites that are not contained in $\hbox{CR}(\h)$ such that those sites have Euclidean distance with sites in $\hbox{CR}(\h)$ less or equal than $\sqrt{2}$ (see Figure \ref{fig:figesempio} on the right-hand side). Note that the external frame of $\hbox{CR}(\h)$ contains only non occupied sites.

We denote $\partial^-\hbox{CR}(\h)$ the {\it internal frame of} $\hbox{CR}(\h)$ as the union of squares $1\times1$ centered at sites that are contained in $\hbox{CR}(\h)$ such that those sites have Euclidean distance with sites not in $\hbox{CR}(\h)$ less or equal than $\sqrt{2}$. If this distance is equal to $\sqrt{2}$, we say that the unit square is a {\it frame-angle} $c^{\alpha\alpha'}(\h)$ in $\partial^-\hbox{CR}(\h)$, where $\alpha\alpha'\in\{ne,nw,se,sw\}$, with $n=\hbox{north}$, $s=\hbox{south}$, etc. Note that the internal frame of $\hbox{CR}(\h)$ is a geometrical object contained in $\R^2$ that can contain both occupied and non occupied sites (see Figure \ref{fig:figesempio} on the left-hand side). We partition the set $\partial^-\hbox{CR}(\h)$ without frame-angles in {\it two horizontal and two vertical rows} $r^{\alpha}(\h)$, with $\alpha\in\{n,w,e,s\}$.

Moreover, we set
\be{cr}
\ba{ll}
\hbox{CR}^-(\h)=\hbox{CR}(\h)\setminus\partial^-\hbox{CR}(\h),\\
\hbox{CR}^+(\h)=\hbox{CR}(\h)\cup\partial^+\hbox{CR}(\h).
\ea
\ee

See Figure \ref{fig:figesempio} for an example.

\br{}
Note that, for example, the frame-angles $c^{ne}(\h)$ and $c^{en}(\h)$ are the same, but this distinction will be useful in Definitions \ref{movepart} and \ref{trenino}.
\er

\item[$\bullet$]
A {\it vertical} (respectively {\it horizontal}) {\it bar} $B^{\a}(\h)$ of a single cluster of $\h$ with length $k$ is a $1\times k$ (respectively $k\times1$) rectangle contained in $C(\h_{cl})$, with $\alpha\in\{n,w,e,s\}, \ k\geq1$, such that each square $1\times1$ of the bar is attached only to one square of $C(\h_{cl})\setminus B^\a(\h)$ (see Figure \ref{fig:figesempio} on the left-hand side). In the cases in which it is not specified if the bar is vertical or horizontal we call it simply {\it bar}. If $k=1$ we say that the bar is a {\it protuberance}.

\br{sommabarre}
Note that two bars $B^{\a}(\h)$ and $B^{\a'}(\h)$, with $\a,\a'\in\{n,s,w,e\}$, can possibly intersect in the frame-angle $c^{\a\a'}(\h)$. If this is the case, we get $|B^{\a}(\h)\cup B^{\a'}(\h)|=|B^{\a}(\h)|+|B^{\a'}(\h)|-1$.
\er

\ei

\medskip
\noindent
{\bf 5. Motions along the border}

\noindent
Recall definitions of $|\cdot|$ and $||\cdot||$ (see Section \ref{moddepdef} point 1). In the following, we give the precise notion of translation by 1 of a bar, for example to the left or to the right, while keeping all the squares of the bar attached to the cluster below.

\bd{translation}
Given $\h$ and a bar $B^\a(\h)$ of length $k$, with $\a\in\{n,s,e,w\}$, we say that it is possible to \emph{translate the bar $B^\a(\h)$} if 
\be{condtrasl}
k=|B^\a(\h)|<|\partial^+ B^\a(\h)|.
\ee

\noindent
We define the \emph{$1$-translation of a bar $B^\a(\h)$} of length $k$ as a sequence of configurations $(\h_1,...,\h_k)$ such that $\h_1=\h$ and $\h_i$ is obtained from $\h_{i-1}$  translating by $1$ a unit square along the rectangle $\partial^+ B^\a(\h)\cap C(\h_{cl})$ for any $2\leq i\leq k$.

\ed

\noindent
In Figure \ref{fig:traslazione1} (resp.\ Figure \ref{fig:traslazione2})we depict a $1$-translation of a horizontal (resp.\ vertical) bar at cost $U_1$ (resp.\ $U_2$).

\setlength{\unitlength}{0.96pt}
\begin{figure}
	\centering
	\begin{picture}(400,40)(0,40)
		\thinlines
		\put(20,50){\line(1,0){50}}
		\put(20,80){\line(1,0){20}}
		\put(40,80){\line(0,1){5}}
		\put(40,85){\line(1,0){20}}
		\put(60,85){\line(0,-1){5}}
		\put(60,80){\line(1,0){10}}
		\put(20,50){\line(0,1){30}}
		\put(70,50){\line(0,1){30}}
		\put(75,65){\vector(1,0){35}}
		\begin{footnotesize} \put(85,55){$+U_1$} \end{footnotesize}
		
		\thinlines
		\put(115,50){\line(1,0){50}}
		\put(115,80){\line(1,0){20}}
		\put(135,80){\line(0,1){5}}
		\put(135,85){\line(1,0){15}}
		\put(150,85){\line(0,-1){5}}
		\put(150,80){\line(1,0){5}}
		\put(155,80){\line(0,1){5}}
		\put(155,85){\line(1,0){5}}
		\put(160,85){\line(0,-1){5}}
		\put(160,80){\line(1,0){5}}
		\put(115,50){\line(0,1){30}}
		\put(165,50){\line(0,1){30}}
		\put(170,65){\line(1,0){5}}
		\put(176,65){\line(1,0){5}}
		\put(182,65){\line(1,0){5}}
		\put(188,65){\line(1,0){5}}
		\put(194,65){\line(1,0){5}}
		\put(200,65){\vector(1,0){5}}
		
		\begin{footnotesize}
			\put(175,55){$0,...,0$}
		\end{footnotesize}
		
		\thinlines
		\put(210,50){\line(1,0){50}}
		\put(210,80){\line(1,0){20}}
		\put(230,80){\line(0,1){5}}
		\put(230,85){\line(1,0){5}}
		\put(235,85){\line(0,-1){5}}
		\put(235,80){\line(1,0){5}}
		\put(240,80){\line(0,1){5}}
		\put(240,85){\line(1,0){15}}
		\put(255,85){\line(0,-1){5}}
		\put(255,80){\line(1,0){5}}
		\put(210,50){\line(0,1){30}}
		\put(260,50){\line(0,1){30}}
		\put(265,65){\vector(1,0){35}}
		\begin{footnotesize}
			\put(273,55){$-U_1$}
		\end{footnotesize}
		
		\thinlines
		\put(305,50){\line(1,0){50}}
		\put(305,80){\line(1,0){25}}
		\put(330,80){\line(0,1){5}}
		\put(330,85){\line(1,0){20}}
		\put(350,85){\line(0,-1){5}}
		\put(350,80){\line(1,0){5}}
		\put(305,50){\line(0,1){30}}
		\put(355,50){\line(0,1){30}}
	\end{picture}
	\vskip 0 cm
	\caption{$1$-translation of the horizontal bar $B^n(\h)$ at cost $U_1$.}
	\label{fig:traslazione1}
\end{figure}
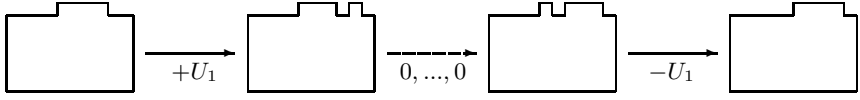

\setlength{\unitlength}{0.95pt}
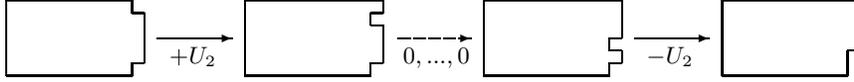
\begin{figure}
	\centering
	\begin{picture}(400,40)(0,40)
		\thinlines
		\put(20,50){\line(1,0){50}}
		\put(20,80){\line(1,0){50}}
		\put(70,75){\line(1,0){5}}
		\put(75,75){\line(0,-1){20}}
		\put(75,55){\line(-1,0){5}}
		\put(60,80){\line(1,0){10}}
		\put(20,50){\line(0,1){30}}
		\put(70,50){\line(0,1){5}}
		\put(70,75){\line(0,1){5}}
		\put(80,65){\vector(1,0){30}}
		\begin{footnotesize} \put(85,55){$+U_2$} \end{footnotesize}
		\thinlines
		\put(115,50){\line(1,0){50}}
		\put(115,80){\line(1,0){55}}
		\put(170,80){\line(0,-1){5}}
		\put(170,75){\line(-1,0){5}}
		\put(165,75){\line(0,-1){5}}
		\put(165,70){\line(1,0){5}}
		\put(170,70){\line(0,-1){15}}
		\put(170,55){\line(-1,0){5}}
		\put(160,80){\line(1,0){5}}
		\put(115,50){\line(0,1){30}}
		\put(165,50){\line(0,1){5}}
		\put(176,65){\line(1,0){5}}
		\put(182,65){\line(1,0){5}}
		\put(188,65){\line(1,0){5}}
		\put(194,65){\line(1,0){5}}
		\put(200,65){\vector(1,0){5}}
		
		\begin{footnotesize}
			\put(178,55){$0,...,0$}
		\end{footnotesize}
		
		\thinlines
		\put(210,50){\line(1,0){50}}
		\put(210,80){\line(1,0){55}}
		\put(265,80){\line(0,-1){15}}
		\put(265,65){\line(-1,0){5}}
		\put(260,65){\line(0,-1){5}}
		\put(260,60){\line(1,0){5}}
		\put(265,60){\line(0,-1){5}}
		\put(265,55){\line(-1,0){5}}
		\put(255,80){\line(1,0){5}}
		\put(210,50){\line(0,1){30}}
		\put(260,50){\line(0,1){5}}
		\put(270,65){\vector(1,0){30}}
		\begin{footnotesize}
			\put(275,55){$-U_2$}
		\end{footnotesize}
		
		\thinlines
		\put(305,50){\line(1,0){50}}
		\put(305,80){\line(1,0){55}}
		\put(360,80){\line(0,-1){20}}
		\put(355,60){\line(1,0){5}}
		\put(305,50){\line(0,1){30}}
		\put(355,50){\line(0,1){10}}
	\end{picture}
	\vskip 0 cm
	\caption{$1$-translation of the vertical bar $B^e(\h)$ at cost $U_2$.}
	\label{fig:traslazione2}
\end{figure}

In the following, we give the precise notion of sliding a unit square from row $r^{\a}(\h)$ to $r^{\a'}(\h)$ passing through the frame angle $c^{\a\a'}(\h)$.

\bd{movepart}
Given $\h$, let $\alpha\alpha'$ such that $c^{\a\a'}(\h)$ is a frame-angle. We say that it is possible to {\rm slide a unit square around a frame-angle $c^{\alpha\alpha'}(\h)\subseteq\partial^-\hbox{CR}(\h)$} from a row $r^{\alpha}(\h)\subseteq\partial^-\hbox{CR}(\h)$ to a row $r^{\alpha'}(\h)\subseteq\partial^-\hbox{CR}(\h)$ via a frame-angle $c^{\alpha\alpha'}(\h)$ if
\be{condcorner}
|c^{\alpha\alpha'}(\h)|=0, \quad |r^{\alpha}(\h)|\geq1, \quad 1\leq|r^{\alpha'}(\h)|<||r^{\alpha'}(\h)||+1.
\ee

\noindent
Let $\a''\neq\a'$ such that $c^{\a\a''}(\h)$ is a frame-angle. See Figure \ref{fig:trenino} for an example. We define a {\rm sliding of a unit square around a frame-angle $c^{\alpha\alpha'}(\h)\subseteq\partial^-\hbox{CR}(\h)$} as the composition of a sequence of $1$-translations of the bar $B^\a(\h)$ from $r^{\alpha}(\h)\cup c^{\a\a''}(\h)$ to $r^{\alpha}(\h)\cup c^{\alpha\alpha'}(\h)$, namely $(\h^1,...,\h^k)$, and the $1$-translation of a bar $B^{\a'}(\h)=C(\h_{cl}^k)\cap (r^{\alpha'}(\h)\cup c^{\alpha\alpha'}(\h))$ from $r^{\alpha'}(\h)\cup c^{\alpha\alpha'}(\h)$ to $r^{\alpha'}(\h)\cup c^{\a'\a'''}(\h)$, where $\a'''\neq\a$ is such that $c^{\a'\a'''}(\h)$ is a frame-angle.
\ed

\setlength{\unitlength}{0.9pt}
\begin{figure}
	\begin{picture}(400,40)(0,40)
		\thinlines
		\put(20,50){\line(1,0){40}}
		\put(20,80){\line(1,0){25}}
		\put(45,80){\line(0,1){5}}
		\put(45,85){\line(1,0){15}}
		\put(60,85){\line(0,-1){5}}
		\put(60,80){\line(1,0){5}}
		\put(65,80){\line(0,-1){20}}
		\put(20,50){\line(0,1){30}}
		\put(60,50){\line(0,1){10}}
		\put(60,60){\line(1,0){5}}
		\put(70,65){\vector(1,0){15}}
		\begin{scriptsize} \put(68,55){$+U_1$} \end{scriptsize}
		
		\thinlines
		\put(90,50){\line(1,0){40}}
		\put(90,80){\line(1,0){25}}
		\put(115,80){\line(0,1){5}}
		\put(115,85){\line(1,0){10}}
		\put(125,85){\line(0,-1){5}}
		\put(125,80){\line(1,0){5}}
		\put(130,80){\line(0,1){5}}
		\put(130,85){\line(1,0){5}}
		\put(135,85){\line(0,-1){25}}
		\put(135,60){\line(-1,0){5}}
		\put(90,50){\line(0,1){30}}
		\put(130,50){\line(0,1){10}}
		\put(140,65){\vector(1,0){15}}
		
		\begin{scriptsize}
			\put(145,55){$0$}
		\end{scriptsize}
		
		\thinlines
		\put(160,50){\line(1,0){40}}
		\put(160,80){\line(1,0){25}}
		\put(185,80){\line(0,1){5}}
		\put(185,85){\line(1,0){5}}
		\put(190,85){\line(0,-1){5}}
		\put(190,80){\line(1,0){5}}
		\put(195,80){\line(0,1){5}}
		\put(195,85){\line(1,0){10}}
		\put(205,85){\line(0,-1){25}}
		\put(205,60){\line(-1,0){5}}
		\put(160,50){\line(0,1){30}}
		\put(200,50){\line(0,1){10}}
		\put(210,65){\vector(1,0){15}}
		\begin{scriptsize}
			\put(208,55){$-U_1$}
		\end{scriptsize}
		
		\thinlines
		\put(230,50){\line(1,0){40}}
		\put(230,80){\line(1,0){30}}
		\put(260,80){\line(0,1){5}}
		\put(260,85){\line(1,0){5}}
		\put(265,85){\line(1,0){10}}
		\put(275,85){\line(0,-1){25}}
		\put(275,60){\line(-1,0){5}}
		\put(230,50){\line(0,1){30}}
		\put(270,50){\line(0,1){10}}
		\put(280,65){\vector(1,0){15}}
		\begin{scriptsize}
			\put(278,55){$+U_2$}
		\end{scriptsize}

		\thinlines
		\put(300,50){\line(1,0){40}}
		\put(300,80){\line(1,0){30}}
		\put(330,80){\line(0,1){5}}
		\put(330,85){\line(1,0){5}}
		\put(335,85){\line(1,0){10}}
		\put(345,85){\line(0,-1){20}}
		\put(345,65){\line(-1,0){5}}
		\put(340,65){\line(0,-1){5}}
		\put(340,60){\line(1,0){5}}
		\put(345,60){\line(0,-1){5}}
		\put(345,55){\line(-1,0){5}}
		\put(300,50){\line(0,1){30}}
		\put(340,50){\line(0,1){5}}
		\put(350,65){\line(1,0){2}}
		\put(353,65){\line(1,0){2}}
		\put(356,65){\line(1,0){2}}
		\put(359,65){\line(1,0){2}}
		\put(362,65){\vector(1,0){3}}
		\begin{scriptsize}
			\put(355,55){$0$}
		\end{scriptsize}
		
		\thinlines
		\put(370,50){\line(1,0){40}}
		\put(370,80){\line(1,0){30}}
		\put(400,80){\line(0,1){5}}
		\put(400,85){\line(1,0){5}}
		\put(405,85){\line(1,0){10}}
		\put(415,85){\line(0,-1){5}}
		\put(415,80){\line(-1,0){5}}
		\put(410,80){\line(0,-1){5}}
		\put(410,75){\line(1,0){5}}
		\put(415,75){\line(0,-1){20}}
		\put(415,55){\line(-1,0){5}}
		\put(370,50){\line(0,1){30}}
		\put(410,50){\line(0,1){5}}
		\put(420,65){\vector(1,0){15}}
		\begin{scriptsize}
			\put(418,55){$-U_2$}
		\end{scriptsize}
		
		\put(440,50){\line(1,0){40}}
		\put(440,80){\line(1,0){30}}
		\put(470,80){\line(0,1){5}}
		\put(470,85){\line(1,0){5}}
		\put(475,85){\line(1,0){5}}
		\put(480,85){\line(0,-1){5}}
		\put(480,80){\line(1,0){5}}
		\put(485,80){\line(0,-1){5}}
		\put(485,75){\line(0,-1){20}}
		\put(485,55){\line(-1,0){5}}
		\put(440,50){\line(0,1){30}}
		\put(480,50){\line(0,1){5}}

	\end{picture}
	\vskip 0 cm
	\caption{Sliding of a unit square around the frame angle $c^{ne}(\h)$ at cost $U_1$. In this case $\a=n$, $\a'=e$, $\a''=w$ and $\a'''=s$.}
	\label{fig:trenino}
\end{figure}
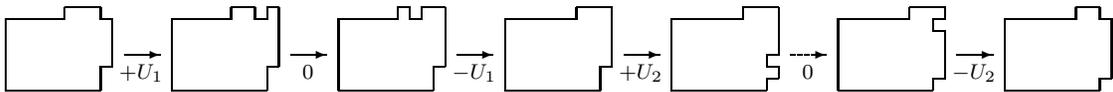

The definition above is used only to define the following sliding of a bar from row $r^{\a}(\h)$ to $r^{\a'}(\h)$ passing through the frame angle $c^{\a\a'}(\h)$, that corresponds to iteratively apply the sliding of a unit square around a frame-angle.

\bd{trenino}
Given $\h$, let $\alpha\alpha'$ such that $c^{\a\a'}(\h)$ is a frame-angle. Before sliding a bar around a frame-angle, we translate the bars $B^{\a}(\h)$ and $B^{\a'}(\h)$ at distance 1 to the frame-angle $c^{\a\a'}(\h)$ obtaining a configuration $\h'$. We say that it is possible to \emph{slide a bar $B^\a(\h')$ around a frame-angle $c^{\alpha\alpha'}(\h')\subseteq\partial^-\hbox{CR}(\h')$} if it is possible to move all the unit squares in $B^\a(\h')$ around a frame-angle $c^{\alpha\alpha'}(\h')$ from a row $r^{\alpha}(\h')\cup c^{\alpha\alpha''}(\h')$ to a row $r^{\alpha'}(\h')\cup c^{\a'\a'''}(\h')$, where $\a''\neq\a'$ and $\a'''\neq\a$ are such that $c^{\a\a''}(\h')$ and $c^{\a'\a'''}(\h')$ are frame-angles. Namely,
\be{condtrenino}
|B^\a(\h')|+|r^{\alpha'}(\h')|\leq||r^{\alpha'}(\h')||+1.
\ee

\noindent
Moreover, we define a \emph{sliding of a bar $B^\a(\h')$ around a frame-angle $c^{\alpha\alpha'}(\h')$} as the sequence of $|B^\a(\h')|$ slidings of unit squares around a frame-angle $c^{\alpha\alpha'}(\h')$.
\ed

\noindent
See the path described in Figure \ref{fig:columntorow}, that connects the configuration $\h$ to the configuration $(12)$ for an example of sliding of the bar $B^e(\h)$ around the frame-angle $c^{en}(\h)$, with $\h$ as the configuration $(3)$.

\subsection{Gate for strongly anisotropic interactions}
\label{sani}

In this Section we impose $U_1>2U_2$ in (\ref{hamilt}), i.e., we consider strongly anisotropic interactions between nearest neighboring sites. Recall the definition of $\e$ given in (\ref{defepsilon}). We will consider $0<\e\ll U_2$, where $\ll$ means sufficiently smaller; for instance $\e\leq\frac{U_2}{100}$ is enough. Many interesting quantities that follow have lower index $sa$ to remind that they refer to strongly anisotropic interactions. In order to state our main results for the gates in the strongly anisotropic regime we need the following definitions. We define the critical vertical length $l_2^*$ as 
\be{critinteri}
l_2^*:=\left\lceil {U_2 \over U_1+U_2 -\D}\right\rceil,
\ee

\noindent
where $\lceil \;\rceil$ denotes the integer part plus 1. We set the critical value of $s_{sa}$ as
\be{defs^*} 
s_{sa}^*:=3l_2^*-1. 
\ee

\noindent
Now we need the following definitions.

\bd{dstrong}
\bi
\item[(a)] We define $\bar\cQ_{sa}$ as the set of configurations having one cluster anywhere in $\L_0$ consisting of a $(2l_2^*-3)\times l_2^*$ rectangle with a single protuberance attached to one of the shortest sides. Similarly, we define $\widetilde\cQ_{sa}$ as the set of configurations having one cluster anywhere in $\L_0$ consisting of a $(2l_2^*-3)\times l_2^*$ rectangle with a single protuberance attached to one of the longest sides. 
\item[(b)] We define
\be{gammastrong}
\gs:=U_1l_2^*+2U_2l_2^*+U_1-U_2-2\e(l_2^*)^2+3\e l_2^*-2\e.
\ee
\item[(c)] We define the volume of the clusters in $\bar\cQ_{sa}$ as
\be{defnsa}
n_{sa}^c:=l_2^*(2l_2^*-3)+1,
\ee

\noindent
and
\be{defDsa}
\ba{lll}
\bar\cD_{sa}:=\{\h'\in\cV_{n_{sa}^c}| \ \exists \ \h\in\bar\cQ_{sa}: H(\h)=H(\h') \hbox{ and } \Phi_{|\cV_{n_{sa}^c}}(\h,\h')\leq H(\h)+U_1\},\\
\widetilde\cD_{sa}:=\{\h'\in\cV_{n_{sa}^c}| \ \exists \ \h\in\widetilde\cQ_{sa}: H(\h)=H(\h') \hbox{ and } \Phi_{|\cV_{n_{sa}^c}}(\h,\h')\leq H(\h)+U_1\}.
\ea
\ee

\noindent
{\rm Note that the last condition in (\ref{defDsa}) is the same as requiring that $\Phi_{|\cV_{n_{sa}^c}}(\h,\h')<\gs+H(\vuoto)=\gs$. We encourage the reader to consult Proposition \ref{card}, where we give the geometrical description of the set $\bar\cD_{sa}$ and $\widetilde\cD_{sa}$. Roughly speaking, one can think of $\bar\cD_{sa}$ as the set of configurations consisting of a rectangular cluster with four bars attached to its four sides, whose lengths satisfy precise conditions.}

\item[(d)] We define
\be{c*sa}
\csgeo:=\bar\cD_{sa}^{fp}.
\ee
\ei
\ed

\noindent
The reason why only the set $\bar\cD_{sa}$ is relevant for the set $\csgeo$ will be clarified later (see Lemma \ref{dtilde}). Note that
\be{}
\ba{lll}
H(\csgeo)&=&H(\bar\cD_{sa}^{fp})=H(\bar\cD_{sa})+\D=H(\cQ_{sa})+\D\\
&=&U_1l_2^*+2U_2l_2^*-U_1-3U_2-\e l_2^*(2l_2^*-3)+2\D\\
&=&U_1l_2^*+2U_2l_2^*+U_1-U_2-2\e(l_2^*)^2+3\e l_2^*-2\e\\
&=&\gs.
\ea
\ee

\noindent
See Figure \ref{fig:P} for an example of configurations in $\csgeo$.

\setlength{\unitlength}{0.89pt}
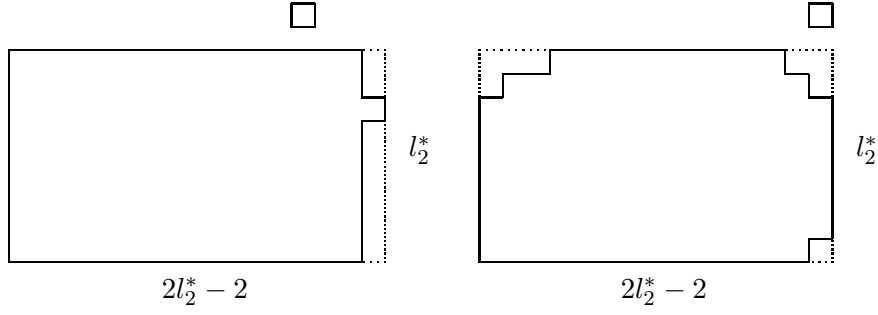
\begin{figure}
	\centering
	\begin{picture}(400,80)(0,30)
		
		\thinlines
		\thinlines
		\qbezier[51](20,0)(100,0)(180,0)
		\qbezier[51](20,90)(100,90)(180,90)
		\qbezier[51](20,0)(20,45)(20,90)
		\qbezier[51](180,0)(180,45)(180,90)
		\thinlines
		\put(20,0){\line(1,0){150}}
		\put(20,90){\line(1,0){150}}
		\put(20,0){\line(0,1){90}}
		\put(170,0){\line(0,1){60}}
		\put(170,70){\line(0,1){20}}
		\put(170,60){\line(1,0){10}}
		\put(180,60){\line(0,1){10}}
		\put(170,70){\line(1,0){10}}
		\put(140,100){\line(1,0){10}}
		\put(150,100){\line(0,1){10}}
		\put(140,100){\line(0,1){10}}
		\put(140,110){\line(1,0){10}}
		\thinlines  
		\put(190,45){$l_2^*$}
		\put(85,-15){$2l_2^*-2$}

		\thinlines
		\qbezier[51](220,0)(300,0)(370,0)
		\qbezier[51](220,90)(300,90)(370,90)
		\qbezier[51](220,0)(220,45)(220,90)
		\qbezier[51](370,0)(370,45)(370,90)
		\thinlines
		\put(220,0){\line(1,0){140}}
		\put(220,0){\line(0,1){70}}
		\put(370,10){\line(0,1){60}}
		\put(370,70){\line(-1,0){10}}
		\put(360,70){\line(0,1){10}}
		\put(360,80){\line(-1,0){10}}
		\put(350,80){\line(0,1){10}}
		\put(350,90){\line(-1,0){100}}
		\put(250,90){\line(0,-1){10}}
		\put(250,80){\line(-1,0){20}}
		\put(230,80){\line(0,-1){10}}
		\put(230,70){\line(-1,0){10}}
		\put(360,10){\line(1,0){10}}
		\put(360,0){\line(0,1){10}}
		\put(360,110){\line(1,0){10}}
		\put(370,110){\line(0,-1){10}}
		\put(370,100){\line(-1,0){10}}
		\put(360,100){\line(0,1){10}}
		
		\put(380,45){$l_2^*$}
		\put(280,-15){$2l_2^*-2$}

	\end{picture}
	\vskip 1.5 cm
	\caption{Critical configurations in $\csgeo$ in the strongly anisotropic case. Moreover, if we remove the free particle we obtain on the left a configuration in $\bar\cQ_{sa}$ and on the right a configuration in $\bar\cD_{sa}\setminus\bar\cQ_{sa}$.}
	\label{fig:P}
\end{figure}

\br{enq}
Note that $H(\bar\cQ_{sa})<H(\widetilde\cQ_{sa})$, indeed

\be{energiaQsa}
\ba{lll}
H(\bar\cQ_{sa})=\gs-\D, \\
H(\widetilde\cQ_{sa})=\gs-\D+U_1-U_2.
\ea
\ee

\er

\noindent
The first main result of Section \ref{sani} is the following.

\bt{thgate}{\rm{(Gate for strongly anisotropic interactions).}}
The set $\csgeo$ is a gate for the transition from $\vuoto$ to $\pieno$.
\et

\noindent
We refer to Section \ref{strong1} for the proof of the Theorem \ref{thgate}.

\br{p2gate}
In Theorem \ref{thgate} we sharpen the previous result obtained in \cite[Theorem 2.4]{BN}. Indeed in \cite{BN} the authors proved that $\cP_1\cup\cP_2$ is a gate (see (\ref{defPsa}) and (\ref{defP2}) for the definitions of $\cP_{sa,0}=\cP_1$ and $\cP_2$ respectively), while here we refine that result by proving that $\csgeo$ is a gate. In particular, we emphasize that $\csgeo\subset\cP_1$ and therefore there is an important improvement in the statement since our gate is much smaller than the one found in \cite{BN}.
\er

In order to give the result regarding the geometric characterization of $\cG_{sa}(\vuoto,\pieno)$, we need some definitions. For any $i=0,1$ we define $\cP_{sa,i}\subseteq\cS_{sa}(\vuoto,\pieno)$ that consists of configurations with a single cluster and no free particle, a fixed number of vacancies, that is not monotone with circumscribed rectangles obtained from the one of the configurations in $\bar\cD_{sa}$ via increasing by one the horizontal or vertical length. More precisely,
\be{defPsa}
\ba{ll}
\cP_{sa,i}:= &\{\h:\, n(\h)= 0,\, v(\h)=2l_2^*+il_2^*-i(i+1)-2,\, g_1'(\h)=i,\,g_2'(\h)=1-i, \h_{cl} \hbox{ is}\\
& \hbox{connected},\,\hbox{with circumscribed rectangle in }
\cR(2l_2^*-i-1,l_2^*+i)\}, \ i=0,1.
\ea
\ee
\noindent 
See Figure \ref{fig:P2} for an example of configurations in $\cP_{sa,0}$ (on the left-hand side) and in $\cP_{sa,1}$ (on the right-hand side).

\setlength{\unitlength}{0.89pt}
\begin{figure}
	\centering
	\begin{picture}(400,90)(0,30)
		\thinlines
		\qbezier[51](20,0)(100,0)(170,0)
		\qbezier[51](20,100)(100,100)(170,100)
		\qbezier[51](20,0)(20,45)(20,100)
		\qbezier[51](170,0)(170,45)(170,90)
		\thinlines
		\put(20,0){\line(1,0){140}}
		\put(20,90){\line(1,0){120}}
		\put(20,0){\line(0,1){90}}
		\put(150,90){\line(1,0){10}}
		\put(160,90){\line(0,1){10}}
		\put(160,100){\line(1,0){10}}
		\put(170,90){\line(0,1){10}}
		\put(170,10){\line(0,1){90}}
		\put(160,10){\line(1,0){10}}
		\put(160,0){\line(0,1){10}}
		\put(140,90){\line(0,1){10}}
		\put(150,90){\line(0,1){10}}
		\put(140,100){\line(1,0){10}}
		\thinlines 
		\put(180,50){$l_2^*$}
		\put(80,-15){$2l_2^*-1$}
		
		\thinlines
		\thinlines
		\qbezier[51](220,110)(290,110)(360,110)
		\qbezier[51](220,0)(290,0)(360,0)
		\qbezier[51](220,0)(220,55)(220,110)
		\qbezier[51](360,0)(360,55)(360,110)
		\thinlines
		\put(230,0){\line(1,0){120}}
		\put(220,10){\line(0,1){100}}
		\put(220,10){\line(1,0){10}}
		\put(230,10){\line(0,-1){10}}
		\put(350,0){\line(0,1){50}}
		\put(350,50){\line(1,0){10}}
		\put(360,50){\line(0,1){10}}
		\put(360,60){\line(-1,0){10}}
		\put(350,60){\line(0,1){10}}
		\put(350,70){\line(1,0){10}}
		\put(360,70){\line(0,1){10}}
		\put(360,80){\line(-1,0){10}}
		\put(350,80){\line(0,1){20}}
		\put(350,100){\line(-1,0){120}}
		\put(230,100){\line(0,1){10}}
		\put(230,110){\line(-1,0){10}}
		\thinlines  
		\put(370,50){$l_2^*+1$}
		\put(275,-15){$2l_2^*-2$}
	\end{picture}
	\vskip 1.5 cm
	\caption{Critical configurations in the strongly anisotropic case: on the left hand-side is represented a configuration in $\cP_{sa,0}$ and on the right hand-side a configuration in $\cP_{sa,1}$.}
	\label{fig:P2}
\end{figure}
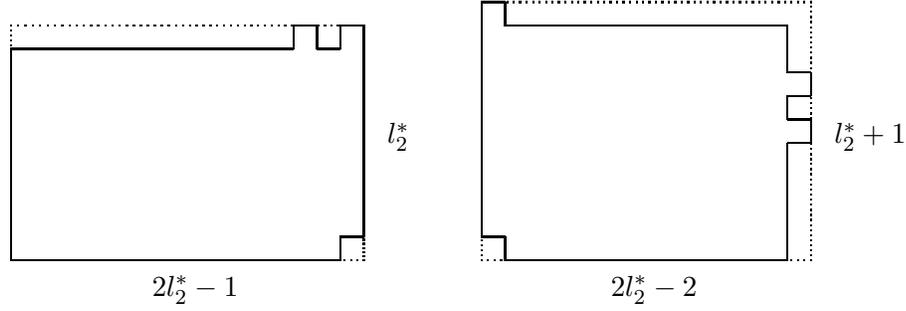

The set $\cG_{sa}(\vuoto,\pieno)$ contains all the configurations that are in the sets defined in (\ref{defPsa}) with the following further conditions. First, we define the subsets $\cA_0^{\a'}$ (resp.\ $\cA_1^{\a}$) of the saddles in $\cP_{sa,0}$ (resp.\ $\cP_{sa,1}$) that contains only one occupied unit square in either a vertical (resp.\ horizontal) row or in one of its two adjacent frame-angles. More precisely,
\be{A1}
\cA_0^{\a'}:=\{\h\in\cP_{sa,0}: |r^{\a'}(\h)\cup c^{\a'\bar\a}(\h)\cup c^{\a'\widetilde\a}(\h)|=1\},
\ee

\noindent
for any $\a'\in\{w,e\}$ and $\bar\a,\widetilde\a\in\{n,s\}$ such that $\bar\a\neq\widetilde\a$, and
\be{A2}
\cA_1^{\a}:=\{\h\in\cP_{sa,1}: |r^{\a}(\h)\cup c^{\a\a''}(\h)\cup c^{\a\a'''}(\h)|=1\},
\ee
for any $\a\in\{n,s\}$ and $\a'',\a'''\in\{n,s\}$ such that $\a'\neq\a'''$. Note that in Figure \ref{fig:P2} the configuration on the right-hand side is in $\cA_1^{n}$.

Next, we define the subsets $\cA_{k}^{\a,\a'}$ of the saddles in $\cP_{sa,0}$ that are obtained from $\h\in\cP_{sa,0}$ during the sliding of the bar $B^{\a'}(\h)$ around the frame-angle $c^{\a'\a}(\h)$. More precisely,
\be{defpsatrenino}
\ba{ll}
\cA_k^{\a,\a'}:=\{\h\in\cP_{sa,0}&:|r^\a(\h)|=k-1, |r^{\a'}(\h)|=l_2^*-k, |c^{\a'\a}(\h)|=1,\\
&(r^{\a}(\h)\cup c^{\a'\a}(\h))\cap\h_{cl}=r^{\a,1}_{cl}\dot{\cup}r^{\a,2}_{cl} \hbox{ with } d(r^{\a,1}_{cl},r^{\a,2}_{cl})=2\},
\ea
\ee

\noindent
where $\a\in\{n,s\}, \ \a'\in\{w,e\}$, $r_{cl}^{\a,1}$, $r_{cl}^{\a,2}$ are two disjoint connected components in $r^{\a(\h)}\cup c^{\a'\a}(\h)$ and $k=2,...,l_2^*$. Note that the conditions in (\ref{defpsatrenino}) guarantee that these configurations are obtained during a sliding of a bar around a frame-angle, that is identified by the indeces $\a$ and $\a'$. Note that in this case there is not the index $k'$ as in \cite[eq.\ (4.4)]{BN2} and \cite[eq.\ (4.11)]{BN2} for the isotropic and weakly anisotropic cases respectively, because in the strongly anisotropic case less sliding on the border of the droplet are allowed. Indeed, in this case $l_2^*-1$ denotes the length of the bar that we are sliding and thus we can consider $k'=l_2^*-1$ fixed. The index $k$ counts the number of particles that are in $r^{\a}(\h)\cup c^{\a'\a}(\h)$ during the sliding and can be less or equal than $l_2^*$. Referring to Figure \ref{fig:columntorow}, configuration (7) (resp.\ (11)) is an example of configuration that belongs to $\cA_2^{n,e}$ (resp.\ $\cA_{l_2^*}^{n,e}$). Note that the set $\cA_k^{\a,\a'}$ contains $k-1$ configurations for any $\a$, $\a'$ and $k$, indeed these configurations are crossed during the sliding of the bar $B^{\a'}(\h)$ around the frame-angle $c^{\a'\a}(\h)$, with $\h$ the configuration obtained by $\cR(2l_2^*-1,l_2^*-1)$ adding a protuberance to one of its longest sides (in Figure \ref{fig:columntorow} $\h$ corresponds to the configuration (3) and  $\a=n$, $\a'=e$). Thus we set
\be{insiemeA}
\cA_{k}^{\a,\a'}=\{\x_1^{\a,\a'},...,\x_{k-1}^{\a,\a'}\}.
\ee

\setlength{\unitlength}{1.08pt}
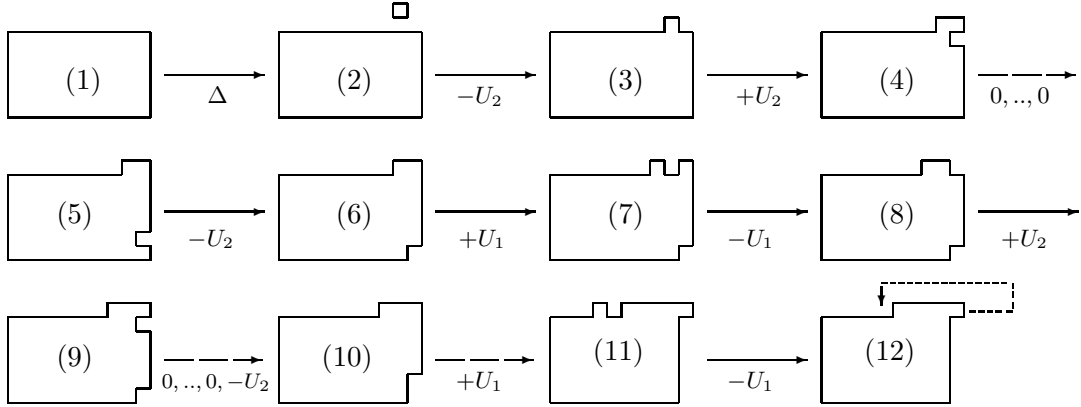
\begin{figure}
	\begin{picture}(400,60)(0,30)
		\thinlines
		\put(20,50){\line(1,0){50}}
		\put(20,80){\line(1,0){50}}
		\put(20,50){\line(0,1){30}}
		\put(70,50){\line(0,1){30}}
		\thinlines  
		\put(40,60){$(1)$}
		\put(75,65){\vector(1,0){35}}
		\begin{footnotesize} \put(90,55){$\D$} \end{footnotesize}
		\thinlines
		\put(115,50){\line(1,0){50}}
		\put(115,80){\line(1,0){50}}
		\put(115,50){\line(0,1){30}}
		\put(165,50){\line(0,1){30}}
		\put(160,85){\line(0,1){5}}
		\put(155,85){\line(1,0){5}}
		\put(155,85){\line(0,1){5}}
		\put(155,90){\line(1,0){5}}
		\thinlines 
		\put(135,60){$(2)$}
		\put(170,65){\vector(1,0){35}}
		\begin{footnotesize}
			\put(177,55){$-U_2$}
		\end{footnotesize}
		
		\thinlines
		\put(210,50){\line(1,0){50}}
		\put(260,50){\line(0,1){30}}
		\put(210,50){\line(0,1){30}}
		\put(210,80){\line(1,0){40}}
		\put(255,80){\line(0,1){5}}
		\put(255,80){\line(1,0){5}}
		\put(250,80){\line(0,1){5}}
		\put(250,85){\line(1,0){5}}
		\thinlines 
		\put(230,60){$(3)$}
		\put(265,65){\vector(1,0){35}}
		\begin{footnotesize}
			\put(275,55){$+U_2$}
		\end{footnotesize}
		
		\thinlines
		\put(305,50){\line(1,0){50}}
		\put(305,50){\line(0,1){30}}
		\put(305,80){\line(1,0){40}}
		\put(345,80){\line(0,1){5}}
		\put(345,85){\line(1,0){10}}
		\put(355,50){\line(0,1){25}}
		\put(350,75){\line(1,0){5}}
		\put(350,75){\line(0,1){5}}
		\put(350,80){\line(1,0){5}}
		\put(355,80){\line(0,1){5}}
		\thinlines 
		\put(325,60){$(4)$}
		\put(360,65){\line(1,0){10}}
		\put(372,65){\line(1,0){10}}
		\put(384,65){\vector(1,0){10}}
		\begin{footnotesize}
			\put(364,55){$0,..,0$}
		\end{footnotesize}

		\thinlines
		\put(20,0){\line(1,0){50}}
		\put(20,0){\line(0,1){30}}
		\put(70,0){\line(0,1){5}}
		\put(20,30){\line(1,0){40}}
		\put(60,30){\line(0,1){5}}
		\put(60,35){\line(1,0){10}}
		\put(65,5){\line(1,0){5}}
		\put(65,5){\line(0,1){5}}
		\put(65,10){\line(1,0){5}}
		\put(70,10){\line(0,1){25}}
		\thinlines 
		\put(37,13){$(5)$}
		\put(75,17){\vector(1,0){35}}
		\begin{footnotesize}
			\put(83,5){$-U_2$}
		\end{footnotesize}

		\thinlines
		\put(115,0){\line(1,0){45}}
		\put(115,0){\line(0,1){30}}
		\put(160,0){\line(0,1){5}}
		\put(115,30){\line(1,0){40}}
		\put(160,5){\line(1,0){5}}
		\put(165,5){\line(0,1){30}}
		\put(155,35){\line(1,0){10}}
		\put(155,30){\line(0,1){5}}
		\thinlines 
		\put(135,13){$(6)$}
		\put(170,17){\vector(1,0){35}}
		\begin{footnotesize}
			\put(178,5){$+U_1$}
		\end{footnotesize}
		
		\thinlines
		\put(210,0){\line(1,0){45}}
		\put(210,0){\line(0,1){30}}
		\put(210,30){\line(1,0){35}}
		\put(245,30){\line(0,1){5}}
		\put(245,35){\line(1,0){5}}
		\put(250,30){\line(0,1){5}}
		\put(250,30){\line(1,0){5}}
		\put(255,30){\line(0,1){5}}
		\put(255,35){\line(1,0){5}}
		\put(255,0){\line(0,1){5}}
		\put(255,5){\line(1,0){5}}
		\put(260,5){\line(0,1){30}}
		\thinlines 
		\put(230,13){$(7)$}
		\put(265,17){\vector(1,0){35}}
		\begin{footnotesize}
			\put(272,5){$-U_1$}
		\end{footnotesize}

		\thinlines
		\put(305,0){\line(1,0){45}}
		\put(305,0){\line(0,1){30}}
		\put(305,30){\line(1,0){35}}
		\put(340,30){\line(0,1){5}}
		\put(340,35){\line(1,0){10}}
		\put(350,30){\line(0,1){5}}
		\put(350,30){\line(1,0){5}}
		\put(350,0){\line(0,1){5}}
		\put(350,5){\line(1,0){5}}
		\put(355,5){\line(0,1){25}}
		\thinlines 
		\put(325,13){$(8)$}
		\put(360,17){\vector(1,0){35}}
		\begin{footnotesize}
			\put(368,5){$+U_2$}
		\end{footnotesize}

		\thinlines
		\put(20,-50){\line(1,0){45}}
		\put(20,-50){\line(0,1){30}}
		\put(20,-20){\line(1,0){35}}
		\put(55,-20){\line(0,1){5}}
		\put(55,-15){\line(1,0){15}}
		\put(65,-50){\line(0,1){5}}
		\put(65,-45){\line(1,0){5}}
		\put(70,-45){\line(0,1){20}}
		\put(65,-25){\line(1,0){5}}
		\put(65,-25){\line(0,1){5}}
		\put(65,-20){\line(1,0){5}}
		\put(70,-20){\line(0,1){5}}
		\thinlines 
		\put(37,-37){$(9)$}
		\put(75,-35){\line(1,0){10}}
		\put(87,-35){\line(1,0){10}}
		\put(99,-35){\vector(1,0){10}}
		\begin{scriptsize}
			\put(74,-45){$0,..,0,-U_2$}
		\end{scriptsize}

		\thinlines
		\put(115,-50){\line(1,0){45}}
		\put(115,-50){\line(0,1){30}}
		\put(115,-20){\line(1,0){35}}
		\put(150,-20){\line(0,1){5}}
		\put(150,-15){\line(1,0){15}}
		\put(160,-50){\line(0,1){10}}
		\put(160,-40){\line(1,0){5}}
		\put(165,-40){\line(0,1){25}}
		\thinlines 
		\put(130,-37){$(10)$}
		\put(170,-35){\line(1,0){10}}
		\put(182,-35){\line(1,0){10}}
		\put(194,-35){\vector(1,0){10}}
		\begin{footnotesize}
			\put(177,-45){$+U_1$}
		\end{footnotesize}

		\thinlines
		\put(210,-50){\line(1,0){45}}
		\put(210,-50){\line(0,1){30}}
		\put(210,-20){\line(1,0){15}}
		\put(225,-20){\line(0,1){5}}
		\put(225,-15){\line(1,0){5}}
		\put(230,-15){\line(0,-1){5}}
		\put(230,-20){\line(1,0){5}}
		\put(235,-20){\line(0,1){5}}
		\put(235,-15){\line(1,0){25}}
		\put(255,-50){\line(0,1){30}}
		\put(255,-20){\line(1,0){5}}
		\put(260,-20){\line(0,1){5}}
		
		\thinlines 
		\put(225,-35){$(11)$}
		\put(265,-35){\vector(1,0){35}}
		\begin{footnotesize}
			\put(272,-45){$-U_1$}
		\end{footnotesize}

		\thinlines
		\put(305,-50){\line(1,0){45}}
		\put(305,-50){\line(0,1){30}}
		\put(305,-20){\line(1,0){25}}
		\put(330,-20){\line(0,1){5}}
		\put(330,-15){\line(1,0){20}}
		\put(350,-50){\line(0,1){30}}
		\put(350,-20){\line(1,0){5}}
		\put(350,-15){\line(1,0){5}}
		\put(355,-20){\line(0,1){5}}
		\put(357,-18){\line(1,0){2}}
		\put(360,-18){\line(1,0){2}}
		\put(363,-18){\line(1,0){2}}
		\put(366,-18){\line(1,0){2}}
		\put(370,-18){\line(1,0){2}}
		\put(372,-17){\line(0,1){2}}
		\put(372,-13){\line(0,1){2}}
		\put(372,-10){\line(0,1){2}}
		\put(370,-8){\line(-1,0){2}}
		\put(367,-8){\line(-1,0){2}}
		\put(364,-8){\line(-1,0){2}}
		\put(361,-8){\line(-1,0){2}}
		\put(358,-8){\line(-1,0){2}}
		\put(355,-8){\line(-1,0){2}}
		\put(352,-8){\line(-1,0){2}}
		\put(349,-8){\line(-1,0){2}}
		\put(346,-8){\line(-1,0){2}}
		\put(343,-8){\line(-1,0){2}}
		\put(340,-8){\line(-1,0){2}}
		\put(337,-8){\line(-1,0){2}}
		\put(334,-8){\line(-1,0){2}}
		\put(331,-8){\line(-1,0){2}}
		\put(328,-8){\line(-1,0){2}}
		\put(326,-9){\line(0,-1){2}}
		\put(326,-11){\vector(0,-1){5}}
		\thinlines 
		\put(320,-35){$(12)$}
		
	\end{picture}
	\vskip 3 cm
	\caption{Transition column to row for $\cR(2l_2^*-1,l_2^*-1)$ in the strongly anisotropic regime: the configuration (1) has energy equal to $\gs-\D+U_2-U_1$ and thus the configurations (7) and (11) have energy equal to $\gs$. In (12) we indicate with a dashed arrow the detachment of the protuberance at cost $U_1$ and afterwards movement of the free particle until it connects to the cluster that decreases the energy by $U_1+U_2$.}
	\label{fig:columntorow}
\end{figure}

\noindent
Now we are able to give the second main result of Section \ref{sani}.

\bt{gstrong}{\rm{(Union of minimal gates for strongly anisotropic interactions).}}
We obtain the following description for $\cG_{sa}(\vuoto,\pieno)$:
\be{}
\cG_{sa}(\vuoto,\pieno)=\csgeo\cup\displaystyle\bigcup_{\a}\bigcup_{\a'}\bigcup_{k=2}^{l_2^*}\cA_k^{\a,\a'}\cup\displaystyle\bigcup_{\a'}\cA_0^{\a'}\cup\displaystyle\bigcup_{\a}\cA_1^{\a}
\ee
\et

\noindent
We refer to Section \ref{strong2} for the proof of Theorem \ref{gstrong}.

\br{simpgates}
With the strategy carried out in \cite{HOS} and the argument explained in Remark \ref{simpmodel}, Theorem \ref{gstrong} can be directly extended to the simplified model.
\er

\subsection{Main results: sharp asymptotics for strongly anisotropic interactions}
\label{sharpestimates}
For a model-independent discussion concerning the prefactor we refer to \cite[Section 10.1]{BN2}. Theorem \ref{sharptimestrong} investigates the prefactor for the strongly anisotropic case. This analysis for the isotropic case is given in \cite[Theorem 1.4.4]{BHN}, while for the weakly anisotropic case is given in \cite[Theorem 4.7]{BN2}. For the proof of Theorems \ref{sharptimestrong} and \ref{gapspettrale} we refer to Section \ref{sharpasymptotics}.

\bt{sharptimestrong}
There exists a constant $K_{sa}=K_{sa}(\L,l_2^*)$ such that
\be{nuclsharpstrong}
\E_{\,\square}(\tau_\blacksquare) = K_{sa}e^{\Gamma_{sa}^*\beta}\,[1+o(1)],
\qquad \b\to\infty,
\ee

\noindent
with
\be{migliore}
\frac{1}{\Theta^{sa}_2}\leq K_{sa}\leq \frac{1}{\Theta^{sa}_1}, 
\ee

\noindent
where $\Theta_{1}^{sa}$ and $\Theta^{sa}_2$ are defined in (\ref{lowerbound}) and (\ref{upperbound}) respectively. Moreover, as $\L\to\Z^2$,
\be{Casymp}
K_{sa}(\L,l_2^*) \ra \frac{1}{4\pi N_{sa}}\,\frac{\log|\L|}{|\L|}
\ee
with
\be{Nid}
N_{sa}=\displaystyle\sum_{k=1}^4\binom{4}{k}\binom{l_2^*+k-2}{2k-1}
\ee
the cardinality of $\bar\cD_{sa}=\bar\cD_{sa}(\L,l_2^*)$ modulo shifts.
\et

\br{}
For the strongly anisotropic case we obtain a sharp estimate of $K_{sa}$ in (\ref{migliore}). Nevertheless, the asymptotic behavior of the the prefactor $K_{sa}$ as $\L\ra\Z^2$ (see (\ref{Casymp})) is the same as $K_{wa}$ (see \cite[eq.\ (4.14)]{BN2}) and $K_{is}$ (see \cite[eq.\ (1.4.9)]{BHN}).
\er

\br{}
Concerning the asymptotics for the prefactor $K_{sa}$ given in Theorem \ref{sharptimestrong}, the sequential limits $\b\ra\infty$ and $\L\ra\Z^2$ are not the physical relevant ones, but since for Kawasaki dynamics we are not able to obtain an explicit expression for the prefactor (see (\ref{migliore})), we give its asymptotic behavior for $\L\ra\Z^2$. The more interesting limit is the joint one $\L=\L^\b\ra\Z^2$ and $\b\ra\infty$, but this is a much harder problem that is out of the scope of the present paper.
\er

\br{}
Note that \cite[Theorem 1.4.3(iii)]{BHN} and \cite[Theorem 4.8]{BN2} concerning the uniform entrance distribution in the gate does not hold for the strongly anisotropic case due to the two possible entrance mechanisms in $\csgeo$ (see Lemma \ref{entratastrong}).
\er

\noindent
We define the {\it mixing time} as
\be{}
t_{mix}(\e):=\min\{n\geq0: \max_{x\in\cX}||P^n (x,\cdot)-\mu(\cdot)||_{TV}\leq\e\},
\ee

\noindent
where $||\nu-\nu'||_{TV}:=\frac{1}{2}\sum_{x\in\cX}|\nu(x)-\nu'(x)|$ for any two probability distributions $\nu,\nu'$ on $\cX$. The {\it spectral gap} of the Markov chain is defined as
\be{}
\r:=\-a^{(2)}
\ee

\noindent
where $1=a^{(1)}>a^{(2)}\geq...\geq a^{|\cX|}\geq-1$ are the eigenvalues of the matrix $(P(x,y))_{x,y\in\cX}$ defined in (\ref{defkaw}).

\bt{gapspettrale}
For any $\e\in(0,1)$
\be{}
\displaystyle\lim_{\b\ra\infty}\frac{1}{\b}\log t_{mix}(\e)=\G_{sa}^*=\lim_{\b\ra\infty}-\frac{1}{\b}\log \r
\ee

\noindent
Furthermore, there exist two constants $0<c_1\leq c_2<\infty$ independent of $\b$ such that for every $\b>0$
\be{}
c_1e^{-\b\G_{sa}^*}\leq\r\leq c_2e^{-\b\G_{sa}^*}
\ee
\et

\noindent
Theorem \ref{gapspettrale} holds also for the isotropic and weakly anisotropic cases (see \cite[Theorem 4.10]{BN2}).

\section{Useful model-dependent definitions and tools}
\label{dependentdef}
In this Section (and only here) we set 
\be{definizioneqbarsa}
\cQ_{sa}=\bar\cQ_{sa}, \qquad \cD_{sa}=\bar\cD_{sa},
\ee
\noindent
to show the similarites between the results with the isotropic model. Moreover, since we are considering the strongly anisotropic model and some properties are in common with the isotropic and weakly anisotropic models, we choose the lower index $int\in\{is,wa,sa\}$ to make clear in the notation which of the three models we are referring to.

\subsection{Geometric description of the protocritical droplets}
In \cite[Theorem 1.4.1]{BHN} the authors obtain the geometric description of the set $\cD_{is}$ as $\cD_{is}=\bar\cD_{is}\cup\widetilde\cD_{is}$. In this Section we derive the geometric description of the analogous sets for the strongly anistropic models $\bar\cD_{sa}$ and $\widetilde\cD_{sa}$ following the argument proposed in \cite{BHN}. The geometric description of the sets $\bar\cD_{wa}$ and $\widetilde\cD_{wa}$ is given in \cite[Proposition 7.1]{BN2}. Recall definition (\ref{defDsa}).

\bp{card}{\rm{(Geometric description of $\widetilde{\cD}_{sa}$ and $\bar\cD_{sa}$).}}
We obtain the following geometric description of $\widetilde\cD_{sa}$ and $\bar\cD_{sa}$:
\bi
\item[(a)]
$\widetilde\cD_{sa}=\{\h\in\cX :\, n(\h)= 0,\, v(\h)=2l_2^*-4,\,
\h_{cl} \hbox{ is connected and monotone},\\
\hbox{   \qquad \quad     } 1\leq|r^\a(\h)\cup c^{\a\a'}(\h)|\leq2, |r^\a(\h)|\leq1, \hbox{ with } \a\in\{n,s\},\ \a'\in\{w,e\}, \hbox{ and }\\
\hbox{   \qquad \quad     } \hbox{ circumscribed rectangle in } \cR(2l_2^*-3,l_2^*+1)\}$,
\item[(b)] $\bar\cD_{sa}$ is the set of configurations having one cluster $\h$ anywhere in $\L_0$ consisting of a $(2l_2^*-4)\times(l_2^*-2)$ rectangle with four bars $B^\a(\h)$, with $\a\in\{n,w,e,s\}$, attached to its four sides satisfying 
\be{condbard}
1\leq|B^w(\h)|,|B^e(\h)|\leq l_2^*, \qquad l_2^*-1\leq |B^n(\h)|,|B^s(\h)|\leq 2l_2^*-2,
\ee

\noindent
and
\be{condbarstrong}
\displaystyle\sum_\a|B^{\a}(\h)|-\displaystyle\sum_{\a\a'\in\{nw,ne,sw,se\}}|c^{\a\a'}(\h)|=5l_2^*-7.
\ee

\ei
\ep

\br{barrestrong}
Let $\h\in\bar\cD_{sa}$. 

\bi
\item[(i)] Note that (\ref{condbarstrong}) takes into account the number of occupied unit squares in $\partial^-\hbox{CR}(\h)$ due to Remark \ref{sommabarre}. We deduce that at most three frame-angles of CR$(\h)$ can be occupied, otherwise $|\partial^-\hbox{CR}(\h)|=6l_2^*-8>5l_2^*-7$, which is absurd. 

\item[(ii)] Since $|B^{s}(\h)|+|B^{e}(\h)|\leq3l_2^*-6+k-|c^{nw}(\h)|$, we get
\be{} 
\ba{ll}
|B^{n}(\h)|+|B^{w}(\h)|=5l_2^*-7-(|B^{s}(\h)|+|B^{e}(\h)|)+k\geq2l_2^*-1+|c^{nw}(\h)|.
\ea
\ee

\noindent
By symmetry, we generalize the inequality above for any $\a\in\{n,s\}$ and $\a'\in\{w,e\}$: we get $|B^{\a}(\h)|+|B^{\a'}(\h)|\geq2l_2^*-1+|c^{\a\a'}(\h)|$.
\ei
\er

\begin{proof*}{\bf of Proposition \ref{card}}
	(a) We introduce only in the proof of this result the following geometrical definition to make the argument more clear
	\be{}
	\ba{ll}
	\widetilde\cD_{sa}^{geo}:=\{\h\in\cX :\, n(\h)= 0,\, v(\h)=2l_2^*-4,\,
	\h_{cl} \hbox{ is connected and monotone},\\
	\hbox{   \qquad \quad     } 1\leq|r^\a(\h)\cup c^{\a\a'}(\h)|\leq2, |r^\a(\h)|\leq1, \hbox{ with } \a\in\{n,s\},\ \a'\in\{w,e\}, \hbox{ and }\\
	\hbox{   \qquad \quad     } \hbox{ circumscribed rectangle in } \cR(2l_2^*-3,l_2^*+1)\}.
	\ea
	\ee
	
	\noindent
	The proof will be given in two steps:
	
	\bi
	\item[(i)] $\widetilde\cD_{sa}^{geo}\subseteq\widetilde\cD_{sa}$;
	\item [(ii)] $\widetilde\cD_{sa}^{geo}\supseteq\widetilde\cD_{sa}$.
	\ei
	
	\noindent
	{\bf Proof of (i).} To prove (i) we must show that for all $\h\in\widetilde\cD_{sa}^{geo}$,
	\bi
	\item[(i1)] $H(\h)=H(\widetilde\cQ_{sa})$;
	\item[(i2)] $\exists \ \o:\widetilde\cQ_{sa}\ra\h$, i.e., $\o=(\o_1,...,\o_k=\h)$ such that $\displaystyle\max_i H(\o_i)\leq H(\widetilde\cQ_{sa})+U_1$, with $|\o_i|=n_{sa}^c$ for all $i=1,...,k$ (see (\ref{defnsa}) for the definition of $n_{sa}^c$).
	\ei

	{\bf Proof of (i1).} Any $\h\in\widetilde\cD_{sa}^{geo}$ satisfies $n(\h)=0$, $|C(\h)|=(2l_2^*-3)(l_2^*+1)-v(\h)=n_{sa}^c$, and $g_1(\h)=2l_2^*-3$ and $g_2(\h)=l_2^*+1$ since the configuration is monotone. Thus by (\ref{Hcont}) we deduce that $H$ is constant on $\widetilde\cD_{sa}^{geo}$. Since $\widetilde\cQ_{sa}\subseteq\widetilde\cD_{sa}^{geo}$, this completes the proof of (i1).

	\medskip
	{\bf Proof of (i2).} Consider $\z\in\widetilde\cQ_{sa}$ and $\h\in\widetilde\cD_{sa}^{geo}$. If $\h\in\widetilde\cQ_{sa}\cap\widetilde\cD_{sa}^{geo}$, i.e., $|r^{\a}(\h)\cup c^{\a\a'}(\h)|=1$ for some $\a\in\{n,s\}$ and $\a'\in\{w,e\}$, then $\h$ can be obtained from $\z$ by moving the protuberance at zero cost along the side which is attached to if the protuberance in $\z$ is on the same side as the protuberance in $\h$, otherwise $\h$ can be obtained detaching the protuberance at cost $U_2$ and reattaching it to the other side at cost $-U_2$. If $\h\in\widetilde\cD_{sa}^{geo}\setminus\widetilde\cQ_{sa}$, i.e., $|r^{\a}(\h)\cup c^{\a\a'}(\h)|=2$ with $|r^{\a}(\h)|=1$ for some $\a\in\{n,s\}$ and $\a'\in\{w,e\}$, again we have two cases. If the protuberance in $\z$ is contained in $r^{\a}(\h)$ (is in the same side of the rectangle), we deduce that $\h$ can be obtained from $\z$ by moving the protuberance at zero cost until it arrives at distance one from $c^{\a\a'}(\z)$ and then translate the bar $B^{\a'}(\z)$ towards the frame-angle $c^{\a'\a}(\z)$ at cost $U_2$. Otherwise, if the protuberance in $\z$ is contained in $r^{\a''}(\h)$ with $\a''\in\{n,s\}\setminus\{\a\}$ (is in the opposite side of the rectangle), the path is constructed as before, provided that first the protuberance is detached at cost $U_2$ and reattached to the other side at cost $-U_2$. This concludes the proof of (i2).

	\medskip
	\noindent
	{\bf Proof of (ii).} By (i2), we know that all the configurations in $\widetilde\cD_{sa}^{geo}$ are connected via $U_1$-path to $\widetilde\cQ_{sa}$. Since $\widetilde\cQ_{sa}\subseteq\widetilde\cD_{sa}\cap\widetilde\cD_{sa}^{geo}$, in order to prove (ii) it suffices to show that following $U_1$-paths it is not possible to exit $\widetilde\cD_{sa}^{geo}$. Let $\h'\in\widetilde\cD_{sa}^{geo}$, thus by (i1) and (\ref{energiaQsa}) we get $H(\h')=\gs-\D+U_1-U_2$. Note that no particle can arrive because we impose that the number of particles is constant to $n_{sa}^c$ thus, if $\h'\in\widetilde\cD_{sa}^{geo}\setminus\widetilde\cQ_{sa}$, the unique possibility is to translate the bar $B^{\a'}(\h')$, with $\a'\in\{w,e\}$, at cost $U_2$ giving rise to a configuration that is in $\widetilde\cD_{sa}^{geo}\cap\widetilde\cQ_{sa}$. Then it is possible either yo move the protuberance at zero cost, or to detach the protuberance at cost $U_2$ and then necessarily reattach it at cost $-U_2$, giving rise to a configuration that is still in $\widetilde\cD_{sa}^{geo}\cap\widetilde\cQ_{sa}$. Note that no other moves are allowed, since  it is not possible to complete the sliding of a vertical bar around a frame-angle because the $1$-translation of a horizontal bar costs $U_1$. Indeed the latter implies that the energy reaches the value $\gs-\D+2U_1-U_2>\gs$, thus the path described is not a $U_1$-path. If $\h'\in\widetilde\cD_{sa}^{geo}\cap\widetilde\cQ_{sa}$, with paths similar to the ones described above (possibly in different order), we can conclude case (a).

	(b) We denote by $\bar\cD_{sa}^{geo}$ the geometric set with the properties specified in point (b) that we introduce to make the argument more clear. The proof will be given in two steps:
	
	\bi
	\item[(i)] $\bar\cD_{sa}^{geo}\subseteq\bar\cD_{sa}$;
	\item [(ii)] $\bar\cD_{sa}^{geo}\supseteq\bar\cD_{sa}$.
	\ei
	
	\noindent
	{\bf Proof of (i).} To prove (i) we must show that for all $\h\in\bar\cD_{sa}^{geo}$,
	\bi
	\item[(i1)] $H(\h)=H(\bar\cQ_{sa})$;
	\item[(i2)] $\exists \ \o:\bar\cQ_{sa}\ra\h$, i.e., $\o=(\o_1,...,\o_k=\h)$ such that $\displaystyle\max_i H(\o_i)\leq H(\bar\cQ_{sa})+U_1$, with $|\o_i|=n_{sa}^c$ for all $i=1,...,k$ and $\o_1\in\bar\cQ_{sa}$.
	\ei

	{\bf Proof of (i1).} Any $\h\in\bar\cD_{sa}^{geo}$ satisfies $n(\h)=0$, $|C(\h)|=(2l_2^*-2)(l_2^*-2)+5l_2^*-7=n_{sa}^c$, and $g_1(\h)=2l_2^*-2$ and $g_2(\h)=l_2^*$ since the configuration is monotone. Thus by (\ref{Hcont}) we deduce that $H$ is constant on $\bar\cD_{sa}^{geo}$. Since $\bar\cQ_{sa}\subseteq\bar\cD_{sa}^{geo}$, this completes the proof of (i1).
	
	\medskip
	{\bf Proof of (i2).} Consider $\z\in\bar\cQ_{sa}$ and $\h\in\bar\cD_{sa}^{geo}$. Here, without loss of generality, we assume that the protuberance is in $r^w(\z)$. 
	Then we have
	\bi
	\item[-] $|B^w(\z)|=1$;
	\item[-] $|B^n(\z)|=|B^s(\z)|=2l_2^*-3$;
	\item[-] $|B^e(\z)|=l_2^*$;
	\item[-] $|c^{ne}(\z)|=|c^{se}(\z)|=1$.
	\ei
	
	\noindent
	Using the sliding of a unit square around a frame-angle described in Figure \ref{fig:trenino} (see Definition \ref{movepart}), we move, one by one, $|B^n(\z)|-|B^n(\h)|$ particles around the frame-angle $c^{nw}(\z)$. After that we move $|B^e(\z)|-|B^e(\h)|+|B^s(\z)|-|B^s(\h)|$ particles around the frame-angle $c^{sw}(\z)$. Finally, we move $|B^e(\z)|-|B^e(\h)|$ particles around the frame-angle $c^{es}(\z)$. The result is the configuration $\h\in\bar\cD_{sa}^{geo}$. This concludes the proof of (i2).
	
	\medskip
	\noindent
	{\bf Proof of (ii).} By (i2), we know that all configurations in $\bar\cD_{sa}^{geo}$ are connected via $U_1$-path to $\bar\cQ_{sa}$. Since $\bar\cQ_{sa}\subseteq\bar\cD_{sa}\cap\bar\cD_{sa}^{geo}$, in order to prove (ii) it suffices to show that following $U_1$-paths it is not possible to exit $\bar\cD_{sa}^{geo}$. We call a path {\it clustering} if all the configurations in the path consist of a single cluster and no free particles. Below we will prove that for any $\h\in\bar\cD_{sa}^{geo}$ and any $\h'$ connected to $\h$ by a clustering $U_1$-path, the following conditions hold:
	
	\bi
	\item[(A)] $\hbox{CR}(\h')=\hbox{CR}(\h)$;
	\item[(B)] $\h'\supseteq\hbox{CR}^-(\h)$.
	\ei

	\medskip
	{\bf Proof of (A).} Starting from any $\h\in\cX$, it is geometrically impossible to modify $\hbox{CR}(\h)$ without detaching a particle, that contradicts the hypotheses of clustering $U_1$-path.
	
	\medskip
	{\bf Proof of (B).} Fix $\h\in\bar\cD_{sa}^{geo}$. The proof is done in two steps.
	
	\medskip
	\noindent
	{\bf 1.} First, we consider clustering $U_1$-paths along which we do not move a particle from $\hbox{CR}^-(\h)$. Along such paths we only encounter configurations in $\bar\cD_{sa}^{geo}$ or configurations obtained from $\bar\cD_{sa}^{geo}$ by breaking one of the bars in $\partial^-\hbox{CR}(\h)$ into two pieces at cost $U_1$ (resp.\ $U_2$) if the bar is horizontal (resp.\ vertical). This holds because there is no particle outside $\hbox{CR}(\h)$ that can lower the cost. 
	
	If the broken bar is horizontal, then only moves at zero cost are admissible, so any particle can be detached. This implies that the unique way to regain $U_1$ and complete the $U_1$-path is to restore the bar.
	
	If the broken bar is vertical, then the admissible moves in a $U_1$-path are those with cost less or equal than $U_1-U_2$. Again any particle can not be detached, indeed its cost is at least $U_1$. The moves at cost $U_2$ are possible, thus it is possible to break another vertical bar. From now on, depending on  $U_1-2U_2>0$, it is possible to break other vertical bars. More precisely, let $U_1=nU_2+\d$, with $n\geq2$ and $0<\d<U_2$ fixed, thus it is possible to break other $n-2$ vertical bars in addition to the previous two. When this sequence of moves is completed, the unique way to complete the $U_1$-path is to restore all the broken bars. Thus we have proved that $\h'\supseteq\hbox{CR}^-(\h)$.

	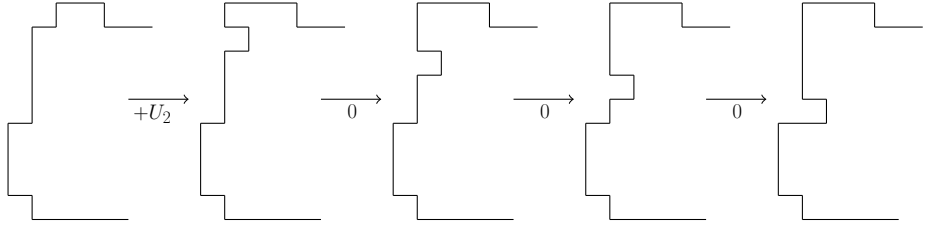
\begin{figure}[htp]
		\centering
		\begin{tikzpicture}[scale=0.32,transform shape]

			\draw (0,0)--(4,0);
			\draw (0,0)--(0,1);
			\draw (0,1)--(-1,1);
			\draw (-1,1)--(-1,4);
			\draw(-1,4)--(0,4);
			\draw (0,4)--(0,8);
			\draw (0,8)--(1,8);
			\draw(1,8)--(1,9);
			\draw (1,9)--(3,9);
			\draw (3,9)--(3,8);
			\draw (3,8)--(5,8);
			\draw[->] (4,5)--(6.5,5);
			\node at (5,4.4) {{\Huge{$+U_2$}}};

			\draw (8,0)--(12,0);
			\draw (8,0)--(8,1);
			\draw (8,1)--(7,1);
			\draw (7,1)--(7,4);
			\draw(7,4)--(8,4);
			\draw (8,4)--(8,7);
			\draw (8,7)--(9,7);
			\draw (9,7)--(9,8);
			\draw (9,8)--(8,8);
			\draw(8,8)--(8,9);
			\draw (8,9)--(11,9);
			\draw (11,9)--(11,8);
			\draw (11,8)--(13,8);
			\draw[->] (12,5)--(14.5,5);
			\node at (13.3,4.5) {{\Huge{0}}};

			\draw (16,0)--(20,0);
			\draw (16,0)--(16,1);
			\draw (16,1)--(15,1);
			\draw (15,1)--(15,4);
			\draw(15,4)--(16,4);
			\draw (16,4)--(16,6);
			\draw (16,6)--(17,6);
			\draw (17,6)--(17,7);
			\draw (17,7)--(16,7);
			\draw(16,7)--(16,9);
			\draw (16,9)--(19,9);
			\draw (19,9)--(19,8);
			\draw (19,8)--(21,8);
			\draw[->] (20,5)--(22.5,5);
			\node at (21.3,4.5) {{\Huge{0}}};

			\draw (24,0)--(28,0);
			\draw (24,0)--(24,1);
			\draw (24,1)--(23,1);
			\draw (23,1)--(23,4);
			\draw(23,4)--(24,4);
			\draw (24,4)--(24,5);
			\draw (24,5)--(25,5);
			\draw (25,5)--(25,6);
			\draw (25,6)--(24,6);
			\draw(24,6)--(24,9);
			\draw (24,9)--(27,9);
			\draw (27,9)--(27,8);
			\draw (27,8)--(29,8);
			\draw[->] (28,5)--(30.5,5);
			\node at (29.3,4.5) {{\Huge{0}}};

			\draw (32,0)--(36,0);
			\draw (32,0)--(32,1);
			\draw (32,1)--(31,1);
			\draw (31,1)--(31,4);
			\draw(31,4)--(33,4);
			\draw (33,4)--(33,5);
			\draw (33,5)--(32,5);
			\draw (32,5)--(32,9);
			\draw (32,9)--(35,9);
			\draw (35,9)--(35,8);
			\draw (35,8)--(37,8);

		\end{tikzpicture}
		
		\vskip 0 cm
		\caption{Creation and motion of the recess at cost 0.}
		\label{fig:recess}
	\end{figure}
	
	\noindent
	{\bf 2.}  Consider now a general clustering $U_1$-path along which we move a particle from a corner of $\hbox{CR}^-(\h)$. It is not allowed to move at cost $U_1+U_2$, because it exceeds $U_1$, thus the overshoot $U_2$ must be regained by letting the particle slide next to a bar that is attached to a side of $\hbox{CR}^-(\h)$ (see Figure \ref{fig:recess}). If the particle moves vertically (resp.\ horizontally), we regain $U_1$ (resp.\ $U_2$). Since there are never two bars attached to the same side, we can at most regain $U_1$, thus it is not possible to move a particle from $\hbox{CR}^-(\h)$ other than from a corner. If the corner particle has been moved vertically (increasing the energy by $U_2$), the same moves (if possible) are allowed on another corner. Depending on the difference $U_1-2U_2>0$, it is possible to break some vertical bars. More precisely, let $U_1=nU_2+\d$, with $n\geq2$ and $0<\d<U_2$ fixed, it is possible to break $n-2$ vertical bars. From now on, only moves at cost at most zero are admissible. There are no protuberances present anymore, because only the configurations in $\bar\cQ_{sa}$ have a protuberance. Thus no particle outside $\hbox{CR}^-(\h)$ can move, except those just detached from $\hbox{CR}^-(\h)$. These particles can move back, in which case we return to the same configuration $\h$ (see Figure \ref{fig:recess}). In fact, all possible moves at zero cost consist in moving the recess just created in $\hbox{CR}^-(\h)$ along the same side of $\hbox{CR}^-(\h)$, until it reaches the top of the bar, after which it cannot advance anymore at zero cost (see Figure \ref{fig:recess}). All these moves do not change the energy, except the last one that returns the particle to its original position and regains $U_1$. This concludes the proof of (B).
	
	From (A), we deduce that $\hbox{CR}(\h')=\cR(2l_2^*-2,l_2^*)$. From (A) and (B), we deduce that the number of particles that are in $\partial^-\hbox{CR}(\h)$ is equal to the number of particles that are in $\partial^-\hbox{CR}(\h')$, thus (\ref{condbarstrong}), $1\leq|B^w(\h')|,|B^e(\h')|\leq l_2^*$ and $1\leq |B^n(\h')|,|B^s(\h')|\leq 2l_2^*-2$ hold. In order to prove that following clustering $U_1$-paths it is not possible to exit $\bar\cD_{sa}^{geo}$, we have to prove the lower bound in (\ref{condbard}) for the lengths $|B^n(\h')|$ and $|B^s(\h')|$. We set
	\be{} 
	k=\displaystyle\sum_{\a\a'\in\{nw,ne,sw,se\}}|c^{\a\a'}(\h')|. 
	\ee
	
	\noindent
	Since $|B^w(\h')|+|B^e(\h')|\leq2l_2^*-4+k$, by (\ref{condbarstrong}) we get
	\be{formulabarre}
	|B^n(\h')|+|B^s(\h')|=5l_2^*-7-(|B^w(\h')|+|B^e(\h')|)+k\geq3l_2^*-3.
	\ee
	
	\noindent
	Since $|B^s(\h')|\leq2l_2^*-2$, (\ref{formulabarre}) implies
	\be{}
	|B^n(\h')|\geq3l_2^*-3-|B^s(\h')|\geq l_2^*-1.
	\ee
	
	\noindent
	By symmetry we can similarly argue for the length $|B^s(\h')|$. This implies that following $U_1$-paths it is not possible to exit $\bar\cD_{sa}^{geo}$ . The argument goes as follows. Detaching a particle costs at least $U_1+U_2$ unless the particle is a protuberance, in which case the cost is $U_1$. The only configurations in $\bar\cD_{sa}^{geo}$ having a protuberance are those in $\bar\cQ_{sa}$. If we detach the protuberance from a configuration in $\bar\cQ_{sa}$, then we obtain a $(2l_2^*-3)\times l_2^*$ rectangle with a free particle. Since in the sequel only moves at zero cost are allowed, it is only possible to move the free particle. Since in a $U_1$-path the particle number is conserved, the only way to regain $U_1$ and complete the $U_1$-path is to reattach the free particle to a vertical side of the rectangle, thus return to $\bar\cQ_{sa}$. This implies that for any $\h\in\bar\cD_{sa}^{geo}$ and any $\h'$ connected to $\h$ by a $U_1$-path we must have that $\h'\in\bar\cD_{sa}^{geo}$. This concludes the proof.
\end{proof*}

\subsection{Definitions}
\label{sitigood}

We set
\be{L*}
L_{sa}^*:=	L-l_2^*.
\ee

\noindent
For $\h\in\cC_{sa}^*$, we associate $(\hat\h,x)$ with $\hat\h\in\cD_{sa}$ protocritical droplet and $x\in\L$ the position of the free particle. We denote by $\cC_{sa}^G(\hat\h)$ (resp.\ $\cC_{sa}^B(\hat\h)$) the configurations that can be reached from $(\hat\h,x)$ by a path that moves the free particle towards the cluster and attaches the particle in $\partial^-CR(\hat\h)$ (resp.\ $\partial^+CR(\hat\h)$). In Figure \ref{fig:sitiGB} on the left-hand side we depict explicitly the good and bad sites for a specific $\hat\h$. Let 
\be{sitiGB}
\cC_{sa}^G=\displaystyle\bigcup_{\hat\h\in\cD_{sa}}\cC_{sa}^G(\hat\h), \qquad \cC_{sa}^B=\displaystyle\bigcup_{\hat\h\in\cD_{sa}}\cC_{sa}^B(\hat\h).
\ee

\noindent
For $\h\in\cC_{sa}^*$, let $\hat\h\in\cD_{sa}$ be the configuration obtained from $\h$ by removing the free particle. For $A\subseteq\L$ and $x\in\L$, recall that $d(x,A)$ denotes the lattice distance between $x$ and $A$. As in \cite[Section 3.5]{BHN}, we need the following definitions.

\bd{cornis}
Let $\L_4$ be $\L$ without its four frame-angles. We define, recursively, 
\be{}
B_1(\hat\h):=\{x\in\L_4| \ x\notin\hat\h, \ d(x,\hat\h)=1\}
\ee
and
\be{} 
\ba{lll} 
B_2(\hat\h):=\{x\in\L_4| \ x\notin\hat\h, \ d(x,B_1(\hat\h))=1\},  \\
\bar B_2(\hat\h):=B_2(\hat\h),
\ea
\ee
and
\be{}
\ba{lll}
B_3(\hat\h):=\{x\in\L_4| \ x\notin B_1(\hat\h), \ d(x,B_2(\hat\h))=1\},  \\
\bar B_3(\hat\h):=B_3(\hat\h)\cup \{\bar B_2(\hat\h)\cap \partial^-\L_4\},
\ea
\ee

\setlength{\unitlength}{0.16cm}

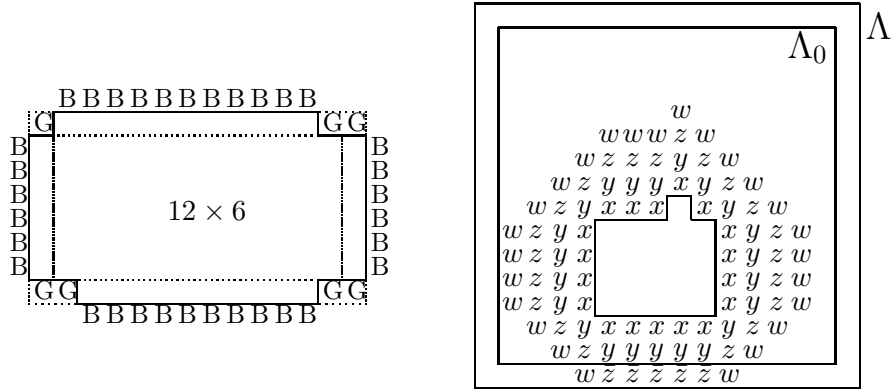
\begin{figure}
	\begin{picture}(-2,30)(-15,0)
		
		
		
		\thinlines
		\qbezier[51](0,5)(14,5)(28,5)
		\qbezier[51](28,5)(28,13)(28,21)
		\qbezier[51](0,21)(14,21)(28,21)
		\qbezier[51](0,5)(0,13)(0,21)
		
		\thinlines
		\qbezier[51](2,7)(13,7)(26,7)
		\qbezier[51](26,7)(26,13)(26,19)
		\qbezier[51](2,19)(13,19)(26,19)
		\qbezier[51](2,7)(2,13)(2,19)
		
		\put(0,7){\line(1,0){4}}
		\put(4,5){\line(1,0){20}}
		\put(24,7){\line(1,0){4}}
		\put(24,19){\line(1,0){4}}
		\put(2,21){\line(1,0){22}}
		\put(0,19){\line(1,0){2}}
		
		\put(4,5){\line(0,1){2}}
		\put(24,5){\line(0,1){2}}
		\put(28,7){\line(0,1){2}}
		\put(28,9){\line(0,1){10}}
		\put(24,19){\line(0,1){2}}
		\put(2,19){\line(0,1){2}}
		\put(0,7){\line(0,1){12}}
		
		\put(0.4,5.4){\small G}
		\put(2.4,5.4){\small G}
		
		\put(24.4,5.4){\small G}
		\put(26.4,5.4){\small G}
		\put(24.4,19.4){\small G}
		\put(26.4,19.4){\small G}
		\put(0.4,19.4){\small G}

		\put(4.4,3.4){\small B}
		\put(6.4,3.4){\small B}
		\put(8.4,3.4){\small B}
		\put(10.4,3.4){\small B}
		\put(12.4,3.4){\small B}
		\put(14.4,3.4){\small B}
		\put(16.4,3.4){\small B}
		\put(18.4,3.4){\small B}
		\put(20.4,3.4){\small B}
		\put(22.4,3.4){\small B}
		
		\put(28.4,7.4){\small B}
		\put(28.4,9.4){\small B}
		\put(28.4,11.4){\small B}
		\put(28.4,13.4){\small B}
		\put(28.4,15.4){\small B}
		\put(28.4,17.4){\small B}
		\put(2.4,21.4){\small B}
		\put(4.4,21.4){\small B}
		\put(6.4,21.4){\small B}
		\put(8.4,21.4){\small B}
		\put(10.4,21.4){\small B}
		\put(12.4,21.4){\small B}
		\put(14.4,21.4){\small B}
		\put(16.4,21.4){\small B}
		\put(18.4,21.4){\small B}
		\put(20.4,21.4){\small B}
		\put(22.4,21.4){\small B}
		\put(-1.6,7.4){\small B}
		\put(-1.6,9.4){\small B}
		\put(-1.6,11.4){\small B}
		\put(-1.6,13.4){\small B}
		\put(-1.6,15.4){\small B}
		\put(-1.6,17.4){\small B}

		\put(11.5,12){$12 \times 6$}

		\put(37,-2){\line(1,0){32}}
		\put(37,-2){\line(0,1){32}}
		\put(69,-2){\line(0,1){32}}
		\put(37,30){\line(1,0){32}}
		\put(39,0){\line(1,0){28}}
		\put(39,0){\line(0,1){28}}
		\put(39,28){\line(1,0){28}}
		\put(67,28){\line(0,-1){28}}
		\put(69.5,27){\Large{$\L$}}
		\put(63,25.5){\Large{$\L_0$}}
		\put(47,4){\line(1,0){10}}
		\put(47,4){\line(0,1){8}}
		\put(47,12){\line(1,0){6}}
		\put(53,12){\line(0,1){2}}
		\put(53,14){\line(1,0){2}}
		\put(55,12){\line(0,1){2}}
		\put(55,12){\line(1,0){2}}
		\put(57,12){\line(0,-1){8}}
		\put(57.5,10.5){$x$}
		\put(57.5,8.5){$x$}
		\put(57.5,6.5){$x$}
		\put(57.5,4.5){$x$}
		\put(55.5,2.5){$x$}
		\put(53.5,2.5){$x$}
		\put(51.5,2.5){$x$}
		\put(49.5,2.5){$x$}
		\put(47.5,2.5){$x$}
		\put(45.5,4.5){$x$}
		\put(45.5,6.5){$x$}
		\put(45.5,8.5){$x$}
		\put(45.5,10.5){$x$}
		\put(47.5,12.5){$x$}
		\put(49.5,12.5){$x$}
		\put(51.5,12.5){$x$}
		\put(53.5,14.5){$x$}
		\put(55.5,12.5){$x$}
		\put(57.5,12.6){$y$}
		\put(59.5,10.6){$y$}
		\put(59.5,8.6){$y$}
		\put(59.5,6.6){$y$}
		\put(59.5,4.6){$y$}
		\put(57.5,2.6){$y$}
		\put(55.5,0.7){$y$}
		\put(53.5,0.7){$y$}
		\put(51.5,0.7){$y$}
		\put(49.5,0.7){$y$}
		\put(47.5,0.7){$y$}
		\put(45.5,2.6){$y$}
		\put(43.5,4.6){$y$}
		\put(43.5,6.6){$y$}
		\put(43.5,8.6){$y$}
		\put(43.5,10.6){$y$}
		\put(45.5,12.6){$y$}
		\put(47.5,14.5){$y$}
		\put(49.5,14.5){$y$}
		\put(51.5,14.5){$y$}
		\put(53.5,16.5){$y$}
		\put(55.5,14.5){$y$}
		\put(57.5,14.5){$z$}
		\put(59.5,12.5){$z$}
		\put(61.5,10.5){$z$}
		\put(61.5,8.5){$z$}
		\put(61.5,6.5){$z$}
		\put(61.5,4.5){$z$}
		\put(59.5,2.5){$z$}
		\put(57.5,0.5){$z$}
		\put(55.5,-1.5){$\bar z$}
		\put(53.5,-1.5){$\bar z$}
		\put(51.5,-1.5){$\bar z$}
		\put(49.5,-1.5){$\bar z$}
		\put(47.5,-1.5){$\bar z$}
		\put(45.5,0.5){$z$}
		\put(43.5,2.5){$z$}
		\put(41.5,4.5){$z$}
		\put(41.5,6.5){$z$}
		\put(41.5,8.5){$z$}
		\put(41.5,10.5){$z$}
		\put(43.5,12.5){$z$}
		\put(45.5,14.5){$z$}
		\put(47.5,16.5){$z$}
		\put(49.5,16.5){$z$}
		\put(51.5,16.5){$z$}
		\put(53.5,18.5){$z$}
		\put(55.5,16.5){$z$}
		\put(57.3,16.5){$w$}
		\put(59.3,14.5){$w$}
		\put(61.3,12.5){$w$}
		\put(63.3,10.5){$w$}
		\put(63.3,8.5){$w$}
		\put(63.3,6.5){$w$}
		\put(63.3,4.5){$w$}
		\put(61.3,2.5){$w$}
		\put(59.3,0.5){$w$}
		\put(57.3,-1.5){$w$}
		\put(45.3,-1.5){$w$}
		\put(43.3,0.5){$w$}
		\put(41.3,2.5){$w$}
		\put(39.3,4.5){$w$}
		\put(39.3,6.5){$w$}
		\put(39.3,8.5){$w$}
		\put(39.3,10.5){$w$}
		\put(41.3,12.5){$w$}
		\put(43.3,14.5){$w$}
		\put(45.3,16.5){$w$}
		\put(47.3,18.5){$w$}
		\put(49.3,18.5){$w$}
		\put(51.3,18.5){$w$}
		\put(53.3,20.5){$w$}
		\put(55.3,18.5){$w$}

	\end{picture}
	\vskip 0.5 cm
	\caption{On the left-hand side we represent good sites ($G$) and bad sites ($B$) for $l_2^*=8$. On the right-hand side we depict with $x$ the sites in $B_1(\hat\h)$, with $y$ the sites in $\bar{B}_2(\hat\h)$, with $z$ and $\bar z$ the sites in $\bar{B}_3(\hat\h)$ and with $\bar z$ and $w$ the sites in $\bar{B}_4(\hat\h)$.}
	\label{fig:sitiGB} 
\end{figure}

and for $i=4,5,...,L_{sa}^*$
\be{}
\ba{lll}
B_i(\hat\h):=\{x\in\L_4| \ x\notin B_{i-2}(\hat\h), \ d(x,B_{i-1}(\hat\h))=1\},  \\
\bar B_i(\hat\h):=B_i(\hat\h)\cup \{\bar B_{i-1}(\hat\h)\cap \partial^-\L_4\}.
\ea
\ee
\ed

In words, $B_1(\bar\h)$ is the ring of sites in $\L_4$ at distance $1$ from $\hat\h$, while $\bar B_i(\hat\h)$ is the ring of sites in $\L_4$ at distance $i$ from $\hat\h$ union all the sites in $\partial^-\L_4$ at distance $1<j<i$ from $\hat\h$ ($i=2,3,...,L_{sa}^*$) (see Figure \ref{fig:sitiGB} on the right-hand side). Note that, depending on the location of $\hat\h$ in $\L$, the $\bar B_i(\hat\h)$ coincide for large enough $i$. The maximal number of rings is $L_{sa}^*$.  

Now we need to introduce specific sets that will be crucial later on.

\bd{c*i}
We define
\be{}
\cC_{sa}^*(i):=\{(\hat\h,x): \ \hat\h\in\cD_{sa}, \ x\in\bar B_i(\hat\h)\}, \quad i=2,3,...,L_{sa}^*.
\ee
\ed

\noindent
First, note that the sets $\cC_{sa}^*(i)$ are not disjoint. 

\br{}
From the definitions of the set $\csgeo$, we deduce that
\be{}
\cC_{sa}^*=\displaystyle\bigcup_{i=2}^{L_{sa}^*}\cC_{sa}^*(i).
\ee
\er

\noindent
For this discussion in the case $int\in\{is,wa\}$ we refer to \cite[Section 7.2]{BN2}.

\subsection{Useful lemmas for the gates}
\label{lemmi3modelli}
In this Section we give some useful Lemmas that help us to characterize the gates. Here we state Lemma \ref{trenini} for the case $int=sa$, but it holds also for the case $int\in\{is,wa\}$ (see \cite[Lemma 7.3]{BN2}). The proof is the same done in \cite[Lemma 7.3]{BN2} for the weakly anisotropic case.

\bl{trenini}
Starting from $\cC_{sa}^*\setminus\cQ_{sa}^{fp}$, if the free particle is attached to a bad site obtaining $\h^B\in\cC^B_{sa}$, the only transitions that does not exceed the energy $\G_{sa}^*$ are either detaching the protuberance, or a sequence of $1$-translations of a bar or slidings of a bar around a frame-angle. Moreover, we get:

\bi
\item[(i)] if it is possible to slide a bar around a frame-angle, then the saddles that are crossed are essential; 
\item[(ii)] if it is not possible to slide a bar around a frame-angle, then the path must come back to the starting configuration and the saddles that are crossed are unessential.
\ei
\el

\noindent
Lemmas \ref{coltorow} and \ref{dtilde} are valid also in the case $int=wa$ (see \cite[Lemma 7.5]{BN2} and \cite[Lemma 7.6]{BN2} respectively), while Lemma \ref{trasless} has a corresponding version for the case $int=wa$ (see \cite[Lemma 7.4]{BN2}).

\bl{trasless}
Starting from $\h^B\in\cC_{sa}^B$, the saddles obtained by a $1$-translation of a bar are essential and in $\cA_{0}^{\a'}\cup\cA_1^{\a}$. Moreover, all the saddles in  $\cA_{0}^{\a'}\cup\cA_1^{\a}$ can be obtained from this $\h^B$ by a $1$-translation of a bar. 
\el

\bpr
Note that $H(\h^B)=\G_{sa}^*-U_2$ (resp.\ $H(\h^B)=\G_{sa}^*-U_1$) if the free particle has been attached to an horizontal (resp.\ vertical) bar. In the first case, in order to avoid exceeding the energy value $\G_{sa}^*$ it is possible to translate only the vertical bars. These saddles are in $\cA_{1}^{\a}$. In the latter case, it is possible to translate both vertical and horizontal bars. If the translated bar is horizontal, the saddles that are crossed are in $\cA_{0}^{\a'}$. If the translated bar is vertical, the configurations obtained do not reach the level $\G_{sa}^*$, thus they are not saddles. To conclude, all the configurations in $\cA_{0}^{\a'}\cup\cA_{1}^{\a}$ can be obtained from a configuration $\h^B$ via a $1$-translation of a bar.

It remains to prove that the saddles in $\cA_{0}^{\a'}\cup\cA_{1}^{\a}$ are essential. This part of the proof is analogue to the corresponding one done for $int=is$ in \cite[Lemma 8.2]{BN2}.
\epr

\bl{coltorow}
Starting from a configuration $\h\in\cC_{sa}^*$, it is not possible to slide a vertical bar around a frame-angle without exceeding the energy $\G_{sa}^*$.
\el

\noindent
With the following Lemma we can justify the definition of $\csgeo$ given in (\ref{c*sa}).

\bl{dtilde}
Starting from $\widetilde\cD_{sa}$, the dynamics either passes through $\bar\cD_{sa}$ or it is not possible that a free particle is created without exceeding the energy level $\G_{sa}^*$.
\el

\noindent
The proof of Lemmas \ref{coltorow} and \ref{dtilde} are analogue to the ones done for \cite[Lemma 7.5]{BN2} and \cite[Lemma 7.6]{BN2} respectively in the case $int=wa$ by replacing $\gw$ (resp.\ $\widetilde\cD_{wa}$) with $\gs$ (resp.\ $\widetilde\cD_{sa}$).

\subsection{Model-dependent strategy}
\label{S6.4}
Our goal is to characterize the union of all the minimal gates for the strongly anisotropic interactions. To this end, due to \cite[Theorem 5.1]{MNOS}, we will characterize all the essential saddles for the transition from the metastable to the stable state. In this Section we apply the model-independent strategy explained in \cite[Section 3.1]{BN2} in order to identify some unessential saddles. Let $int=sa$. We apply \cite[eq.\ (3.1)]{BN2} both for $\s=\vuoto$, $\cA=\{\pieno\}$ and $\G=\G_{sa}^*$ defining $\cC_{\pieno}^{\vuoto}(\G_{sa}^*)$, and for $\s=\pieno$, $\cA=\{\vuoto\}$ and $\G=\G_{sa}^*-H(\pieno)$ defining $\cC_{\vuoto}^{\pieno}(\G_{sa}^*-H(\pieno))$. We chose this notation in order to emphasize the dependence on $\G_{sa}^*$. First, we prove the required model-dependent inputs (iii)-(a) and (iii)-(b) in \cite[Section 3.1]{BN2} (see Proposition \ref{pathstrong}(i) and Proposition \ref{pathstrong}(ii)). Second, by Theorem \ref{thgate}, we know that $\cC_{sa}^*$ is a gate for the transition from $\vuoto$ to $\pieno$. Thus we apply the model-independent strategy explained in \cite[Section 3.1]{BN2} to Kawasaki dynamics by taking $m=\vuoto$, $\cX^s=\{\pieno\}$, $\cW(m,\cX^s)=\cC_{sa}^*$, $\cL^B=\cC_{sa}^B$ and $\cL^G=\cC_{sa}^G$. In Proposition \ref{c*contenuto} we prove that $\cC_{sa}^*\subseteq\cG_{sa}(\vuoto,\pieno)$, that allows us to study the essentiality only of the saddles that are not in $\cC_{sa}^*$.

In order to apply \cite[Propositions 3.3, 3.5]{BN2}, we need to characterize the sets $K_{sa}$ and $\widetilde{K}_{sa}$ (see \cite[eq.\ (3.2)]{BN2} and \cite[eq.\ (3.3)]{BN2} respectively for the definitions). This is done in Proposition \ref{Kstrong}. Due to this result, our strategy consists in partitioning the saddles that are not in $\cC_{sa}^*$ in three types: the saddles that are in the boundary of $\cC_{\pieno}^{\vuoto}(\G_{sa}^*)$, i.e., $\s\in\partial\cC_{\pieno}^{\vuoto}(\G_{sa}^*)\cap(\cS_{int}(\vuoto,\pieno)\setminus\cC_{sa}^*)$, the saddles that are in the boundary of $\cC_{\vuoto}^{\pieno}(\G_{sa}^*-H(\pieno))$ and not in $\widetilde{K}_{sa}$, i.e., $\z\in\partial\cC_{\vuoto}^{\pieno}(\G_{sa}^*-H(\pieno))\cap(\cS_{sa}(\vuoto,\pieno)\setminus(\cC_{sa}^*\cup\widetilde{K}_{sa}))$, and the remaining saddles $\x\in
\cS_{sa}(\vuoto,\pieno)\setminus(\partial\cC_{\pieno}^{\vuoto}(\G_{sa}^*)\cup(\partial\cC_{\vuoto}^{\pieno}(\G_{sa}^*-H(\pieno))\setminus\widetilde{K}_{sa})\cup\cC_{int}^*)$. By \cite[Propositions 3.3, 3.5]{BN2}, we obtain Corollary \ref{corstrategy} that states that the saddles of the first and second types are respectively unessential. In Proposition \ref{selle3bis} we highlight some of the saddles of type three that are unessential. We need to distinguish these two cases due to the different entrance in $\cC_{int}^*$ for $int\in\{is,wa\}$ and $int=sa$ (see \cite[Lemma 7.13]{BN2} and Lemma \ref{entratastrong} respectively). Note that in Proposition \ref{selle3bis} the set $\partial\cC_{\vuoto}^{\pieno}(\G_{sa}^*-H(\pieno))\setminus\widetilde{K}_{sa}$ reduces to $\partial\cC_{\vuoto}^{\pieno}(\G_{sa}^*-H(\pieno))$ due to Proposition \ref{Kstrong}(ii). For the case $int\in\{is,wa\}$ this strategy is presented in \cite[Section 7.4]{BN2}. Finally, we identify the essential saddles of the third type in Proposition \ref{strongselless}.

\subsubsection{Main Propositions}
In this Subsection we give the main results for our model-dependent strategy. 

The next proposition shows that when the dynamics reaches $\cC_{sa}^G$ it has gone ``over the hill", while when it reaches $\cC_{sa}^B$ the energy has to increase again to the level $\gs$ to visit $\vuoto$ or $\pieno$. An analogue version for $int=is$ is proven in \cite[Proposition 2.3.9]{BHN} and for $int=wa$ is proven in \cite[Proposition 7.7]{BN2}, while here we extend that result to $int=sa$ following a similar argument.

\bp{pathstrong}
The following hold:
\bi
\item[(i)]  If $\h\in\cC_{sa}^G$, then there exists a path $\o:\h\ra\pieno$ such that $\max_{\z\in\o}H(\z)<\G_{sa}^*$. 
\item[(ii)] If $\h\in\cC_{sa}^B$, then there are no $\o:\h\ra\vuoto$ or $\o:\h\ra\pieno$ such that $\max_{\z\in\o}H(\z)<\G_{sa}^*$.
\ei
\ep

\begin{proof*}
	The proof is analogue to the one done in \cite[Proposition 7.7]{BN2} for the weakly anisotropic case by using the reference path for the nucleation constructed in \cite[Section 3.2]{BN}. 
\end{proof*}

\noindent
The next Proposition holds also in the case $int\in\{is,wa\}$ (see \cite[Proposition 7.8]{BN2}). We refer to Subsection \ref{proofproposizionidep} for the proof of Propositions \ref{c*contenuto}, \ref{Kstrong} and \ref{selle3bis}.

\bp{c*contenuto}
$\cC_{sa}^*\subseteq\cG_{sa}(\vuoto,\pieno)$.
\ep

\noindent
For the corresponding result of Proposition \ref{c*contenuto} for $int\in\{is,wa\}$ we refer to \cite[Proposition 7.8]{BN2}.

\bp{Kstrong}
The following hold:
\bi
\item[(i)] $K_{sa}=\emptyset$; 
\item[(ii)]$\widetilde{K}_{sa}\cap\partial\cC_{\vuoto}^{\pieno}(\G_{sa}^*-H(\pieno))=\emptyset$.
\ei
\ep

\noindent
For the corresponding result of Proposition \ref{Kstrong} for $int\in\{is,wa\}$ we refer to \cite[Proposition 7.9]{BN2}.

\bc{corstrategy}
The following hold:
\bi
\item[(i)]The saddles of the first type $\s\in\partial\cC_{\pieno}^{\vuoto}(\G_{sa}^*)\cap(\cS_{sa}(\vuoto,\pieno)\setminus\cC_{sa}^*)$ are unessential;
\item[(ii)] The saddles of the second type $\z\in\partial\cC_{\vuoto}^{\pieno}(\G_{sa}^*-H(\pieno))\cap(\cS_{sa}(\vuoto,\pieno)\setminus\cC_{sa}^*)$ are unessential.
\ei
\ec

\begin{proof*}
	Combining \cite[Propositions 3.3, 3.5]{BN2} and \ref{Kstrong} we get the claim.
\end{proof*}

\noindent
For the corresponding result of Corollary \ref{corstrategy} for $int\in\{is,wa\}$ we refer to \cite[Corollary 7.10]{BN2}.

\bp{selle3bis}
Any saddle $\x$ that is neither in $\cC_{sa}^*$, nor in $\bigcup_{k,\a,\a'}\cA_{k}^{\a,\a'}$, nor in the boundary of the cycles $\cC_{\pieno}^{\vuoto}(\G_{sa}^*)$ or $\cC_{\vuoto}^{\pieno}(\G_{sa}^*-H(\pieno))$, i.e.,  $\x\in\cS_{sa}(\vuoto,\pieno)\setminus(\partial\cC_{\pieno}^{\vuoto}(\G_{sa}^*)\cup\partial\cC_{\vuoto}^{\pieno}(\G_{sa}^*-H(\pieno))\cup\cC_{sa}^*\cup\bigcup_{k,\a,\a'}\cA_{k}^{\a,\a'})$, such that $\t_{\x}<\t_{\cC^B_{sa}}$ is unessential. Therefore it is not in $\cG_{sa}(\vuoto,\pieno)$.
\ep

\noindent
For the corresponding result of Proposition \ref{selle3bis} for $int\in\{is,wa\}$ we refer to \cite[Proposition 7.11]{BN2}.

\subsubsection{Useful Lemmas for the model-dependent strategy}
\label{lemminuovi}
In this Subsection we give some useful lemmas about the entrance in the gate and some properties of the sets $\cC_{sa}^*(i)$ with $i=3,...,L_{sa}^*$. We stress that the behavior for $int\in\{is,wa\}$ is very different from that observed for $int=sa$, indeed we note that the weakly anisotropic model has some characteristics similar to the isotropic and some similar to the strongly anisotropic model. For the corresponding results obtained in the case $int\in\{is,wa\}$ we refer to \cite[Subsection 7.4.2]{BN2}. Recall (\ref{definizioneqbarsa}) and (\ref{c*sa}) for the definitions of $\cQ_{sa}$, $\cD_{sa}$ and $\csgeo$. The next lemma generalizes \cite[Proposition 2.3.8]{BHN}, proved for $int=is$, to the case $int=sa$ following similar arguments. In the case $int=wa$, this result is given in \cite[Lemma 7.12]{BN2}. The proof of Lemma \ref{motion} is analogue to the one done in \cite[Lemma 7.12]{BN2} for the weakly anisotropic case.

\bl{motion}
\bi
\item[(i)] Starting from $\cC_{sa}^*\setminus\cQ_{sa}^{fp}$, the only transitions that do not raise the energy are motions of the free particle in the region where the free particle is at lattice distance $\geq3$ from the protocritical droplet.
\item[(ii)] Starting from $\cQ_{sa}^{fp}$, the only transitions that do not raise the energy are motions of the free particle in the region where the free particle is at lattice distance $\geq3$ from the protocritical droplet and motions of the protuberance along the side of the rectangle where it is attached. When the lattice distance is 2, either the free particle can be attached to the protocritical droplet or the protuberance can be detached from the protocritical droplet and attached to the free particle, to form a rectangle plus a dimer. From the latter configuration the only transition that does not raise the energy is the reverse move.
\item[(iii)] Starting from $\cC_{sa}^*$, the only configurations that can be reached by a path that lowers the energy and does not decrease the particle number, are those where the free particle is attached to the protocritical droplet.
\ei
\el

\bl{C*2ess}
The saddles in $\cC_{sa}^*(2)$ are essential.
\el

\noindent
The proof of Lemma \ref{C*2ess} is analogue to the one done in \cite[Lemma 7.16]{BN2} for the case $int\in\{is,wa\}$. The next Lemma investigates how the entrance in $\cC_{sa}^*$ occurs. For the corresponding result that holds in the case $int\in\{is,wa\}$ see \cite[Lemma 7.13]{BN2}. We encourage the reader to inspect the difference between Lemma \ref{entratastrong} and \cite[Lemma 7.13]{BN2}, indeed the entrance in the gate in the strongly anisotropic case is peculiar and different with respect the isotropic and weakly anisotropic ones. Recall (\ref{defpsatrenino}) and (\ref{insiemeA}).

\bl{entratastrong}
Any $\o\in(\vuoto\ra\pieno)_{opt}$ enters $\cC_{sa}^*$ in one of the following ways:
\bi
\item[(i)] $\o$ passes first through $\bar\cQ_{sa}$, then possibly through $\bar\cD_{sa}\setminus\bar\cQ_{sa}$, and finally reaches $\cC_{sa}^*$;
\item[(ii)] $\o$ passes through the configuration $\cR(2l_2^*-1,l_2^*-1)$, then a free particle is created and moved towards the rectangle until it is attached to an horizontal side $\a\in\{n,s\}$. Then for some $\a'\in\{w,e\}$ the path $\o$ passes through the sets $\cA_{k}^{\a,\a'}$ for all $k=2,...,l_2^*$, and finally reaches $\cC_{sa}^*(2)$.
\ei
\el

\noindent
Since in Lemma \ref{entratastrong} we have proved that there are two possible ways to reach $\cC_{sa}^*(2)$, where the possibility (i) is analogue to the cases $int=is$ and $int=wa$, to find the minimal gates for $sa$ for any $i=3,...,L_{sa}^*$ we need to consider $\cC_{sa}^*(i)$ union some particular saddles belonging to the paths described in (ii).

\bl{mingatecornici}
For any $i=3,...,L_{sa}^*$ and $k=2,...,l_2^*$ the set $\cC_{sa}^*(i)\cup\bigcup_{\a,\a'}\{\x_{j(k)}^{\a,\a'}\}$ is a minimal gate for all $1\leq j(k)\leq k-1$, where $\{\x_{j(k)}^{\a,\a'}\}$ are the elements in $\cA_{k}^{\a,\a'}$ defined in (\ref{defpsatrenino}) and (\ref{insiemeA}).
\el

\br{}
We encourage the reader to inspect the difference between the statement of Lemma \ref{mingatecornici} and \cite[Lemma 7.14]{BN2}: the sets $\cC_{int}^*(i)$, with $i=3,...,L_{int}^*$, are minimal gates if $int\in\{is,wa\}$, while $\cC_{sa}^*(i)$ are not minimal gates for any $i=3,...,L_{sa}^*$. 
\er

\bl{strongminimalgates}
For the strongly anisotropic interactions, we have 
\be{}
\cC_{sa}^*\cup\bigcup_{k=2}^{l_2^*}\bigcup_{\a,\a'}\cA_k^{\a,\a'}\subseteq\cG_{sa}(\vuoto,\pieno)
\ee
\el

\subsubsection{Proof of Propositions}
\label{proofproposizionidep}

\begin{proof*}{\bf of Proposition \ref{c*contenuto}}
	The statement of the Proposition follows by Lemma \ref{strongminimalgates}.
\end{proof*}

\begin{proof*}{\bf of Proposition \ref{Kstrong}}
	The proof of (i) is analogue to the one done in \cite[Proposition 7.9]{BN2} for $int\in\{is,wa\}$.
	
	Now we prove (ii).  Let $\bar\h\in\widetilde{K}_{sa}\cap\partial\cC_{\vuoto}^{\pieno}(\G_{sa}^*-H(\pieno))$. By the definition of the set $\widetilde{K}_{sa}$ we know that there exist $\h\in\cC_{sa}^*$ and $\o=\o_1\circ\o_2$ from $\h$ to $\pieno$ with the properties described in \cite[eq.\ (3.3)]{BN2}. We know that $\h$ is composed by the union of a protocritical droplet $\hat\h\in\cD_{sa}$ and a free particle. Since $\o_1\cap\cC_{sa}^*=\{\h\}$, we note that $\h\in\cC_{sa}^*(2)$, otherwise the free particle has to cross at least $\bar{B}_2(\hat\h)$ and $\bar{B}_3(\hat\h)$, the latter in the configuration $\h'\in\cC_{sa}^*$, with $\h'\neq\h$, which contradicts the conditions in \cite[eq.\ (3.3)]{BN2}. Therefore, starting from $\h$, by the optimality of the path we deduce that the unique admissible move is to attach the free particle to the cluster. If $\bar\h$ is obtained from $\h$ by attaching the free particle in a good site giving rise to a configuration in $\cC_{sa}^G(\hat\h)$, by Proposition \ref{pathstrong}(i) we know that $\o_1\cap\cC_{\vuoto}^{\pieno}(\G_{sa}^*-H(\pieno))\neq\emptyset$, that contradicts \cite[eq.\ (3.3)]{BN2}, thus it is not possible to find $\o_1$ and $\o_2$, therefore $\bar\h\notin\widetilde{K}_{sa}$, which is in contradiction with the assumption. 
	
	Assume now that $\bar\h$ is obtained from $\h$ by attaching the free particle in a bad site giving rise to a configuration in $\cC_{sa}^B(\hat\h)$. If $\h\in\cQ_{sa}^{fp}$, then by Lemma \ref{motion}(ii) the unique admissible move is the reverse one, thus we may assume that $\h\in\cC_{sa}^*\setminus\cQ_{sa}^{fp}$ and that the path does not go back to $\h$, otherwise we can iterate this argument for a finite number of steps since the path has to reach $\pieno$. Starting from $\h$, by Lemma \ref{trenini} we know that $\bar\h$ is obtained either via a sequence of $1$-translations of a bar or via a sliding of a bar around a frame-angle. If a sequence of $1$-translations takes place, by the optimality of the path we deduce that the unique possibility is either detaching the protuberance or sliding a bar around a frame-angle. In the first case the configuration that is obtained is in $\cC_{sa}^*$ and thus $\bar\h\notin\widetilde{K}_{sa}$, which contradicts the assumption. by (\ref{condtrenino}), Proposition \ref{card}(b) (in particular conditions in (\ref{condbard})) and Lemma \ref{coltorow} we deduce that the only possibility to slide a bar around a frame-angle is that the bar is horizontal and it has length exactly $l_2^*-1$. Thus the configurations visited by the path $\o$ during this sliding are $\bar\h_1,...,\bar\h_m\in\cA_k^{\a,\a'}$ for some $\a\in\{n,s\}$, $\a'\in\{w,e\}$ and $k=2,...,l_2^*-1$, while the last configuration $\widetilde\h$ obtained when the last particle of the bar is detached is composed by the union of $\cR(2l_2^*-1,l_2^*-1)$ and a free particle (see the time-reversal of the path described in Figure \ref{fig:columntorow}, in particular $\bar\h_m$ is the configuration (7) and $\widetilde\h$ is the configuration (2)). Therefore $\widetilde\h$ belongs to the set $\cB$ defined in \cite[eq. (3.29)]{BN} since $s(\widetilde\h)=s_{sa}^*-1$ and $p_2(\widetilde\h)=l_2^*-1$. Thus by \cite[Theorem 3.7]{BN} we deduce that $\widetilde\h\notin\cC_{\vuoto}^{\pieno}(\G_{sa}^*-H(\pieno))$ and therefore  $\bar\h_m\notin\partial\cC_{\vuoto}^{\pieno}(\G_{sa}^*-H(\pieno))$, that implies $\widetilde{K}_{sa}\cap\partial\cC_{\vuoto}^{\pieno}(\G_{sa}^*-H(\pieno))=\emptyset$. 
\end{proof*}

\begin{proof*}{\bf of Proposition \ref{selle3bis}}
	The proof is analogue to the one done for \cite[Proposition 7.11]{BN2} in the case $int\in\{is,wa\}$, but in this case we use Lemmas \ref{motion}(ii) and \ref{entratastrong} and $\h$ is the union of a cluster $\hat\h\in\cQ_{sa}$ and a free particle at distance 2 from the cluster. Moreover, $\x_i$ is the union of a rectangle $(2l_2^*-3)\times l_2^*$ with an horizontal dimer.
\end{proof*}

\subsubsection{Proof of Lemmas}
\label{prooflemmiapp}

\begin{proof*}{\bf of Lemma \ref{entratastrong}}
	By Theorem \ref{thgate} we know that any $\o\in(\vuoto\ra\pieno)_{opt}$ passes through $\cC_{sa}^*$. We denote by $\h$ this configuration, that is composed by the union of a protocritical droplet $\hat\h\in\bar\cD_{sa}$ and a free particle in the site $x$. Note that there exists $i=2,...,L_{sa}^*$ such that either $x\in B_i(\hat\h)$ if $d(\partial^-\L_4,\hat\h)>i$ or $x\in\bar B_i(\hat\h)$ if $d(\partial^-\L_4,\hat\h)\leq i$. We set $\o=(\vuoto,\o_1,...,\o_k,\h)\circ\bar{\o}$, where $\bar\o$ is a path that connects $\h$ to $\pieno$ such that $\max_{\s\in\o}H(\s)\leq\gs$. In order to analyze the entrance in $\cC_{sa}^*$ we consider the time-reversal of the path $\o$. Since $H(\h)=\G_{sa}^*$, the move from $\h$ to $\o_k$ must have a non-positive cost and thus the unique admissible moves are:
	\begin{description}
		\item[(i)] either moving the free particle at zero cost; 
		\item[(ii)] or removing a free particle at cost $-\Delta$; 
		\item[(iii)] or attaching the free particle at cost $-U_1$ (resp.\ $-U_2$) or $-U_1-U_2$.
	\end{description}
	
	\medskip
	\noindent
	{\bf Case (i).} In this case we obtain that the configuration $\o_k$ is still in $\cC_{sa}^*$, thus it is analogue to $\h$ and therefore we can iterate the argument by taking this configuration as $\h$.
	
	\medskip
	\noindent
	{\bf Case (ii).} In this case $H(\o_k)=\G_{sa}^*-\D$ and $\o_k\in\bar\cD_{sa}$. If $\o_k\in\bar\cQ_{sa}$ we get the claim, thus in the sequel we assume that $\o_k\in\bar\cD_{sa}\setminus\bar\cQ_{sa}$. Since the path $\o$ starts from $\vuoto$, there exist $k_1<k_2<k$ such that $|\o_{k_1}|=|\o_k|-1$ and there is a free particle in $\o_{k_2}$, i.e., $n(\o_{k_2})=1$. Starting from $\o_k$ and considering the time-reversal of the path $\o$, in order to obtain a free particle in $\o_{k_2}$ we note that the minimal cost for detaching a particle is $U_1+U_2$ giving rise to the energy value greater or equal than $\G_{sa}^*-\D+U_1+U_2>\G_{sa}^*$, which is in contradiction with the optimality of the path. Thus the unique possibility is detaching the protuberance from a configuration in $\bar\cQ_{sa}$ at cost $U_1$. 
	This implies that $\o_k$ is obtained via a $U_1$-path starting from a configuration in $\bar\cQ_{sa}$.

	\medskip
	\noindent
	{\bf Case (iii).} First, we consider the case where from $\h$, again considering the time-reversal, we attach a particle at cost $-U_1$ giving rise to the configuration $\o_k$, i.e., $H(\o_k)=\gs-U_1$. Since the path $\o$ starts from $\vuoto$, there exists $k_1<k$ such that $|\o_{k_1}|=|\o_{k}|-1$, that implies that there exists a configuration $\o_{\bar{k}}$ with a free particle during the transition from $\o_{k_1}$ to $\o_{k}$ (see Figure \ref{fig:columntorow} where $\o_k$ is configuration (12) and $\o_{\bar{k}}$ is configuration (2)). If $\o_{\bar{k}}\in\cC_{sa}^*$, we can iterate the argument by taking this configuration as $\h$. Otherwise if $\o_{\bar{k}}\notin\cC_{sa}^*$, we deduce that $(\o_{\bar{k}})_{cl}\notin\bar\cD_{sa}$. Starting from $\o_{\bar{k}}$, since the activation of a sequence of $1$-translations of bars of configurations in $\bar\cD_{sa}$ gives rise to configurations that are in $\bar\cD_{sa}$, the unique possibility in order not to exceed $\gs$ is that $\o_k$ is obtained from $\o_{\bar{k}}$ via a sliding of a bar, say $B^{\a'}(\o_k)$, around a frame-angle, say $c^{\a'\a}(\o_k)$.
	In order to do that, by (\ref{condtrenino}) and Proposition \ref{card}(b), we deduce that the unique possibility to match the two conditions in (\ref{condbard}) is that during the transition the path $\o$ crosses $\cS_{sa}(\vuoto,\pieno)$ through the sets $\cA_{k}^{\a,\a'}$ for any $k=2,...,l_2^*$, with $\a\in\{n,s\}$, $\a'\in\{w,e\}$. In Figure \ref{fig:columntorow} we represent this transition with $\a=n$, $\a'=e$ and the configurations (11)-(3) for the sliding. At the end of this sliding we obtain the configuration $\cR(2l_2^*-1,l_2^*-1)$ union a free particle. This configuration must be obtained from $\cR(2l_2^*-1,l_2^*-1)$ via adding a free particle, otherwise the path $\o$ is not optimal. 
	
	Second, we consider the case where from $\h$ we attach a particle at cost $-U_2$ giving rise to the configuration $\o_k$, i.e., $H(\o_k)=\gs-U_2$. We argue in a similar way as above, but the difference is that in this case the sliding of a bar around a frame-angle at cost $U_2$ is not allowed by Lemma \ref{coltorow}.

	Third, we consider the case where from $\h$ we attach a particle at cost $-U_1-U_2$ giving rise to the configuration $\o_k$, i.e., $H(\o_k)=\gs-U_1-U_2$. Since $\h_{cl}\in\bar\cD_{sa}$, the unique possibility is that $\o_{{k}}\in\cC_{sa}^G(\hat\h)$, therefore by Lemma \ref{pathstrong}(i) we get $\o_k\in\cC_{\vuoto}^{\pieno}(\gs-H(\pieno))$. Since by Theorem \ref{thgate} we know that $\csgeo$ is a gate for the transition, we deduce that there exists $k_1<k$ such that $\o_{k_1}\in\cC_{sa}^*$. Thus we can iterate the argument by taking this configuration as $\h$.
\end{proof*}

\begin{proof*}{\bf of Lemma \ref{mingatecornici}}
	Let $i\in\{3,...,L_{sa}^*\}$, $k\in\{1,...,n\}$ and $1\leq j(k)\leq k-1$. 
	
	First, we prove that $\cC_{sa}^*(i)\cup\bigcup_{\a,\a'}\{\x_{j(k)}^{\a,\a'}\}$ is a gate. By Theorem \ref{thgate} we know that any $\o\in(\vuoto\ra\pieno)_{opt}$ crosses $\csgeo$. If the path $\o$ enters $\csgeo$ without crossing the set $\cP_{sa,0}$, then by Lemma \ref{entratastrong}(i) we know that $\o$ has to pass through $\cC_{sa}^*(i)$. If the path $\o$ enters $\csgeo$ after crossing the set $\cP_{sa,0}$, then by Lemma \ref{entratastrong}(ii) we know that $\o\cap\bigcup_{\a,\a'}\{\x_{j(k)}^{\a,\a'}\}\neq\emptyset$. 
	
	Now we prove that $\cC_{sa}^*(i)\cup\bigcup_{\a,\a'}\{\x_{j(k)}^{\a,\a'}\}$ is a minimal gate by showing that for any $\h\in\cC_{sa}^*(i)\cup\bigcup_{\a,\a'}\{\x_{j(k)}^{\a,\a'}\}$ the set $(\cC_{sa}^*(i)\cup\bigcup_{\a,\a'}\{\x_{j(k)}^{\a,\a'}\})\setminus\{\h\}$ is not a gate: there exists $\o\in(\vuoto\ra\pieno)_{opt}$ such that $\o\cap((\cC_{sa}^*(i)\cup\bigcup_{\a,\a'}\{\x_{j(k)}^{\a,\a'}\})\setminus\{\h\})=\emptyset$. We consider separately the cases $\h\in\bigcup_{\a,\a'}\{\x_{j(k)}^{\a,\a'}\}$ and $\h\in\cC_{sa}^*(i)$.
	
	\medskip
	\noindent
	{\bf Case 1.} Let $\h\in\bigcup_{\a,\a'}\{\x_{j(k)}^{\a,\a'}\}$, thus $\h=\x_{j(k)}^{\bar\a,\bar\a'}$ for some $\bar\a\in\{n,s\}$ and $\bar\a'\in\{w,e\}$. We can define $\o$ as the reference path defined in \cite[Section 3.2]{BN} that crosses the configurations $\x_1^{\bar\a,\bar\a'}$,...,$\x_{k-1}^{\bar\a,\bar\a'}$, then it enters $\cC_{sa}^*(2)$ and finally the free particle is attached in a good site without passing through $\cC_{sa}^*(i)$ with $i=3,...,L_{sa}^*$ (see Figure \ref{fig:columntorow}). From this configuration, the path proceeds towards $\pieno$ as the one in Proposition \ref{pathstrong}(i). The constructed $\o$ is optimal and $\o\cap\bigcup_{\a,\a'}\{\x_{j(k)}^{\a,\a'}\}=\{\x_{j(k)}^{\bar\a,\bar\a'}\}$, thus this case is concluded.
	
	\medskip
	\noindent
	{\bf Case 2.} Let $\h\in\csgeo(i)$. We take an arbitrary path starting from $\vuoto$ and that enters $\cC_{sa}^*(i)$ in $\h=(\hat\h,z)$, where $\hat\h\in\bar\cD_{sa}$ is the protocritical droplet and $z$ is the position of the free particle at distance $i$ from the cluster. Then the path proceeds by moving the free particle from $z$ to $\hat\h$ such that, the distance between the free particle and $\hat\h$ at the first step is strictly decreasing, and at the later steps is not increasing. Finally the free particle is attached in  a good site $x\in\partial^-CR(\hat\h)$ giving rise to a configuration in $\cC_{sa}^G(\hat\h)$. From this configuration, the path proceeds towards $\pieno$ as the one in Proposition \ref{pathstrong}(i). Since the constructed $\o\in(\vuoto\ra\pieno)_{opt}$ and $\o\cap\cC_{sa}^*(i)=\{\h\}$, the proof is completed.
\end{proof*}

\begin{proof*}{\bf of Lemma \ref{strongminimalgates}}
	By Lemma \ref{C*2ess} we know that the saddles in $\cC_{sa}^*(2)$ are essential and thus are in the set $\cG_{sa}(\vuoto,\pieno)$ due to \cite[Theorem 5.1]{MNOS}. Furthermore, by Lemma \ref{mingatecornici} we know that $\cC_{sa}^*(i)\cup\bigcup_{\a,\a'}\{\x_{j(k)}^{\a,\a'}\}$ is a minimal gate for any $i=3,...,L_{sa}^*$ and $j(k)=1,...,k-1$, with $k=2,...,l_2^*$. Therefore we get
	\be{} 
	\cG_{sa}(\vuoto,\pieno)\supseteq\cC_{sa}^*(2)\cup\displaystyle\bigcup_{i=3}^{L_{sa}^*}\bigcup_{k=2}^{l_2^*}\bigcup_{j(k),\a,\a'}(\cC_{sa}^*(i)\cup\{\x_{j(k)}^{\a,\a'}\})=\cC_{sa}^*\cup\displaystyle\bigcup_{k=2}^{l_2^*}\displaystyle\bigcup_{\a,\a'}\cA_k^{\a,\a'}
	\ee
\end{proof*}

\section{Proof of the main results: strongly anisotropic case}
\label{proofstrong}
In this Section we give the proof of the main Theorems \ref{thgate} and \ref{gstrong} (see Sections \ref{strong1} and \ref{strong2} respectively).

\subsection{Proof of the main Theorem  \ref{thgate}} 
\label{strong1}
In this Section we give the proof of the main Theorem \ref{thgate}. Now we recall the definition of the set $\cP_2$ given in \cite{BN} as
\be{defP2}
\ba{ll}
\cP_2:= \{\h&:\, n(\h)= 1,\, v(\h)=l_2^*-1,\,
\h_{cl} \hbox{ is connected},\, \hbox{monotone},\\
&\hbox{ with circumscribed rectangle in }
\cR(2l_2^*-2,l_2^*)\}.
\ea
\ee

\noindent
In particular, in order to state that the set $\csgeo$ is a gate for the transition from $\vuoto$ to $\pieno$, we need the following 

\bl{gatestrong}
If $\o\in(\vuoto\ra\pieno)_{opt}$ is such that $\o\cap\cP_2$, then $\o\cap\csgeo\neq\emptyset$.
\el

\noindent
We postpone the proof of Lemma \ref{gatestrong} after the proof of the main Theorem \ref{thgate}.

\bp{p2}
If $\o\in(\vuoto\ra\pieno)_{opt}$ is such that $\o\cap\cP_{sa,0}\neq\emptyset$, then $\o\cap\csgeo\neq\emptyset$.
\ep

\noindent
We postpone the proof of Proposition \ref{p2} after the proof of Lemma \ref{gatestrong}.

\medskip
\begin{proof*}{\bf of the main Theorem \ref{thgate}}
	By \cite[Theorem 2.4]{BN} taking $\cP_1=\cP_{sa,0}$, we know that the set $\cP_{sa,0}\cup\cP_2$ is a gate for the transition from $\vuoto$ to $\pieno$. By Lemma \ref{gatestrong} we know that every $\o\in(\vuoto\ra\pieno)_{opt}$ that crosses $\cP_2$ then crosses $\csgeo$, thus we deduce that the set $\cP_{sa,0}\cup\csgeo$ is a gate. Furthermore, by Proposition \ref{p2} we obtain that every path $\o\in(\vuoto\ra\pieno)_{opt}$ that crosses $\cP_{sa,0}$ then crosses also $\csgeo$. This implies that every optimal path $\o$ from $\vuoto$ to $\pieno$ is such that $\o\cap\csgeo\neq\emptyset$, thus $\csgeo$ is a gate.
\end{proof*}

\begin{proof*}{\bf of Lemma \ref{gatestrong}}
	Consider $\o\in(\vuoto\ra\pieno)_{opt}$. If $\o\cap\csgeo\neq\emptyset$, we get the claim. Thus we can reduce our analysis to the case in which the path $\o$ reaches the set $\cP_2$ in a configuration $\h\in\cP_2\setminus\csgeo$. We set $\o=(\vuoto,\o_1,...,\o_k,\h)\circ\bar\o$, where $\bar\o$ is a path that connects $\h$ to $\pieno$ such that $\max_{\s\in\o}H(\s)\leq\gs$. We are interested in the time-reversal of the path. Since $\h\in\cP_2\setminus\csgeo$, we know that it is composed by the union of a cluster $\hbox{CR}^-(\h)=\cR(2l_2^*-4,l_2^*-2)$, such that at least one frame-angle of $\hbox{CR}^-(\h)$ is empty, a free particle and four bars attached to the four sides of $\hbox{CR}^-(\h)$ in such a way that $\h$ contains $n^c_{sa}+1$ particles (see (\ref{defnsa}) for the definition of $n_{sa}^c$). Suppose that $\hbox{CR}^-(\h)$ contains $x$ empty frame-angles, with $1\leq x\leq4$. See Figure \ref{fig:figesempiostrong}(a) to visualize the configuration $\h$ in the case $x=1$. Since $H(\h)=\gs$, the move from $\h$ to $\o_k$ must have a non-positive cost and thus the unique admissible moves are:
	
	\begin{description}
		\item[(i)] either moving the free particle at zero cost;
		\item[(ii)] or removing the free particle;
		\item[(iii)] or attaching the free particle at cost $-U_1$ (see Figure \ref{fig:figesempiostrong}(b)) or $-U_2$, or $-U_1-U_2$.
	\end{description}
	
	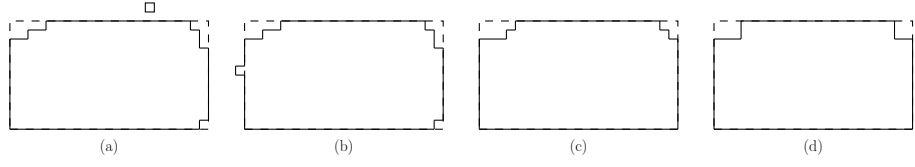
\begin{figure}
		\centering
		\begin{tikzpicture}[scale=0.24,transform shape]

			\draw (16,2.5)--(16,3);
			\draw (16,3) -- (16.5,3);
			\draw (16.5,3) -- (16.5,7);
			\draw (16.5,7) -- (16,7);
			\draw (16,7) -- (16,8);
			\draw (16,8)--(15.5,8);
			\draw(15.5,8)--(15.5,8.5);
			\draw (15.5,8.5)--(7.5,8.5);
			\draw (7.5,8.5)--(7.5,8);
			\draw [dashed] (5.5,2.5) rectangle (16.5,8.5);
			\draw (5.5,2.5) -- (16,2.5);
			\draw (5.5,2.5)--(5.5,7.5);
			\draw (5.5,7.5)--(6.5,7.5);
			\draw (7.5,8)--(6.5,8);
			\draw (6.5,8)--(6.5,7.5);
			\draw (13,9) rectangle (13.5,9.5);
			\node at (11,1.5){\Huge{(a)}};

			\draw (29,2.5)--(29,3);
			\draw (29,3) -- (29.5,3);
			\draw (29.5,3) -- (29.5,7);
			\draw (29.5,7) -- (29,7);
			\draw (29,7) -- (29,8);
			\draw (29,8)--(28.5,8);
			\draw(28.5,8)--(28.5,8.5);
			\draw (28.5,8.5)--(20.5,8.5);
			\draw (20.5,8.5)--(20.5,8);
			\draw [dashed] (18.5,2.5) rectangle (29.5,8.5);
			\draw (18.5,2.5) -- (29,2.5);
			\draw (18.5,2.5)--(18.5,5.5);
			\draw (18.5,5.5)--(18,5.5);
			\draw (18,5.5)--(18,6);
			\draw (18,6)--(18.5,6);
			\draw (18.5,6)--(18.5,7.5);
			\draw (18.5,7.5)--(19.5,7.5);
			\draw (20.5,8)--(19.5,8);
			\draw (19.5,8)--(19.5,7.5);
			\node at (24,1.5){\Huge{(b)}};

			\draw (42.5,2.5)--(42.5,3);
			\draw (42,2.5) -- (42.5,2.5);
			\draw (42.5,3) -- (42.5,7.5);
			\draw (42.5,7.5) -- (42,7.5);
			\draw (42,7.5) -- (42,8);
			\draw (42,8)--(41.5,8);
			\draw(41.5,8)--(41.5,8.5);
			\draw (41.5,8.5)--(33.5,8.5);
			\draw (33.5,8.5)--(33.5,8);
			\draw [dashed] (31.5,2.5) rectangle (42.5,8.5);
			\draw (31.5,2.5) -- (42,2.5);
			\draw (31.5,2.5)--(31.5,7.5);
			\draw (31.5,7.5)--(33,7.5);
			\draw (33.5,8)--(33,8);
			\draw (33,8)--(33,7.5);
			\node at (37,1.5){\Huge{(c)}};

			\draw (55.5,2.5)--(55.5,3);
			\draw (55,2.5) -- (55.5,2.5);
			\draw (55.5,3) -- (55.5,7.5);
			\draw (55.5,7.5) -- (54.5,7.5);
			\draw (54.5,7.5) -- (54.5,8);
			\draw(54.5,8)--(54.5,8.5);
			\draw (54.5,8.5)--(46,8.5);
			\draw (46,8.5)--(46,8);
			\draw [dashed] (44.5,2.5) rectangle (55.5,8.5);
			\draw (44.5,2.5) -- (55,2.5);
			\draw (44.5,2.5)--(44.5,7.5);
			\draw (44.5,7.5)--(46,7.5);
			\draw (46,8)--(46,7.5);
			\node at (50,1.5){\Huge{(d)}};

		\end{tikzpicture}
		
		\vskip 0 cm
		\caption{Here we depict in (a) the configuration $\h$; in (b) the configuration obtained by $\h$ by attaching the free particle at cost $-U_1$ to $B^w(\h)$; in (c) the configuration $\h'$ obtained from $\h$ by attaching the free particle to $c^{se}(\h)$ and then detaching the particle in $c^{nw}(\hbox{CR}^-(\h))$ and attach it to $B^e(\h)$, and in (d) the configuration $\h''$ obtained from $\h'$ by detaching the particle in $c^{ne}(\hbox{CR}^-(\h'))$ attach it to $B^n(\h')$.}
		\label{fig:figesempiostrong}
	\end{figure}

	\noindent
	{\bf Case (i).} In this case the configuration $\o_{k}$ is analogue to $\h$ and therefore we can iterate this argument by taking this configuration as $\h$.
	
	\medskip
	\noindent
	{\bf Case (ii).} In this case $H(\o_k)=\gs-\D$. We may assume that the configuration $\o_{k-1}$ is not obtained by $\o_{k}$ via adding a free particle, otherwise $\o_{k-1}$ is analogue to $\h$ and thus we can iterate the argument by taking this configuration as $\h$. By the optimality of the path, again considering the time-reversal, we deduce that the unique admissible move to obtain $\o_{k-1}$ from $\o_{k}$ is breaking a horizontal (resp.\ vertical) bar at cost $U_1$ (resp.\ $U_2$). Thus it is possible that either a sequence of $1$-translations of a bar or a sliding of a bar around a frame-angle takes place. In the first case, we obtain a configuration that is analogue to $\o_{k-1}$ and thus we can iterate the argument for a finite number of steps, since the path has to reach $\vuoto$. In the latter case, by Remark \ref{barrestrong}(ii) we deduce that the condition (\ref{condtrenino}) is not satisfied and therefore it is not possible to complete any sliding of a bar around a frame-angle. This implies that the unique admissible moves are the reverse ones, thus we obtain a configuration that is analogue to $\o_{k-1}$ and therefore we can iterate the argument for a finite number of steps, since the path has to reach $\vuoto$. In this way we can reduce ourselves to consider the case (iii).
	
	\medskip
	\noindent
	{\bf Case (iii).} (a) We consider the case where from $\h$, again considering the time-reversal, we attach a particle at cost $-U_1$ in $\partial^+\hbox{CR}(\h)$ giving rise to the configuration $\o_k$, i.e., $H(\o_k)=\gs-U_1$ (see Figure \ref{fig:figesempiostrong}(b)). Thus it is possible that either a sequence of $1$-translations of a bar or a sliding of a bar around a frame-angle takes place. In the first case, we obtain a configuration that is analogue to $\o_{k}$ and thus we can iterate the argument for a finite number of steps, since the path has to reach $\vuoto$. In the latter case, by Remark \ref{barrestrong}(ii) we deduce that the condition (\ref{condtrenino}) is not satisfied and therefore it is not possible to complete any sliding of a bar around a frame-angle. This implies that the unique admissible moves are the reverse ones, thus we obtain a configuration that is analogue to $\o_{k}$ and therefore we can iterate the argument for a finite number of steps, since the path has to reach $\vuoto$.
	
	(b) We consider the case where from $\h$, again considering the time-reversal, we attach a particle at cost $-U_2$ in $\partial^+\hbox{CR}(\h)$ giving rise to the configuration $\o_k$, i.e., $H(\o_k)=\gs-U_2$. We argue in a similar way as above.
	
	(c) We consider the case where from $\h$, again considering the time-reversal, we attach a particle at cost $-U_1-U_2$ in $\partial^-\hbox{CR}(\h)$ giving rise to the configuration $\o_k$, i.e., $H(\o_k)=\gs-U_1-U_2$. Thus it is possible either to have a sequence of $1$-translations of a bar, or to have a sliding of a bar around a frame-angle, or to detach a particle at cost $U_1+U_2$. In the first two possibilities, analogously to what has been discussed previously in (a) and (b), the unique admissible moves are the reverse ones and therefore we conclude as above. In the latter possibility, we have that either $\o_{k-1}$ is obtained from $\o_{k}$ by detaching a particle from a bar at cost $U_1+U_2$ or from a corner of $\h$ that is in $\hbox{CR}^-(\h)$. In the first case, the particle can be attached to an empty frame-angle of $\hbox{CR}^-(\h)$ and we can repeat these steps at most $x-1$ times (if $x\geq2$), that implies that there exists $\bar{k}<k-1$ such that $\o_{\bar{k}}$ is composed by the union of a free particle and a rectangle $\cR(2l_2^*-4,l_2^*-2)$ with four bars attached to its four sides in such a way that $\o_{\bar{k}}$ contains $n^c_{sa}+1$ particles, namely $\o_{\bar{k}}\in\csgeo$. In the second case, we may assume that the detached particle is attached to a bar in $\partial^-\hbox{CR}(\h)$ giving rise to a configuration $\h'$ (see Figure \ref{fig:figesempiostrong}(c)), otherwise we obtain a configuration that is analogue to $\h$. Starting from $\h'$, similarly we obtain $\h''$ (see Figure \ref{fig:figesempiostrong}(d)) if $\h'$ has a corner in $\hbox{CR}^-(\h')$. If this is the case, we can proceed in a similar way until we obtain a configuration $\h'''$ that has no corner in $\hbox{CR}^-(\h''')$. Starting from $\h'''$, by the optimality of the path we deduce that the unique admissible moves are the reverse ones and therefore the path goes back to $\h$. This concludes the proof. 
\end{proof*}

\begin{proof*}{\bf of Proposition \ref{p2}}
	Consider any $\o\in(\vuoto\ra\pieno)_{opt}$ such that $\o\cap\cP_{sa,0}\neq\emptyset$. We consider separately the following three cases.

	\noindent
	{\bf Case (i).} Assume that the path $\o$ crosses the set $\cA:=\bigcup_{\a,\a',k}\cA_{k}^{\a,\a'}$, with $\a\in\{n,s\}$, $\a'\in\{w,e\}$ and $k=2,...,l_2^*$, in the configuration $\h$. Since $H(\h)=\gs$, as long as the energy does not exceed $\gs$, it is impossible to create a free particle before further lowering the energy by a quantity greater or equal than $U_1+U_2$. Since $g_2'(\h)=1$, it is possible to connect at cost $-U_1$ the two protuberances or a bar and a protuberance. Moreover, there is no admissible move that costs $-U_2$, since $g_1'(\h)=0$ and there is no free particle that could been attached to an horizontal side of the cluster. Thus the only admissible moves are starting a sliding of a bar around a frame-angle $c^{\a\a'}(\h)$ or $c^{\a'\a}(\h)$, with $\a\in\{n,s\}$, $\a'\in\{w,e\}$, at cost less or equal than $U_1$. We consider separately these two possibilities, that correspond to the two different directions to cross the path described in Figure \ref{fig:columntorow} starting from the configuration $\h$. More precisely, one of these possibilities (that we will analyze in (iA)), gives rise to the configuration (12), while the other (that will be treated in (iB)) corresponds to the time-reversal of the path described in Figure \ref{fig:columntorow} starting from the configuration $\h$.

	{\bf (iA).} In this situation it is possible to obtain one or more saddles $\x_1,...,\x_{n-1}$ such that $\x_i\in\cA$ for all $i=1,...,n-1$ and for the last configuration we have $|r^{\a'}(\x_{n-1})\cup c^{\a\a'}(\x_{n-1})|=1$, with $\a\in\{n,s\}$, $\a'\in\{w,e\}$ and $g_2'(\x_{n-1})=1$ (see configuration (12) in Figure \ref{fig:columntorow}). From this configuration, since $H(\x_{n-1})=\gs$, by the optimality of the path $\o$ it is impossible to detach the protuberance before lowering the energy. Thus the unique admissible moves are either the reverse move or connect the protuberance and the bar at cost $-U_1$ and then detach the protuberance at cost $U_1$ (see the move starting from the configuration (12) in Figure \ref{fig:columntorow} that is described with a dashed arrow). In the latter situation the path reaches a configuration $\x_n\in\csgeo$. Thus we have to consider the possibilities that $\o$ visits $\x_{n}$ and $\o$ does not visit $\x_n$. In the first possibility, since $\o$ passes through $\x_n\in\csgeo$ we get the claim. In the latter possibility, the path $\o$ does not pass through the configuration $\x_n$, but assume that the path $\o$ visits the saddles $\x_i,...,\x_j$ for some $1\leq i\leq j\leq n-1$. We set $\o=(\vuoto,\o_1,...,\o_k,\x_i,\z_i,..,\z_h,\x_{i+1},\z_{h+1},...,\z_{h+m},...,\x_j)\circ\bar{\o}$, where $\z_i,...,\z_h$ are not saddles, but are crossed during the sliding of a bar around a frame-angle connecting $\x_i$ to $\x_{i+1}$ and so on. Moreover, $\bar\o$ is a path that connects $\x_j$ to $\pieno$ such that $\max_{\s\in\o}H(\s)\leq\gs$. Note that the configuration $\x_j$ coincides with $\x_{i+1}$ in the case $j=i+1$. If $i=j$, we set $\o=(\vuoto,\o_1,...,\o_k,\x_i)\circ\bar{\o}$, where $\bar\o$ is a path that connects $\x_i$ to $\pieno$ such that $\max_{\s\in\o}H(\s)\leq\gs$. To prove our statement we investigate the structure of the path $\o$ before entering $\cA$, namely we consider the time-reversal of the path. Since $\x_i\in\cA$ that implies $H(\x_i)=\gs$, we note that the move from $\x_i$ to $\o_k$ must have a non-positive cost. Thus the admissible transitions from $\x_i$ to $\o_k$ are either moving the particle at zero cost or moving a particle at cost $-U_1$. In the first case, we obtain a configuration that is analogue to $\x_i$ and therefore we can iterate for a finite number of steps this argument until we get the situation described in the latter case. In the latter case $H(\o_k)=\gs-U_1$, thus $\o_{k-1}$ can be obtained from $\o_{k}$ by breaking a bar at cost $U_1$ or $U_2$, since it is not possible to detach any particle because its cost is at least $U_1+U_2$. If the cost is $U_1$, we deduce that $\o_{k-1}$ is analogue to the initial configuration $\x_i$ and thus we can iterate this argument for a finite number of steps, because the path has to reach $\vuoto$. If the cost is $U_2$, we can iterate this argument to deduce that starting from $\o_k$ a sliding of a bar around a frame-angle takes place (see the time-reversal of the path described in Figure \ref{fig:columntorow}, where $\o_{k-1}$ can be, for example, the configuration (9)). Since the path has to reach $\vuoto$, this implies that there exists $k_1<k$ such that $\o_{k_1}$ is composed by the union of a rectangle $\cR(2l_2^*-1,l_2^*-1)$ and a protuberance attached to one of the longest sides. Note that $H(\o_{k_1})=\gs-U_1$. Furthermore, since either moving the protuberance along the side or detaching it from the cluster are the only admissible moves with $\max_{\s\in\o}H(\s)\leq\gs$, we note that there exists $k_2<k_1$ such that the configuration $\o_{k_2}$ is composed by the union of $\cR(2l_2^*-1,l_2^*-1)$ and a free particle. Note that $\o_{k_2}$ belongs to the set $\cB$ defined in \cite[Definition 3.5]{BN} because $p_2(\o_{k_2})=l_2^*-1$. Thus by \cite[Theorem 3.7]{BN} we know that $\o$ reaches a configuration in $\cP_2$. We get the claim by using Lemma \ref{gatestrong}.

	{\bf (iB).}  Note that this situation can be treated as in (iA) for the case in which $\o$ does not visit $\x_n$, indeed without loss of generality we may assume that $\h=\x_i$ and then we proceed as above.

	\noindent
	{\bf Case (ii).} Assume that $\o$ crosses the set $\cA_{0}^{\a'}\cup\cA_{1}^{\a}$ in the configuration $\h$, with $\a\in\{n,s\}$ and $\a'\in\{w,e\}$. By Lemma \ref{trasless} we know that $\h$ has been obtained from a configuration $\h^B\in\cC_{sa}^B$, thus there exists a configuration $\bar\h\in\cC_{sa}^*$ such that $\o$ passes through $\bar\h$ before crossing $\h^B$. 
	
	\noindent
	{\bf Case (iii).} Assume that $\o$ crosses the set $\cP_{sa,0}\setminus(\bigcup_{\a,\a',k}\cA_{k}^{\a,\a'}\cup\bigcup_{\a'}\cA_{0}^{\a'}\cup\bigcup_{\a}\cA_{1}^{\a})$ in the configuration $\h$, thus the path $\o$ crosses either $\h^B\in\cC_{sa}^B$ or $\h^G\in\cC_{sa}^G$ before passing through $\h$. In the first case, $\h$ is obtained either via a $1$-translation of a bar or via a sliding of a bar around a frame-angle that in both cases can not be completed because $\h$ is not in $\bigcup_{\a,\a',k}\cA_{k}^{\a,\a'}\cup\bigcup_{\a'}\cA_{0}^{\a'}\cup\bigcup_{\a}\cA_{1}^{\a}$. Therefore the path $\o$, before crossing $\h^B$, passes through a configuration $\bar\h\in\cC_{sa}^*$. In the latter case, we argue similarly. This concludes the proof.
\end{proof*}

\subsection{Proof of the main Theorem  \ref{gstrong}} 
\label{strong2}

In this Section we analyze the geometry of the set $\cG_{sa}(\vuoto,\pieno)$ (recall (\ref{defg})). In particular, we give the proof of the main Theorem \ref{gstrong} by giving in Proposition \ref{strongselless} the geometric characterization of the essential saddles of the third type that are not in $\csgeo$ and that are visited after crossing the set $\cC_{sa}^B$.

\bp{strongselless}
Any saddle $\x$ that is neither in $\cC_{sa}^*$, nor in the boundary of the cycles $\cC_{\pieno}^{\vuoto}(\G_{sa}^*)$ nor $\cC_{\vuoto}^{\pieno}(\G_{sa}^*-H(\pieno))$ such that $\t_{\x}\geq\t_{\cC_{sa}^B}$ can be essential or not. For those essential we obtain the following description:
\be{}
\ba{ll}
\cG_{sa}(\vuoto,\pieno)\cap(\cS_{sa}(\vuoto,\pieno)\setminus(\partial\cC_{\pieno}^{\vuoto}(\gs)\cup\partial\cC_{\vuoto}^{\pieno}(\gs-H(\pieno))\cup\csgeo))=\\
\qquad\qquad\qquad\qquad\qquad\qquad\quad\qquad\qquad\qquad\displaystyle\bigcup_{\a}\bigcup_{\a'}\bigcup_{k=2}^{l_2^*}\cA_k^{\a,\a'}\cup\displaystyle\bigcup_{\a'}\cA_0^{\a'}\cup\displaystyle\bigcup_{\a}\cA_1^{\a}
\ea
\ee
\ep

\br{}
In Proposition \ref{selle3bis} we have proved that the saddles $\x$ of type three in Section \ref{S6.4} that are not in $\bigcup_{k,\a,\a'}\cA_{k}^{\a,\a'}$ and such that $\t_{\x}<\t_{\cC_{sa}^B}$ are unessential. Note that we have not to study separately the essentiality of the saddles $\x\in\bigcup_{k,\a,\a'}\cA_{k}^{\a,\a'}$, since $\bigcup_{k,\a,\a'}\cA_{k}^{\a,\a'}$ is included in the essential saddles $\x$ of type three such that $\t_{\x}\geq\t_{\cC_{sa}^B}$, analyzed in Proposition \ref{strongselless}.
\er

\noindent
We postpone the proof of Proposition \ref{strongselless} after the proof of the main Theorem \ref{gstrong}.

\medskip

\begin{proof*}{\bf of Theorem \ref{gstrong}}
	By Corollary \ref{corstrategy} we know that the saddles of the first and second type, defined in \cite[Definition 3.2]{BN2} and \cite[Definition 3.4]{BN2} respectively, are unessential. By Propositions \ref{selle3bis} and \ref{strongselless} we have the characterization of the essential saddles of the third type in Section \ref{S6.4}. We use Proposition \ref{c*contenuto} to get the claim.
\end{proof*}

\begin{proof*}{\bf of Proposition \ref{strongselless}}
	Consider a configuration $\h\in\csgeo(2)$ such that $\h=(\hat\h,x)$, with $\hat\h\in\bar\cD_{sa}$ and $d(\hat\h,x)=2$. By Proposition \ref{card}(b), note that $\hat\h$ consists of an $(2l_2^*-4)\times(l_2^*-2)$ rectangle with four bars $B^\a(\h)$, with $\a\in\{n,s,w,e\}$, attached to its four sides satisfying
	\be{eq1}
	1\leq|B^{w}(\h)|,|B^e(\h)|\leq l_2^*, \qquad l_2^*-1\leq|B^{n}(\h)|,|B^s(\h)|\leq 2l_2^*-2, 
	\ee
	
	\noindent
	and
	\be{eq2} 
	\displaystyle\sum_\a |B^\a(\h)|-k=5l_2^*-7,
	\ee
	
	\noindent
	with $k=\sum_{\a\a'\in\{nw,ne,sw,se\}}|c^{\a\a'}(\h)|$. Assume that the free particle is attached in a bad site obtaining a configuration $\h'\in\cC_{sa}^B$. Due to \cite[Theorem 5.1]{MNOS}, our strategy consists in characterizing the essential saddles that could be visited after attaching the free particle in a bad site. By Remark \ref{barrestrong}(i) we consider separately the following cases:
	
	\begin{itemize}
		\item[A.] three frame-angles of CR$(\hat\h)$ are occupied;
		\item[B.] two frame-angles of CR$(\hat\h)$ are occupied;
		\item[C.] one frame-angle of CR$(\hat\h)$ is occupied;
		\item[D.] no frame-angle of CR$(\hat\h)$ is occupied.
	\end{itemize}
	
	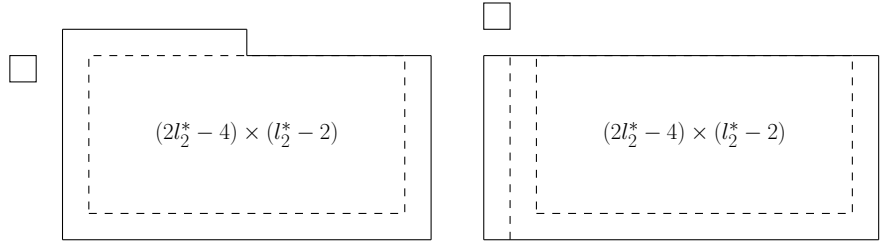
\begin{figure}
		\centering
		\begin{tikzpicture}[scale=0.35,transform shape]

			\draw [dashed] (6,3) rectangle (18,9);
			\node at (12,6) {\Huge{$(2l_2^*-4)\times(l_2^*-2)$}};
			\draw(12,9)--(12,10);
			\draw(12,10)--(5,10);
			\draw(5,10)--(5,2);
			\draw(5,2)--(19,2);
			\draw(19,2)--(19,9);
			\draw(19,9)--(18,9);
			\draw(18,9)--(12,9);
			\draw(3,8) rectangle (4,9);

			\draw [dashed] (23,3) rectangle (35,9);
			\node at (29,6) {\Huge{$(2l_2^*-4)\times(l_2^*-2)$}};
			\draw(22,2)--(36,2);
			\draw(36,2)--(36,9);
			\draw(36,9)--(35,9);
			\draw(35,9)--(21,9);
			\draw(21,9)--(21,2);
			\draw(21,2)--(22,2);
			\draw[dashed] (22,2)--(22,9);
			\draw(21,10) rectangle (22,11);

		\end{tikzpicture}
		
		\vskip 0 cm
		\caption{Case A: on the left-hand side we represent a possible starting configuration $\h\in\csgeo$ and on the right-hand side the configuration $\h'$ obtained from $\h$ after the sliding of the bar $B^n(\h)$ around the frame-angle $c^{nw}(\h')$.}
		\label{fig:fig14}
	\end{figure}

	Note that from case A one can go to the other cases and viceversa, but since the path has to reach $\pieno$ this back and forth must end in a finite number of steps.
	
	\noindent
	{\bf Case A.} Without loss of generality we consider $\h$ as in Figure \ref{fig:fig14} on the left-hand side. If we are considering the case in which a sequence of $1$-translations of a bar is possible and takes place, then by Lemma \ref{trasless} the saddles that are crossed are essential and in $\cA_0^{\a'}\cup\cA_{1}^{\a}$. If a sequence of $1$-translations of a bar takes place in such a way that the last configuration has at most two occupied frame-angles, then the saddles that are visited starting from it will be analyzed in cases B, C and D. Thus we are left to analyze the case in which there is the activation of a sliding of a bar around a frame-angle. In the following we quickly exclude the cases in which the particles is attached to $B^n(\h)$, $B^s(\h)$ or $B^e(\h)$ and then explain the more interesting case in which it is attached to $B^w(\h)$ giving rise to Figure \ref{fig:fig14} on the right-hand side. If the free particle is attached to the bar $B^n(\h)$ (resp.\ $B^s(\h)$), by Lemma \ref{coltorow} we know that it is not possible to complete the sliding of the bar $B^w(\h)$ (resp.\ $B^e(\h)$) around the frame-angle $c^{wn}(\h')$ (resp.\ $c^{es}(\h')$). If the free particle is attached to the bar $B^e(\h)$ or $B^w(\h)$, then it is not possible to slide the bar $B^s(\h)$ around the frame-angle $c^{se}(\h')$ or $c^{sw}(\h')$ respectively, since (\ref{condtrenino}) is not satisifed. In the last two cases by Lemma \ref{trenini}(ii) we know that the saddles that are visited are unessential. This implies that the unique possibility to activate and complete a sliding of a bar around a frame-angle is attaching the free particle to the bar $B^w(\h)$, then sliding the bar $B^n(\h)$ around the frame-angle $c^{nw}(\h')$ when $|B^n(\h)|=l_2^*-1$ and $|B^w(\h)|=l_2^*$, otherwise (\ref{condtrenino}) is not satisifed. The saddles that are possibly visited by the path we described are in $\cA_{k,k'}^{\a,\a'}$ except the last one, thus by Lemma \ref{trenini}(i) they are essential. The last configuration visited during this sliding of a bar is depicted in Figure \ref{fig:fig14} on the right-hand side. This configuration has energy $\gs-U_1+U_2$ and therefore it is not a saddle. Starting from this configuration, by Lemma \ref{entratastrong} we know that the saddles that could be visited are in $\csgeo$ or again in $\cA_{k,k'}^{\a,\a'}$. This concludes case A.

	\begin{figure}
		\centering
		\begin{tikzpicture}[scale=0.23,transform shape]

			\draw [dashed] (6,3) rectangle (18,9);
			\node at (12,6) {\Huge{$(2l_2^*-4)\times(l_2^*-2)$}};
			\draw(16,9)--(16,10);
			\draw(16,10)--(5,10);
			\draw(5,10)--(5,4);
			\draw(5,4)--(6,4);
			\draw(6,4)--(6,3);
			\draw(6,3)--(7,3);
			\draw(7,3)--(7,2);
			\draw(7,2)--(19,2);
			\draw(19,2)--(19,8);
			\draw(19,8)--(18,8);
			\draw(18,8)--(18,9);
			\draw(18,9)--(16,9);
			\draw(3,8) rectangle (4,9);
			\node at (12,0.5) {\Huge{(a)}};

			\draw [dashed] (23,3) rectangle (35,9);
			\node at (29,6) {\Huge{$(2l_2^*-4)\times(l_2^*-2)$}};
			\draw(23,3)--(24,3);
			\draw(24,3)--(24,2);
			\draw(24,2)--(33,2);
			\draw(33,2)--(33,3);
			\draw(33,3)--(35,3);
			\draw(35,3)--(35,4);
			\draw (35,4)--(36,4);
			\draw(36,4)--(36,10);
			\draw(36,10)--(22,10);
			\draw(22,10)--(22,4);
			\draw(22,4)--(23,4);
			\draw(23,4)--(23,3);
			\draw (23,11) rectangle (24,12);
			\node at (29,0.5) {\Huge{(b)}};

			\draw [dashed] (40,3) rectangle (52,9);
			\node at (46,6) {\Huge{$(2l_2^*-4)\times(l_2^*-2)$}};
			\draw(47,9)--(47,10);
			\draw(47,10)--(39,10);
			\draw(39,10)--(39,2);
			\draw(39,2)--(52,2);
			\draw(52,2)--(52,3);
			\draw(52,3)--(53,3);
			\draw(53,3)--(53,9);
			\draw(53,9)--(53,9);
			\draw(53,9)--(47,9);
			\draw(40,11) rectangle (41,12);
			\node at (46,0.5) {\Huge{(c)}};
			
		\end{tikzpicture}
		
		\vskip -0.2 cm
		\caption{Case B(i): in (a) we depict a possible starting configuration $\h\in\csgeo$. Case B(ii): in (b) we depict a possible starting configuration $\h\in\csgeo$. Case B(iii): in (c) we depict a possible starting configuration $\h\in\csgeo$.}
		\label{fig:fig5}
	\end{figure}
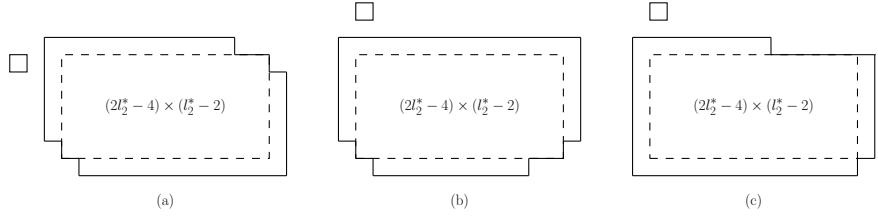

	\noindent 
	{\bf Case B.} If we are considering the case in which a sequence of $1$-translations of a bar is possible and takes place, then by  Lemma \ref{trasless} the saddles that are crossed are essential and in $\cA_0^{\a'}\cup\cA_{1}^{\a}$. We consider separately the following subcases:
	
	\begin{description}
		\item[(i)] The two occupied frame-angles are $c^{\a\a'}(\h)$ and $c^{\a''\a'''}(\h)$, with all the indeces $\a,\a',\a''$ and $\a'''$ different between each other (see Figure \ref{fig:fig5}(a));
		\item[(ii)] The two occupied frame-angles are $c^{\a\a'}(\h)$ and $c^{\a'\a''}(\h)$, with $\a'\in\{n,s\}$ and $\a\neq\a''$ (see Figure \ref{fig:fig5}(b));
		\item[(iii)] The two occupied frame-angles are $c^{\a\a'}(\h)$ and $c^{\a'\a''}(\h)$, with $\a'\in\{e,w\}$ and $\a\neq\a''$ (see Figure \ref{fig:fig5}(c)).
	\end{description}

	{\bf Case B(i).} Without loss of generality we consider $\h$ as in Figure \ref{fig:fig5}(a). We can reduce our proof to the case in which there is no translation of a bar and therefore there is the activation of a sliding of a bar around a frame-angle. If the free particle is attached to the bar $B^n(\h)$ (resp.\ $B^s(\h)$), by Lemma \ref{coltorow} we know that it is not possible to complete the sliding of the bar $B^w(\h)$ (resp.\ $B^e(\h)$) around the frame-angle $c^{wn}(\h')$ (resp.\ $c^{es}(\h')$). By Lemma \ref{trenini}(ii), this implies that the saddles that could be crossed are unessential. If the free particle is attached to the bar $B^w(\h)$ (resp.\ $B^e(\h)$), it is possible to slide the bar $B^n(\h)$ (resp.\ $B^s(\h)$) around the frame-angle $c^{nw}(\h')$ (resp.\ $c^{se}(\h')$) when $|B^n(\h)|<|B^w(\h)|$ (resp. $|B^s(\h)|<|B^e(\h)|$), otherwise (\ref{condtrenino}) is not satisifed. By (\ref{eq1}) and (\ref{eq2}) we note that $|B^n(\h)|<|B^w(\h)|$ (resp. $|B^s(\h)|<|B^e(\h)|$) is not possible and the case B(i) is concluded.

	{\bf Case B(ii).} Without loss of generality we consider $\h$ as in Figure \ref{fig:fig5}(b). If one bar among $B^w(\h)$ and $B^e(\h)$ is full, it is possible to translate $B^s(\h)$ in order to have three occupied frame-angles. This situation has already been analyzed in case A. Thus we can reduce our proof to the case in which there is no translation of a bar and therefore there is the activation of a sliding of a bar around a frame-angle. If the free particle is attached to the bar $B^n(\h)$ (or $B^s(\h)$), by Lemma \ref{coltorow} we know that it is not possible to complete the sliding of a vertical bar around any frame-angle. If the free particle is attached to the bar $B^w(\h)$ or $B^e(\h)$, since the bar $B^n(\h)$ is full, we deduce that (\ref{condtrenino}) is not satisfied. This implies that it is not possible to slide the bar $B^n(\h)$ around the frame-angle $c^{nw}(\h')$ and $c^{ne}(\h')$. In the last two cases by Lemma \ref{trenini}(ii) we know that the saddles that could be visited are unessential. This concludes case B(ii).

	{\bf Case B(iii).} Without loss of generality we consider $\h$ as in Figure \ref{fig:fig5}(c). If the bar $B^n(\h)$ (or $B^s(\h)$) is full, it is possible to translate $B^e(\h)$ to occupy the frame-angle $c^{ne}(\h')$ (or $c^{se}(\h')$). This situation has already been analyzed in case A. Otherwise, it is possible to translate a bar with one occupied frame-angle in order to have two occupied frame-angles in such a way that they have no bar in common. This situation has already been analyzed in case B(i). Thus we can reduce our proof to the case in which there is no translation of a bar and therefore we can consider only the activation of a sliding of a bar around a frame-angle. If the free particle is attached to the bar $B^n(\h)$ (resp.\ $B^s(\h)$), by Lemma \ref{coltorow} we know that it is not possible to complete the sliding of the bar $B^w(\h)$ around the frame-angle $c^{wn}(\h')$ (resp.\ $c^{ws}(\h')$). If the free particle is attached to the bar $B^e(\h)$, we deduce that (\ref{condcorner}) is not satisfied. In the last two cases by Lemma \ref{trenini}(ii) we know that the saddles that are visited are unessential. If the free particle is attached to the bar $B^w(\h)$, it is possible to slide the bar $B^n(\h)$ (resp.\ $B^s(\h)$) around the frame-angle $c^{nw}(\h')$ (resp.\ $c^{sw}(\h')$) when $|B^n(\h)|<|B^w(\h)|$ (resp.\  $|B^s(\h)|<|B^w(\h)|$), otherwise (\ref{condtrenino}) is not satisifed. By (\ref{eq1}) and (\ref{eq2}) we note that $|B^n(\h)|<|B^w(\h)|$ (resp.\  $|B^s(\h)|<|B^w(\h)|$) is not possible and the case B(iii) is concluded.
	
	\begin{figure}
		\centering
		\begin{tikzpicture}[scale=0.32,transform shape]

			\draw [dashed] (14,3) rectangle (26,9);
			\node at (20,6) {\Huge{$(2l_2^*-4)\times(l_2^*-2)$}};
			\draw(24,9)--(24,10);
			\draw(24,10)--(13,10);
			\draw(13,10)--(13,4);
			\draw(13,4)--(14,4);
			\draw(14,4)--(14,3);
			\draw(14,3)--(15,3);
			\draw(15,3)--(15,2);
			\draw(15,2)--(26,2);
			\draw(26,2)--(26,3);
			\draw(26,3)--(27,3);
			\draw(27,3)--(27,9);
			\draw(27,9)--(26,9);
			\draw(26,9)--(24,9);
			\draw(11,8) rectangle (12,9);

			\draw [dashed] (32,3) rectangle (44,9);
			\node at (38,6) {\Huge{$(2l_2^*-4)\times(l_2^*-2)$}};
			\draw(32,3)--(33,3);
			\draw(33,3)--(33,2);
			\draw(33,2)--(44,2);
			\draw(44,2)--(44,3);
			\draw(44,3)--(44,4);
			\draw (44,4)--(45,4);
			\draw(45,4)--(45,9);
			\draw(45,9)--(44,9);
			\draw(44,9)--(43,9);
			\draw(43,9)--(43,10);
			\draw(43,10)--(32,10);
			\draw(32,10)--(32,9);
			\draw(32,9)--(31,9);
			\draw(31,9)--(31,3);
			\draw(31,3)--(32,3);
			\draw(32,4)--(32,3);
			\draw (32,11) rectangle (33,12);

		\end{tikzpicture}
		
		\vskip 0. cm
		\caption{Case C: on the left-hand side we depict a possible starting configuration $\h\in\csgeo$. Case D: on the right-hand side we depict a possible starting configuration $\h\in\csgeo$.}
		\label{fig:fig7}
	\end{figure}
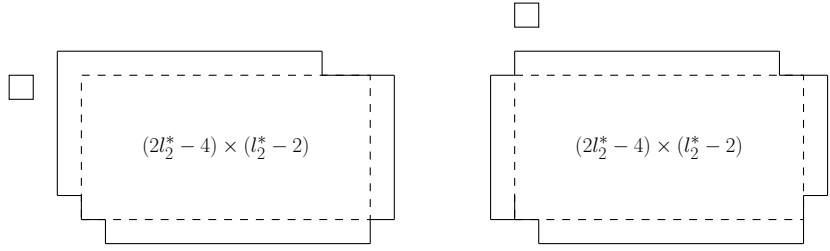

	\noindent
	{\bf Case C.} Without loss of generality we consider $\h$ as in Figure \ref{fig:fig7} on the left-hand side. If we are considering the case in which a sequence of $1$-translations of a bar is possible and takes place, then by Lemma \ref{trasless} the saddles that are crossed are essential and in $\cA_0^{\a'}\cup\cA_{1}^{\a}$. Starting from this configuration it is possible to obtain two occupied frame-angles, that it has been already analyzed in Case B. Thus we can reduce our proof to the case in which there is no translation of a bar and therefore there is the activation of a sliding of a bar around a frame-angle. If the free particle is attached to the bar $B^n(\h)$ (resp.\ $B^s(\h)$), by Lemma \ref{coltorow} we know that it is not possible to complete the sliding of the bar $B^w(\h)$ around the frame-angle $c^{wn}(\h')$ (resp.\ $c^{ws}(\h')$). If the free particle is attached to the bar $B^e(\h)$, we deduce that (\ref{condtrenino}) is not satisfied. In the last two cases by Lemma \ref{trenini}(ii) we know that the saddles that are visited are unessential. If the free particle is attached to the bar $B^w(\h)$, it is possible to slide the bar $B^n(\h)$ around the frame-angle $c^{nw}(\h')$ when $|B^n(\h)|<|B^w(\h)|$, otherwise (\ref{condtrenino}) is not satisfied. By (\ref{eq1}) and (\ref{eq2}) we note that $|B^n(\h)|<|B^w(\h)|$ is not possible and the case C is concluded.
	
	\noindent
	{\bf Case D.} Without loss of generality we consider $\h$ as in Figure \ref{fig:fig7} on the right-hand side. If we are considering the case in which a sequence of $1$-translations of a bar is possible and takes place, then by Lemma \ref{trasless} the saddles that are crossed are essential and in $\cA_0^{\a'}\cup\cA_{1}^{\a}$. Starting from this configuration, it is possible to obtain one or two occupied frame-angles: these situations have been already analyzed in cases C and B respectively. Thus we can reduce our proof to the case in which there is no translation of a bar and therefore there is the activation of a sliding of a bar around a frame-angle. If the free particle is attached to the bar $B^n(\h)$ (resp.\ $B^s(\h)$), by Lemma \ref{coltorow} we know that it is not possible to complete the sliding of the bar $B^w(\h)$ around the frame-angle $c^{wn}(\h')$ (resp.\ $c^{ws}(\h')$). If the free particle is attached to the bar $B^w(\h)$ or $B^e(\h)$, we deduce that (\ref{condcorner}) is not satisfied. In the last two cases by Lemma \ref{trenini}(ii) we know that the saddles that are visited are unessential. This concludes case D.
\end{proof*}

\section{Proof of the sharp asymptotics}
\label{sharpasymptotics}
For the model-independent discussion we refer to \cite[Section 10.1]{BN2}. Following the strategy given in \cite{BHN} for the isotropic case, here we apply this argument for the strongly anisotropic one. For the corresponding strategy in the isotropic and weakly anisotropic cases we refer to \cite[Section 10.2]{BN2}.

\subsection{Application of the potential theory to the strongly anisotropic case}
\label{kawstrong}
In \cite{BH} the authors let the protocritical and critical sets as $\mathscr{P}^*(m,s)$ and $\mathscr{C}^*(m,s)$ respectively (see \cite[Definition 16.3]{BH} for the definition of $\mathscr{P}^*(m,s)$ and $\mathscr{C}^*(m,s)$). Since they differ from our notation, we refer to them as $\mathscr{P}_{PTA}^*(m,s)$ and $\mathscr{C}_{PTA}^*(m,s)$. In \cite{BH} the authors proved \cite[Theorem 16.4]{BH} and \cite[Theorem 16.5]{BH} subject to the two hypotheses
\begin{description}
	\item[(H1)] $\cX^{m}=\{m\}$ and $\cX^{s}=\{s\}$;
	\item[(H2)] $\x'\ra|\{\x\in\mathscr{P}_{PTA}^*(m,s): \ \x\sim\x'\}|$ is constant on $\mathscr{C}_{PTA}^*(m,s)$.
\end{description}
For our model $\cX^{m}_{sa}=\{\vuoto\}$ and $\cX^{s}_{sa}=\{\pieno\}$, thus (H1) holds and $\G^*=\Phi(\vuoto,\pieno)-H(\vuoto)=\G^*_{sa}$. Now we abbreviate $\mathscr{P}_{PTA}^*=\mathscr{P}_{PTA}^*(\vuoto,\pieno)$ and  $\mathscr{C}_{PTA}^*=\mathscr{C}_{PTA}^*(\vuoto,\pieno)$. Moreover, we prove that geometrically $\mathscr{P}_{PTA}^*=\bar\cD_{sa}\cup\bigcup_{\a,\a'}\bar\cA_2^{\a,\a'}$ (see (\ref{definizionebarA}) for the definition of $\bar\cA_2^{\a,\a'}$) and $\mathscr{C}_{PTA}^*=\cC_{sa}^*(L_{sa}^*)\cup\bigcup_{\a,\a'}\cA_2^{\a,\a'}$ (recall (\ref{defpsatrenino}) for the definition of $\cA_2^{\a,\a'}$) with $\a\in\{n,s\}$ and $\a'\in\{e,w\}$. Therefore it is clear that $\cC_{sa}^*\neq\mathscr{C}_{PTA}^*$.  Note that (H2) follows from Lemma \ref{entratastrong}, indeed each configuration in $\mathscr{C}_{PTA}^*$ has exactly one configuration in $\mathscr{P}_{PTA}^*$ from which it can be reached via an allowed move. In particular, the configurations in $\cC_{sa}^*(L_{sa}^*)$ and $\bar\cD_{sa}$ are connected by removing the free particle in $\partial^-\L$, while those in $\cA_2^{\a,\a'}$ and $\bar\cA_{2}^{\a,\a'}$ are connected between each other by attaching the two particles separated by an empty site at cost $-U_1$. Since (H1) and (H2) hold, \cite[Theorem 16.4]{BH} and \cite[Theorem 16.5]{BH} should hold, but for the strongly anisotropic case this is not true. More precisely, this model represents a counterexample of \cite[Theorem 16.4(b)]{BH}, indeed on the one hand \cite[Theorem 16.4(a)]{BH} and \cite[Theorem  16.5]{BH} are valid, but on the other hand \cite[Theorem 16.4(b)]{BH} does not hold. This relies on a peculiar feature of this model: the entrance in $\mathscr{C}_{PTA}^*$ can not be uniform due to the two possibile different entrance mechanisms, as claimed in Lemma \ref{entratastrong}. This depends on the hypothesis (H2), that takes into account only the map from $\mathscr{C}_{PTA}^*$ to $\mathscr{P}_{PTA}^*$ and not the reverse one. Therefore we propose to replace the hypothesis (H2) with
\begin{description}
	\item[(H2')] $\x'\ra|\{\x\in\mathscr{P}_{PTA}^*(m,s): \x\sim\x'\}|$ is constant on $\mathscr{C}_{PTA}^*(m,s)$ and $\x\ra|\{\x'\in\mathscr{C}_{PTA}^*(m,s):\x'\sim\x\}|$ is constant on $\mathscr{P}_{PTA}^*(m,s)$.
\end{description}

\noindent
We are convinced that this could be the correct hypotheses, indeed the analysis of the uniform entrance distribution in $\mathscr{C}_{PTA}^*(m,s)$ has to take into account the number of configurations in $\mathscr{P}_{PTA}^*(m,s)$ that communicate with $\mathscr{C}_{PTA}^*(m,s)$ via one step of the dynamics. Now it is clear that this model does not satisfy (H2'), indeed each configuration in $\bar\cD_{sa}$ has exactly $4L-4$ configurations in $\cC_{sa}^*$ from which it can be reached via an allowed move, while each configuration in $\bar\cA_2^{\a,\a'}$ has only one configuration in $\cA_2^{\a,\a'}$ with this property. Therefore \cite[Theorem 16.4(b)]{BH} does not hold for this model.

Recall \cite[Definition 10.1]{BN2} for the definition of the wells $\cZ_{sa,j}^{\vuoto}$ and $\cZ_{sa,j}^{\pieno}$ and \cite[Definition 3.2]{BN2} and \cite[Definition 3.4]{BN2} for the definition of the saddles $\s_{sa,j}$ of the first type and $\z_{sa,j}$ of the second type respectively. Concerning \cite[Theorem 16.5]{BH}, by \cite[Lemma 16.16]{BH} for the case $int=sa$, we know that $h$ is constant on each wells. For the wells $\cZ_j^{m}$ and $\cZ_j^{s}$ this constant is computed in \cite[Lemma 10.4]{BN2}, indeed \cite[Lemma 16.15]{BH} can be extended for these sets together with the unessential saddles of the first and second type. Thanks to the model-independent discussion given in \cite[Section 10.1]{BN2} and \cite[Lemma 10.4]{BN2}, \cite[eq.\ (10.7)]{BN2} becomes
\be{}
h=
\begin{cases}
	1 &\hbox{on }  \displaystyle\cC_{\pieno}^{\vuoto}(\G^*_{sa})\cup\bigcup_{j=1}^{J_{\vuoto}}(\{\s_{sa,j}\}\cup\cZ_{sa,j}^{\vuoto}), \\
	0 &\hbox{on } \displaystyle\cC^{\pieno}_{\vuoto}(\G^*_{sa}-H(\pieno))\cup\bigcup_{j=1}^{J_{\pieno}}(\{\z_{sa,j}\}\cup\cZ_{sa,j}^{\pieno}),\\
	c_i &\hbox{on } \cX_{sa}(i), i=1,...,\bar{I},
\end{cases}
\ee

\noindent
where $\cX_{sa}(i)$, $i=1,...,\bar{I}$, are all the wells of the transition except $\bigcup_{j=1}^{J_{\vuoto}}\cZ_{sa,j}^{\vuoto}$ and $\bigcup_{j=1}^{J_{\pieno}}\cZ_{sa,j}^{\pieno}$. This implies that the unessential saddles, not characterizing the typical behavior of the process, can not be neglected in the study of the prefactor $K$. However, since they do not communicate with some $\cX_{sa}(i)$ via one step of the dynamics together with the fact that $h(\s_{sa,j})=1$ and $h(\z_{sa,j})=0$ for any $j$, the transitions that involve these unessential saddles do not contribute numerically to the computation of $K$. The variational formula for $\Theta=1/K$ in \cite[eq.\ (10.11)]{BN2} is non-trivial because it depends on the geometry of all the wells and on the form of the function $h$ on the configurations in $\cX_{sa}^{*}\setminus\cX_{sa}^{**}$, namely the saddle configurations.

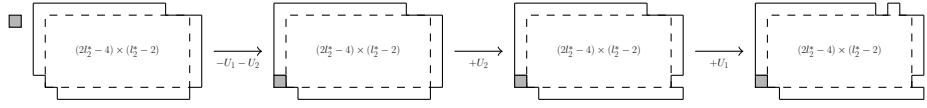
\begin{figure}
	\centering
	\begin{tikzpicture}[scale=0.16,transform shape]

		\draw [dashed] (14,3) rectangle (26,9);
		\node at (20,6) {\Huge{$(2l_2^*-4)\times(l_2^*-2)$}};
		\draw(24,9)--(24,10);
		\draw(24,10)--(13,10);
		\draw(13,10)--(13,4);
		\draw(13,4)--(14,4);
		\draw(14,4)--(14,3);
		\draw(14,3)--(15,3);
		\draw(15,3)--(15,2);
		\draw(15,2)--(26,2);
		\draw(26,2)--(26,3);
		\draw(26,3)--(27,3);
		\draw(27,3)--(27,9);
		\draw(27,9)--(26,9);
		\draw(26,9)--(24,9);
		\draw [fill=grigio2] (11,8) rectangle (12,9);

		\draw[->] (28,6)--(32,6);
		\node at (30,5){\Huge{{$-U_1-U_2$}}};
		\draw [dashed] (34,3) rectangle (46,9);
		\node at (40,6) {\Huge{$(2l_2^*-4)\times(l_2^*-2)$}};
		\draw(44,9)--(44,10);
		\draw(44,10)--(33,10);
		\draw(33,10)--(33,4);
		\draw(33,4)--(34,4);
		\draw(34,4)--(34,3);
		\draw(34,3)--(35,3);
		\draw(35,3)--(35,2);
		\draw(35,2)--(46,2);
		\draw(46,2)--(46,3);
		\draw(46,3)--(47,3);
		\draw(47,3)--(47,9);
		\draw(47,9)--(46,9);
		\draw(46,9)--(44,9);
		\draw [fill=grigio2] (33,3) rectangle (34,4);

		\draw[->] (48,6)--(52,6);
		\node at (50,5){\Huge{{$+U_2$}}};
		\draw [dashed] (54,3) rectangle (66,9);
		\node at (60,6) {\Huge{$(2l_2^*-4)\times(l_2^*-2)$}};
		\draw(64,9)--(64,10);
		\draw(64,10)--(53,10);
		\draw(53,10)--(53,4);
		\draw(53,4)--(54,4);
		\draw(54,4)--(54,3);
		\draw(54,3)--(55,3);
		\draw(55,3)--(55,2);
		\draw(55,2)--(67,2);
		\draw(67,2)--(67,3);
		\draw (67,3)--(66,3);
		\draw (66,3)--(66,4);
		\draw(66,4)--(67,4);
		\draw(67,4)--(67,9);
		\draw(67,9)--(66,9);
		\draw(66,9)--(64,9);
		\draw [fill=grigio2] (53,3) rectangle (54,4);

		\draw[->] (68,6)--(72,6);
		\node at (70,5){\Huge{{$+U_1$}}};
		\draw [dashed] (74,3) rectangle (86,9);
		\node at (80,6) {\Huge{$(2l_2^*-4)\times(l_2^*-2)$}};
		\draw(83,9)--(84,9);
		
		\draw(83,9)--(83,10);
		\draw(83,10)--(73,10);
		\draw(73,10)--(73,4);
		\draw(73,4)--(74,4);
		\draw(74,4)--(74,3);
		\draw(74,3)--(75,3);
		\draw(75,3)--(75,2);
		\draw(75,2)--(87,2);
		\draw(87,2)--(87,3);
		\draw (87,3)--(86,3);
		\draw (86,3)--(86,4);
		\draw(86,4)--(87,4);
		\draw(87,4)--(87,9);
		\draw(87,9)--(86,9);
		\draw(86,9)--(85,9);
		\draw (85,9)--(85,10);
		\draw(85,10)--(84,10);
		\draw (84,10)--(84,9);
		\draw [fill=grigio2] (73,3) rectangle (74,4);

	\end{tikzpicture}
	\vskip 0. cm
	\caption{We depict the transition that, starting from a configuration in $\cC_{sa}^*$, passes trough $\cC_{sa}^G$ and after two moves reaches an unessential saddle $\zeta_{sa}$.}
	\label{fig:sellaG}
\end{figure}

\br{}
In Figure \ref{fig:sellaG} we depict a transition that, starting from a configuration in $\cC_{sa}^*$, gives an unessential saddle $\z_{sa}$ of the third type. 
\er

\subsection{Proof of Theorem \ref{sharptimestrong}}
Using \cite[Lemma 10.4]{BN2}, in order to prove Theorem \ref{sharptimestrong} it remains to analyze in detail the number of possible transitions inside the gates and in-between their boundaries. Finally, we need to count the cardinality of $\bar\cD_{sa}$ modulo shifts. We denote by $N_{sa}$ that quantity.

\subsubsection{Lower bound}
Recall \cite[Definition 10.1]{BN2} for the definition of the wells $\cZ_{sa,j}^{\vuoto}$ and $\cZ_{sa,j}^{\pieno}$ and \cite[Definition 3.2]{BN2} and \cite[Definition 3.4]{BN2} for the definition of the saddles $\s_{sa,j}$ of the first type and $\z_{sa,j}$ of the second type respectively.  Thus we consider the following sets:

\be{XA}
\cX_{sa}^{\vuoto}:=\cC_{\pieno}^{\vuoto}(\gs)\cup\bigcup_{j=1}^{J_{\vuoto}}(\{\s_{sa,j}\}\cup\cZ_{sa,j}^{\vuoto})
\ee
\be{XB}
\cX_{sa}^{\pieno}:=\cC_{\vuoto}^{\pieno}(\gs-H(\pieno))\cup\bigcup_{j=1}^{J_{\pieno}}(\{\z_{sa,j}\}\cup\cZ_{sa,j}^{\pieno})
\ee

\noindent
The lower bound $\Theta_{sa}\geq\Theta_1^{sa}$ is obtained by removing all the transitions that do not involve either a protocritical droplet and a free particle that moves or the path described in Figure \ref{fig:columntorow}. The first type of transitions gives a contribute that can be treated in a similar way as in the lower bound in \cite[Proposition 3.3.4]{BHN}. Indeed we obtain:

\be{caprw}
\ba{lll}
\displaystyle\sum_{\hat\h\in\bar\cD_{sa}}\min_{c_j(\hat\h),j=1,2,3,4}\min_{g:\Lambda^*\ra[0,1] \atop g_{|\partial^G\hat\h} \equiv 0, g_{|\partial_j^B\hat\h}\equiv c_j, j=1,2,3,4, g_{|\partial^+\Lambda}\equiv1}\frac{1}{2}\sum_{x,x'\in\Lambda^+ \atop x\sim x'}[g(x)-g(x')]^2 \\
\qquad  \geq \displaystyle\sum_{\hat\h\in\bar\cD_{sa}} \CAPA^{\Lambda^+}(\partial^+\Lambda,\hbox{CR}(\hat\h))
\ea
\ee

\noindent
where $g(x):=h(\hat\h,x)=h(\h)$ for $\hat\h\in\bar\cD_{sa}$ and $x\in\L\setminus\hbox{CR}^{++}(\hat\h)$, and $\partial^G\hat\h$ denotes the set of good sites in $\partial^-\hbox{CR}(\hat\h)$, $\partial_j^B\hat\h$, $j=1,2,3,4$, denote the four bars of bad sites in $\partial^+\hbox{CR}(\hat\h)$ and

\be{cap+}
\CAPA^{\Lambda^+}(\partial^+\Lambda,F)=\min_{g:\Lambda^+\ra[0,1] \atop g_{|\partial^+\Lambda}\equiv1, g_{|F}\equiv0} \frac{1}{2}\sum_{x,x'\in\Lambda^+ \atop x\sim x'}[g(x)-g(x')]^2, \quad \hbox{for any } F\subseteq\Lambda^+.
\ee

\noindent
For $\a\in\{n,s\}$, $\a'\in\{w,e\}$ and $k=2,...,l_2^*$, we define
\be{definizionebarA}
\ba{ll}
\bar{\cA}_k^{\a,\a'}&:=\{\h: n(\h)=0, v(\h)=2l_2^*-2, |r^\a(\h)|=k-1, |r^\a(\h)|=l_2^*-k, |c^{\a\a'}(\h)|=1, \\
&\quad \ \ \h_{cl} \hbox{ is } \hbox{connected, monotone}, \hbox{with circumscribed rectangle in } \cR(2l_2^*-1,l_2^*)\}
\ea
\ee

\noindent
and
\be{}
\ba{ll}
\widetilde{\cA}_k^{\a,\a'}&:=\{\h: n(\h)=0, v(\h)=2l_2^*-2, |r^\a(\h)|=k, |r^\a(\h)|=l_2^*-k, |c^{\a\a'}(\h)|=0,\\
&\quad \ \ \h_{cl} \hbox{ is } \hbox{connected, monotone}, \hbox{with circumscribed rectangle in } \cR(2l_2^*-1,l_2^*)\}
\ea
\ee

\noindent
Referring to Figure \ref{fig:columntorow}, note that the configuration (6) is in $\bar{\cA}_2^{n,e}$ and the configuration (10) is in $\bar{\cA}_3^{n,e}$, while the configuration (8) is in $\widetilde{\cA}_2^{n,e}$.

Now we analyze the transitions described in Figure \ref{fig:columntorow}. The configuration (6) is in $\bar\cA_2^{n,e}\subseteq\cC_{\pieno}^{\vuoto}(\gs)\subseteq\cX_{sa}^{\vuoto}$ and therefore $h(\bigcup_{\a,\a'}\bar\cA_2^{n,e})=1$. Thanks to \cite[Lemma 3.3.2]{BHN}, we know that $h$ is constant on the wells and thus we analyze the transitions to and from each wells. In particular, note that during the transition from $\widetilde\cA_{k}^{\a,\a'}$ and $\bar\cA_{k+1}^{\a,\a'}$ only configurations with energy strictly smaller than $\gs$ are crossed, thus they belong to the same well and therefore we can set $h$ constant on these configurations. We set $h(\bigcup_{\a,\a'}\cA_2^{\a,\a'})=c_1$, $h\equiv c_2$ on $\bigcup_{\a,\a'}\widetilde{\cA}_2^{\a,\a'}$ and on its wells, and so on until the last set $h(\bigcup_{\a,\a'}\bar{\cA}_{l_2^*}^{\a,\a'})=c_{2l_2^*-2}$ and $h(\h)=c_{2l_2^*-1}$ for any $\h\in\cC_{sa}^*(2)$. Thus we have to minimize w.r.t.\ $c_1,c_2,...,c_{2l_2^*-1}$ the following term:
\be{minimizzare}
(1-c_1)^2+(c_1-c_2)^2+...+(c_{2l_2^*-3}-c_{2l_2^*-2})^2+(c_{2l_2^*-2}-c_{2l_2^*-1})^2+kc_{2l_2^*-1}^2
\ee

\noindent
where $k=l_2^*-1$ is the number of good sites of the configuration (12) in Figure \ref{fig:columntorow}. We prove by induction over $n$ that 
\be{cn}
c_n=\frac{1+K_n c_{n+1}}{K_n+1}, \qquad 1\leq n\leq 2l_2^*-2,
\ee

\noindent
where $K_n$ satisfies the following recurrence relation
\be{sistema}
\begin{cases}
	K_n=(K_{n-1}+1)^2, \quad 2\leq n\leq 2l_2^*-2, \\
	K_1=1.
\end{cases}
\ee

\medskip
\noindent
\underline{$n=1$} We have to prove that $c_1=\frac{1+c_2}{2}$. This can be easily checked by minimizing the function $f(c_1,c_2)=(1-c_1)^2+(c_1-c_2)^2$ with respect to $c_1$. Indeed we get:
\be{}
\frac{\partial f}{\partial c_1}=4c_1-2c_2-2=0 \ \ \Leftrightarrow \ \ c_1=\frac{1+c_2}{2}
\ee

\medskip
\noindent
\underline{$2\leq n\leq 2l_2^*-2$} Assume now that (\ref{cn}) holds for $n-1$ and we prove that it holds also for $n$. We consider the function $f(c_{n-1},c_n,c_{n+1})=(c_{n-1}-c_n)^2+(c_n-c_{n+1})^2$ and we replace the expression of $c_{n-1}$ in terms of $c_n$. Thus we get
\be{}
f(c_n,c_{n+1})=\Bigg(\frac{1+K_{n-1} c_{n}}{K_{n-1}+1}-c_n\Bigg)^2+(c_n-c_{n+1})^2
\ee

\noindent
and therefore
\be{}
\frac{\partial f}{\partial c_n}=2\Bigg(\frac{K_{n-1}}{K_{n-1}+1}-1\Bigg)\Bigg(\frac{1+K_{n-1}c_n}{K_{n-1}+1}-c_n\Bigg)+2(c_n-c_{n+1})
\ee

\noindent
The equation $\frac{\partial f}{\partial c_n}=0$ gives 
\be{}
\frac{c_n-1}{(K_{n-1}+1)^2}=c_{n+1}-c_n
\ee

\noindent
that implies (\ref{cn}) by using (\ref{sistema}). Thus we get the claim. 

Now we have to find the value of the constant $c_{2l_2^*-1}$ that minimizes (\ref{minimizzare}). By considering the function $f(c_{2l_2^*-2},c_{2l_2^*-1})=(c_{2l_2^*-2}-c_{2l_2^*-1})^2+(l_2^*-1)c_{2l_2^*-1}^2$ and proceeding in a similar way as above, we deduce that
\be{c2l*}
c_{2l_2^*-1}=\frac{1}{(l_2^*-1)(K_{2l_2^*-2}+1)^2+1}
\ee

\noindent
Finally, by (\ref{cn}) and (\ref{c2l*}) we get
\be{diff}
1-c_1=\frac{1-c_2}{2} \quad \hbox{and} \quad c_{n-1}-c_n=\frac{1-c_n}{K_{n-1}+1}
\ee

\noindent
Finally, by (\ref{diff}) we deduce that the minimizer of the quantity in (\ref{minimizzare}) is given by
\be{min}
\Big(\frac{1-c_2}{2}\Big)^2+\displaystyle\sum_{n=2}^{2l_2^*-1}\Big(\frac{1-c_n}{K_{n-1}+1}\Big)^2+\frac{l_2^*-1}{((l_2^*-1)(K_{2l_2^*-2}+1)^2+1)^2}
\ee

\noindent
where the coefficients $c_2,...,c_{2l_2^*-1}$ can be explicitly derived from (\ref{cn}) e (\ref{c2l*}). Combining (\ref{caprw}) and (\ref{min}), we get
\be{lowerbound}
\ba{ll}
\Theta_{sa}&\geq \displaystyle\sum_{\hat\h\in\bar\cD_{sa}} \CAPA^{\Lambda^+}(\partial^+\Lambda,\hbox{CR}(\hat\h))+4\Bigg[\Big(\frac{1-c_2}{2}\Big)^2+\displaystyle\sum_{n=2}^{2l_2^*-1}\Big(\frac{1-c_n}{K_{n-1}+1}\Big)^2\\
&\quad+\dfrac{l_2^*-1}{[(l_2^*-1)(K_{2l_2^*-2}+1)^2+1]^2}\Bigg]:=\Theta^{sa}_1
\ea
\ee

\noindent
The first term in the r.h.s.\ of (\ref{lowerbound}) can be treated in a similar way as \cite[Lemma 3.4.1]{BHN} for $\Lambda\ra\Z^2$. Since the remaining part of the r.h.s. of (\ref{lowerbound}) does not depend on the size of the box, that implies that we can neglect its contribute as $\Lambda\ra\Z^2$, we deduce that 
\be{theta1}
\Theta^{sa}_1\ra{4\pi N_{sa}}\frac{|\Lambda|}{\log|\Lambda|} \qquad \hbox{as} \qquad \Lambda\ra\Z^2,
\ee

\noindent
where $N_{sa}$ is computed in Proposition \ref{Nstrong}.

\subsubsection{Upper bound}

We define
\be{C++}
\cC_{sa}^{++}:=\{\h=(\hat\h,x): \ \hat\h\in\bar\cD_{sa}, \ x\in\L\setminus\hbox{CR}^{++}(\hat\h)\}.
\ee

\noindent
and we consider the following test function

\be{testfunction}
h(\h):=
\begin{cases}
	1 &\hbox{ if } \h\in\cX_{sa}^{\vuoto}, \\
	c_i &\hbox{ if } \h\in\cX_{sa}(i), i=1,...,\bar{I} \\
	g(x) &\hbox{ if } \h\in\cC_{sa}^{++}, \\
	0 &\hbox{ if } \h\in\cX_{sa}^{\pieno},
\end{cases}
\ee

\noindent
where $g(x):=h(\hat\h,x)=h(\h)$ for $\hat\h\in\bar\cD_{sa}$ and $x\in\L\setminus\hbox{CR}^{++}(\hat\h)$, i.e., $\h\in\cC_{sa}^{++}$. Thus by \cite[eq.\ (10.1)]{BN2} we get
\be{CAP}
\ba{ll}
\CAPA(\vuoto,\pieno)&\leq(1+o(1))\Bigg(\displaystyle\sum_{\hat\h\in\bar\cD_{sa}}\CAPA^{\L^+}(\partial^+\L,\hbox{CR}^{++}(\hat\h))+\min_{c_1,...,c_I}\\
&\quad\displaystyle\min_{h:\cX^*_{sa}\ra [0,1] \atop h_{|\cX_{sa}^{\vuoto}}=1, h_{|\cX_{sa}^{\pieno}}=0, h_{|\cX_{sa}(i)}=c_i, i=1,...,\bar{I}} \frac{1}{2} \sum_{\h,\h'\in\cX^*}\mu_{\b}(\h,\h') c_{\b}(\h,\h')[h(\h)-h(\h')]^2
\ea
\ee

\noindent
where $\CAPA^{\Lambda^+}(\partial^+\Lambda,F)$ is defined in (\ref{cap+}).
We have to analyze the possible transitions between $\cX_{sa}^{\pieno}$ and $\cX_{sa}(i)\cup\bigcup_j\{\x_{sa,j}\}$, between $\cX_{sa}(i)$ and $\cX_{sa}(j)$ with $i\neq j$ and between $\cX_{sa}(i)\cup\bigcup_j\{\x_{sa,j}\}$ and $\cX_{sa}^{\vuoto}$, where the saddles $\x_{sa,j}$ are neither saddles of the first type nor saddles of the second type (recall \cite[Definitions 3.2, 3.4]{BN2}). 

We set $\cX_{sa}(1)=\bigcup_{\a\in\{n,s\}}\bigcup_{\a'\in\{w,e\}}\cA_{2}^{\a,\a'}$ (an example is given in Figure \ref{fig:columntorow} by the configuration (7)). We set $\cX_{sa}(2)$ as the union of $\widetilde\cA_{2}^{\a,\a'}$, $\bar\cA_{3}^{\a,\a'}$ and the configurations connecting these sets in the path described in Figure \ref{fig:columntorow} that are in the same well. We iterate this construction until the last set $\cX_{sa}(2l_2^*-1)$ as we did for the lower bound. Furthermore, we set 
\be{}
\cX_{sa}(2l_2^*)=\cC_{sa}^B, \qquad \cX_{sa}(2l_2^*+1)=\bigcup_{\a'\in\{e,w\}}\cA_{0}^{\a'}, \qquad \cX_{sa}(2l_2^*+2)=\bigcup_{\a'\in\{e,w\}}\cA_{1}^{\a}
\ee

\noindent
Now we analyze all the transitions that give a non-trivial contribute to (\ref{CAP}).

\medskip
\noindent
{\bf Transitions between $\cX_{sa}^{\vuoto}$ and $\cX_{sa}(i)\cup\bigcup_j\{\x_{sa,j}\}$.} The transition via one step of the dynamics from $\h\in\cX_{sa}^{\vuoto}$ and $\h'\in\cX_{sa}(i)\cup\bigcup_j\{\x_{sa,j}\}$ is possible only if either $\h\in\bar\cD_{sa}$ and $\h'\in\cC_{sa}^*$ or $\h\in\bar\cA_2^{\a,\a'}$ and $\h'\in\cA_{2}^{\a,\a'}$ for some $\a\in\{n,s\}$ and $\a'\in\{w,e\}$. The latter transition contributes 4 times to the quantity in (\ref{CAP}) depending on which frame-angle is involved in the transition.

\medskip
\noindent
{\bf Transitions between $\cX_{sa}(i)$ and $\cX_{sa}(j)$.} We consider the sequence of transitions that forms the path described in Figure \ref{fig:columntorow}, the transitions between $\cC_{sa}^*(2)$ and $\cC_{sa}^B$, between $\cC_{sa}^B$ and $\cX_{sa}(2l_2^*+1)$ and between $\cC_{sa}^B$ and $\cX_{sa}(2l_2^*+2)$.

\medskip
\noindent
{\bf Transitions between $\cX_{sa}(i)$ and $\cX_{sa}^{\pieno}$.} The transition via one step of the dynamics from $\h\in\cX_{sa}(i)$ and $\h'\in\cX_{sa}^{\pieno}$ is possible only if $\h\in\cC_{sa}^*(2)$ and $\h'\in\cC_{sa}^G$.

Collecting all these transitions, by (\ref{testfunction}) and (\ref{CAP}) we get
\be{}
\ba{ll}
\CAPA(\vuoto,\pieno)&\leq\dfrac{e^{-\b\gs}}{Z_{\b}}(1+o(1))\displaystyle\min_{\bar{c},c_1,...,c_{2l_2^*+2}}4\Big[(1-c_1)^2+(c_1-c_2)^2+...+(c_{2l_2^*-2}-c_{2l_2^*-1})^2\\
&\quad+(l_2^*-1)c_{2l_2^*-1}^2\Big]+\displaystyle\sum_{\h\in\cC_{sa}^*(2) \atop \h'\in\cC_{sa}^G}\bar{c}^2+\displaystyle\sum_{\h\in\cC_{sa}^*(2)\atop \h'\in\cX_{sa}(2l_2^*)}(\bar{c}^2-c_{2l_2^*})^2\\
&\quad+2\displaystyle\sum_{\h\in\cX_{sa}(2l_2^*)\atop \h'\in\cX_{sa}(2l_2^*+1)}(c_{2l_2^*}-c_{2l_2^*+1})^2+2\displaystyle\sum_{\h\in\cX_{sa}(2l_2^*)\atop \h'\in\cX_{sa}(2l_2^*+2)}(c_{2l_2^*}-c_{2l_2^*+2})^2\\
&\quad+\displaystyle\sum_{\hat\h\in\bar\cD_{sa}}\CAPA^{\L^+}(\partial^+\L,\hbox{CR}^{++}(\hat\h))
\ea
\ee

\noindent
Rearragging the term, we get
\be{}
\ba{ll}
\CAPA(\vuoto,\pieno)&\leq\dfrac{e^{-\b\gs}}{Z_{\b}}(1+o(1))\Bigg(\displaystyle\min_{c_1,...,c_{2l_2^*-1}}4\Big[(1-c_1)^2+(c_1-c_2)^2\\
&\quad+...+(c_{2l_2^*-2}-c_{2l_2^*-1})^2+(l_2^*-1)c_{2l_2^*-1}^2\Big]+\displaystyle\min_{\bar{c},c_{2l_2^*,...,c_{2l_2^*+2}}}\displaystyle\sum_{\h\in\cC_{sa}^*(2) \atop \h'\in\cC_{sa}^G}\bar{c}^2\\
&\quad\displaystyle+\displaystyle\sum_{\h\in\cC_{sa}^*(2)\atop \h'\in\cX_{sa}(2l_2^*)}(\bar{c}-c_{2l_2^*})^2+2\displaystyle\sum_{\h\in\cX_{sa}(2l_2^*)\atop \h'\in\cX_{sa}(2l_2^*+1)}(c_{2l_2^*}-c_{2l_2^*+1})^2\\
&\quad+2\displaystyle\sum_{\h\in\cX_{sa}(2l_2^*)\atop \h'\in\cX_{sa}(2l_2^*+2)}(c_{2l_2^*}-c_{2l_2^*+2})^2+\sum_{\hat\h\in\bar\cD_{sa}}\CAPA^{\L^+}(\partial^+\L,\hbox{CR}^{++}(\hat\h))\Bigg)
\ea
\ee

\noindent
Now we analyze separately the two following terms:
\be{thetabar}
\bar\Theta_{sa}=\displaystyle\min_{c_1,...,c_{2l_2^*-1}}4\Big[(1-c_1)^2+(c_1-c_2)^2+...+(c_{2l_2^*-2}-c_{2l_2^*-1})^2+(l_2^*-1)c_{2l_2^*-1}^2\Big]
\ee

\noindent
and
\be{thetatilde}
\ba{ll}
\widetilde\Theta_{sa}&=\displaystyle\min_{\bar{c},c_{2l_2^*},...,c_{2l_2^*+2}}\displaystyle\sum_{\h\in\cC_{sa}^*(2) \atop \h'\in\cC_{sa}^G}\bar{c}^2+\displaystyle\sum_{\h\in\cC_{sa}^*(2)\atop \h'\in\cX_{sa}(2l_2^*)}(\bar{c}-c_{2l_2^*})^2+2\displaystyle\sum_{\h\in\cX_{sa}(2l_2^*)\atop \h'\in\cX_{sa}(2l_2^*+1)}(c_{2l_2^*}-c_{2l_2^*+1})^2\\
&\quad+2\displaystyle\sum_{\h\in\cX_{sa}(2l_2^*)\atop \h'\in\cX_{sa}(2l_2^*+2)}(c_{2l_2^*}-c_{2l_2^*+2})^2+\sum_{\hat\h\in\bar\cD_{sa}}\CAPA^{\L^+}(\partial^+\L,\hbox{CR}^{++}(\hat\h))
\ea
\ee

\noindent
Note that the minimum of (\ref{thetabar}) coincides with 4 times the minimum of (\ref{minimizzare}) and therefore $\bar\Theta_{sa}$ can be computed in the same way as in the lower bound (see (\ref{cn}), (\ref{c2l*}) and (\ref{min})). It is easy to check that the minimum of (\ref{thetatilde}) w.r.t.\ $\bar{c},c_{2l_2^*}$,..., $c_{2l_2^*+2}$ is obtained for $\bar{c}=c_{2l_2^*}=c_{2l_2^*+1}=c_{2l_2^*+2}=0$. Thus the term $\widetilde\Theta_{sa}$ becomes
\be{bound1}
\widetilde\Theta_{sa}=\displaystyle\sum_{\hat\h\in\bar\cD_{sa}}\CAPA^{\L^+}(\partial^+\L,\hbox{CR}^{++}(\hat\h))
\ee

\noindent
and therefore
\be{upperbound}
\ba{ll}
\Theta_{sa}&\leq\dfrac{e^{-\b\gs}}{Z_{\b}}(1+o(1))(\bar\Theta+\widetilde\Theta)\\
&=\dfrac{e^{-\b\gs}}{Z_{\b}}(1+o(1))\Bigg[4\Bigg(\Big(\dfrac{1-c_2}{2}\Big)^2+\displaystyle\sum_{n=2}^{2l_2^*-1}\Big(\dfrac{1-c_n}{K_{n-1}+1}\Big)^2\\
&\quad+\dfrac{l_2^*-1}{[(l_2^*-1)(K_{2l_2^*-2}+1)^2+1]^2}\Bigg)+\displaystyle\sum_{\hat\h\in\bar\cD_{sa}}\CAPA^{\L^+}(\partial^+\L,\hbox{CR}^{++}(\hat\h))\Bigg]:=\Theta^{sa}_2
\ea
\ee

\noindent
where the coefficients $c_2,...,c_{2l_2^*-2},c_{2l_2^*-1}$ can be explictly derived from (\ref{cn}) and (\ref{c2l*}), where the sequence $K_n$ is defined in (\ref{sistema}). The first term in the r.h.s. of (\ref{upperbound}) does not depend on the size of the box, thus we can neglect its contribution as $\Lambda\ra\Z^2$. Since the remaining part of the r.h.s. of (\ref{upperbound}) can be treated in a similar way as \cite[Lemma 3.4.1]{BHN} for $\Lambda\ra\Z^2$, we deduce that 
\be{theta2}
\Theta^{sa}_2\ra{4\pi N_{sa}}\frac{|\Lambda|}{\log|\Lambda|} \qquad \hbox{as} \qquad \Lambda\ra\Z^2,
\ee

\noindent
where $N_{sa}$ is computed in Proposition \ref{Nstrong}.

\bp{Nstrong}
\begin{equation*}
	N_{sa}=\displaystyle\sum_{k=1}^4\binom{4}{k}\binom{l_2^*+k-2}{2k-1}.
\end{equation*}
\ep

\begin{proof*}
	We have to count the number of different shapes of the clusters in $\bar\cD_{sa}$. We do this by counting in how many ways $l_2^*-1$ particles can be removed from the four bars of a $(2l_2^*-2)\times l_2^*$ rectangle starting from the corners. We split the counting according to the number $k=1,2,3,4$ of corners from which particles are removed. The number of ways in which we can choose $k$ corners is $\binom{4}{k}$. After we have removed the particles at these corners, we need to remove $l_2^*-1-k$ more particles frome either side of each corner. The number of ways in which this can be done is
	\be{}	
	\ba{lll}
	|\{(m_1,...,m_{2k})\in\N_0^{2k}:\ m_1+...+m_{2k}=l_2^*-1-k\}| \\
	\qquad=|\{(m_1,...,m_{2k})\in\N^{2k}:\ m_1+...+m_{2k}=l_2^*-1+k\}| \\
	\qquad=\displaystyle\binom{l_2^*+k-2}{2k-1}.
	\ea
	\ee

	\noindent
	Thus we get the claim.
\end{proof*}

\subsection{Proof of Theorem \ref{gapspettrale}}

Thanks to \cite[Lemma 3.6]{BNZ}, we deduce that for our model the quantity $\widetilde\G(B)$, with $B\subsetneq\cX$, defined in \cite[eq.\ (21)]{BNZ} is such that $\widetilde\G(\cX\setminus\{\pieno\})=\G_{sa}^*$. Thus Theorem \ref{gapspettrale} follows by the following proposition.

\bp{}{\cite[Proposition 3.24]{BNZ}}
For any $\e\in(0,1)$ and any $s\in\cX^s$
\be{}
\displaystyle\lim_{\b\ra\infty}\frac{1}{\b}\log t_{\b}^{mix}(\e)=\widetilde\G(\cX\setminus\{s\})=\lim_{\b\ra\infty}-\frac{1}{\b}\log \r_{\b}
\ee

\noindent
Furthermore, there exist two constants $0<c_1\leq c_2<\infty$ independent of $\b$ such that for every $\b>0$
\be{}
c_1e^{-\b\widetilde\G(\cX\setminus\{s\})}\leq\r_{\b}\leq c_2e^{-\b\widetilde\G(\cX\setminus\{s\})}
\ee
\ep

\end{document}